\documentclass{techReport}
\usepackage{enumerate}

\usepackage{graphicx}
\graphicspath{{Figures/pdf/}{Figures/}}
\newcommand{\tfig}{0.33}

\usepackage{circuitikz}
\usetikzlibrary{arrows,shapes,backgrounds,patterns}
\usepackage{pifont}
\usepackage{tikzgraphicx}
\tikzstyle{every picture}+=[remember picture]
\tikzstyle{na} = [baseline=-.5ex]
\makeatletter
\newcommand{\gettikzxy}[3]{%
  \tikz@scan@one@point\pgfutil@firstofone#1\relax
  \edef#2{\the\pgf@x}%
  \edef#3{\the\pgf@y}%
}
\makeatother

\renewcommand\ie{\textit{i.e. }}

\usepackage{booktabs}

\newcommand{\medhrule}{\specialrule{0.05em}{0.3em}{0.3em}}
\newcommand{\bighrule}{\specialrule{0.1em}{0.3em}{0.5em}}

\usepackage{amsmath}
\makeatletter
\newcommand{\reqnomode}{\tagsleft@false\let\veqno\@@eqno}
\newcommand{\leqnomode}{\tagsleft@true\let\veqno\@@leqno}
\makeatother

\usepackage{amssymb}
\usepackage{mathrsfs}
\usepackage{mathtools}
\usepackage{xstring}

\newcommand\vphi{\varphi}
\newcommand\cste{a}
\usepackage[f]{esvect}
\newcommand\vvec{\overrightarrow}
\newcommand{\fracpartial}[2]{
    \frac{\partial {#1}}{\partial {#2}}
}
\newcommand\xt{\widetilde{x}}
\newcommand\Ft{\widetilde{F}}
\newcommand\ysol{\bar{y}}
\newcommand\xsol{\bar{x}}
\newcommand\usol{\bar{u}}
\newcommand\zsol{\bar{z}}
\newcommand\psol{\bar{p}}
\newcommand\tfsol{\bar{t}_f}
\newcommand\tsol{\bar{t}\hspace{0.1em}}
\newcommand\etasol{\bar{\eta}}
\newcommand\musol{\bar{\mu}}
\DeclareMathOperator{\sign}{sign}
\DeclareMathOperator{\ad}{ad}
\DeclareMathOperator{\vect}{Span}
\DeclareMathOperator{\rank}{rank}
\DeclareMathOperator{\lie}{Lie}
\newcommand{\hampath}{\texttt{HamPath}}
\newcommand{\bocop}{\texttt{Bocop}}
\newcommand{\hybrj}{\texttt{hybrj}}
\newcommand{\dopri}{\texttt{dopri5}}
\newcommand{\dop}{\texttt{dop853}}
\newcommand{\radau}{\texttt{radau}}
\newcommand{\minpack}{\texttt{minpack}}
\newcommand{\tapenade}{\texttt{tapenade}}
\newcommand{\expmap}[3]{\exp({#2 #3}) (#1)}
\newcommand{\Hom}{h}

\newcommand{\imax}{i_\mathrm{max}}
\newcommand{\vmax}{v_\mathrm{max}}
\newcommand{\af}{\alpha_f}

\usepackage{stmaryrd}

\newcommand{\intervalle}[4]{\mathopen{#1}#2
                                \mathclose{}\mathpunct{},#3
                                \mathclose{#4}}
\newcommand{\intervalleff}[2]{\intervalle{[}{#1}{#2}{]}}
\newcommand{\intervalleof}[2]{\intervalle{(}{#1}{#2}{]}}

\newcommand{\intervalleoo}[2]{\intervalle{(}{#1}{#2}{)}}

\newcommand{\nbSet}[1]{\mathbb{#1}}
\newcommand{\setPositive}{\text{\bf{\tiny+}}}
\newcommand{\setNegative}{\mathbb{\tiny-}}
\newcommand{\setStar}{\text{*}}
\newcommand{\R}{\xR}

\newcommand{\setDeco}[2]{
    \IfEqCase{#2}{
        {s}{\nbSet{#1}^{\setStar}}
        {n}{\nbSet{#1}^{\phantom{\setStar}}_{\setNegative}}
        {p}{\nbSet{#1}^{\phantom{\setStar}}_{\setPositive}}
        {sn}{\nbSet{#1}^{\setStar}_{\setNegative}}
        {sp}{\nbSet{#1}^{\setStar}_{\setPositive}}
    }
}

\newcommand{\Rsp}{ \ensuremath{\setDeco{R}{sp}} }

\newcommand{\Ucal}{\mathcal{U}}
\newcommand{\Acal}{\mathcal{A}}
\newcommand{\Lcal}{\mathcal{L}}
\newcommand{\Tcal}{\mathcal{T}}

\newcommand{\diff}{\xdif}

\newcommand{\abs}[1]{\lvert#1\rvert} 

\newcommand{\norme}[1]{\lVert#1\rVert}

\newcommand{\petito}[1]{o\mathopen{}\left(#1\right)}

\newcommand{\enstq}[2]{\left\{#1\mathrel{}\middle|\mathrel{}#2\right\}}

\newcommand{\prodscal}[2]{\left\langle#1,#2\right\rangle}

\usepackage{yhmath}

\begin{document}
\title{Geometric and numerical methods for a state constrained minimum time control problem of an electric vehicle}
\runningtitle{Minimum time control problem of an electric vehicle}
%

\author{Olivier Cots}
\address{Toulouse Univ., INP-ENSEEIHT-IRIT, UMR CNRS 5505, 2 rue Camichel, 31071 Toulouse, France; \texttt{olivier.cots@enseeiht.fr}}
%
\date{\today}
\begin{abstract}
    In this article, the minimum time control problem of an electric vehicle is modeled as a Mayer problem in optimal control, with affine dynamics with
    respect to the control and with state constraints.
    The candidates as minimizers are selected among a set of extremals, solutions of a Hamiltonian system given by the maximum principle.
    An analysis, with the techniques of geometric control, is used first to reduce the set of candidates and then to construct the numerical methods.
    This leads to a numerical investigation based on indirect methods using the \hampath\ software. Multiple shooting and homotopy techniques
    are used to build a synthesis with respect to the bounds of the boundary sets.
\end{abstract}
\subjclass[2011]{49K15, 49M05, 90C90, 80M50.}
\keywords{Geometric optimal control; state constraints; shooting and homotopy methods; electric car.}
\maketitle

\section*{Introduction}

    In this article, we are interested in the optimal control of an electric car, with a hybrid motor which can be operated in two discrete modes,
    $u(t) \in \{-1,1\}$, leading either to acceleration with energy consumption, or to a breaking-induced recharging of the battery.
    This vehicle can be seen as a specific type of hybrid electric vehicle (HEV), which is a vast field of research. Some recent works on controlling HEVs
    and further references can be found in \cite{JVK2013,SG2007}. Moreover, the particular model we are interested in, see equations \eqref{eq:edo_u},
    have already been studied from the optimal control point of view. In recent papers \cite{Messine2014,Messine2015}, F.~Messine \textit{et al.} solved
    the problem of the minimization of the energy consumption of an electric vehicle described by equations \eqref{eq:edo_u}
    during its displacement. In this article, we intend to analyze a complementary optimal control problem by a different approach, combining
    geometric control and numerical methods based on the maximum principle.
    Besides, as mentioned in \cite{Messine2015}, optimization of driving strategy and optimization-driven assistant systems are on the rise from these
    last few decades, especially because of recent technological advances and successfully operating autonomous cars. One can found recent works,
    using different numerical approaches, such as indirect methods for off-line optimization in \cite{Jan2013},
    or nonlinear model predictive control for real-time optimization in \cite{KB2S2013}.

    The dynamics of the electric vehicle is described with three differential states: the electric current $i$ (in ampere), 
    the position of the car $\alpha$ (in meter)
    and the angular velocity $\omega$ (in radian per second).
    The dynamical system is the following:
    \begin{equation}
        \begin{aligned}
            \displaystyle \frac{\diff      i}{\diff t}(t) & = \frac{1}{L_m} \big( -R_m\, i(t) - K_m\, \omega(t) + V_\mathrm{alim}\, u(t) \big) \\
            \displaystyle \frac{\diff \alpha}{\diff t}(t) & = \frac{r}{K_r}\, \omega(t)\\
            \displaystyle \frac{\diff \omega}{\diff t}(t) & = - \frac{K_r}{r} \, g \, K_f + \frac{K_r^2 K_m}{r^2 M}\, i(t) -
            \frac{1}{2} \, \rho S C_x \, \frac{r}{K_r M} \, \abs{\omega(t)}\, \omega(t). \\
        \end{aligned}
        \label{eq:edo_u}
    \end{equation}
    To prevent the motor from mechanical damages, the current is bounded by $\abs{i(t)} \le \imax$, $\imax>0$.
    The convexification of the control domain leads to consider the control in the whole interval $\intervalleff{-1}{1}$.
    This convexification is justified by the fact that the motor can switch from one discrete mode to another at extremely high frequency.
    We may also fix some constraints on the maximal speed (\ie linear velocity) of the car. The linear velocity (in km/h) is given by the relation
    $v(t) = \omega(t) \times 3.6\,r / K_r$ and we bound it by $\abs{v(t)} \le \vmax$, $\vmax>0$.
    For given boundary conditions in the state space, the problem of interest is the minimization of the transfer time.
    Note that we consider only positive angular velocity so we replace $\abs{\omega(t)}$ by $\omega(t)$ in the dynamics.
    The parameters, given in Table~\ref{table:parameters}, correspond to an electric solar car.
    \begin{table}[ht!]
        \centering
        \begin{tabular}{lll | lll}
            \medhrule
            Parameter           & Description                       & Unit          & Parameter           & Description                       & Unit  \\
            \bighrule
            $C_x$               & Aerodynamic coefficient           &               &
            $R_m$               & Inductor resistance               &  ohms         \\
            $\rho$              & Air density                       &  kg/m$^3$     &
            $K_m$               & Motor torque coefficient          &               \\
            $V_\mathrm{alim}$   & Battery voltage                   &  volt         &
            $r$                 & Radius of the wheels              &  m            \\
            $K_f$               & Friction coefficient              &               &
            $K_r$               & Reduction coefficient             &               \\
            &                     (of the wheels on the road)       &               &
            $R_\mathrm{bat}$    & Resistance of the battery         &  ohms         \\
            $g$                 & Gravity constant                  &  m.s$^{-2}$   &
            $S$                 & Surface in front of the vehicle   &  m$^2$        \\
            $L_m$               & Inductance of the rotor           &  henry        &
            $M$                 & Total mass                        &  kg           \\
            \medhrule
            \\
        \end{tabular}
        \caption{Description of the electric solar car parameters.}
        \label{table:parameters}
    \end{table}

    Let denote by $q \coloneqq (i,\alpha,\omega)$ the state, $M \coloneqq \R^3$ the state space, $q_0 \in M$ the initial point,
    $M_c$ a submanifold of $M$ with boundary, $M_f$ the terminal submanifold and $U \coloneqq \intervalleff{-1}{1}$ the control domain.
    The optimal control problem with control and state constraints is of the form:
    \[
        \min_{u(\cdot),\, t_f} t_f, \quad \dot{q}(t) = f(q(t),u(t)) = f_0(q(t)) + u(t) \, f_1(q(t)),
        \quad t \in \intervalleff{0}{t_f}, \quad u(t) \in U, \quad q(t) \in M_c \subset M, \quad t_f > 0,
    \]
    with boundary conditions
    $
        q(0) = q_0 \in M$ and $q(t_f) \in M_f \subset M.
    $
    The optimal solution can be found as an extremal solution of the maximum principle (with state constraints)
    and analyzed with the recent advanced techniques of geometric optimal control. The in-depth analysis of the minimum
    time problem without state constraints leads to the conclusion that the optimal
    policy is a single positive bang arc, therefore the optimal control is constant and maximum everywhere. There are no subarcs with
    intermediate values on the norm of the control, namely singular arcs.
    If we take into account the state constraints,
    then it introduces more complex structures with boundary ($q(t) \in \partial M_c$) and interior ($q(t) \in \mathring{M_c}$) arcs.
    We fully determine each type of extremals and we give new junction conditions between bang and boundary arcs.
    Besides, the local classification of the bang-bang extremals near the switching surface combined with the analysis of the state constrained problem
    provides a local time minimal synthesis which gives better insight into the structure of the optimal solutions.
    These theoretical results are first used to reduce the set of candidates as minimizers, but they are also necessary to build up the
    numerical methods. This geometric analysis leads to a numerical investigation based on indirect methods using the \hampath\ software.
    For one particular optimal control problem, we have to define the associated Multi-Point Boundary Value Problem which is solved by shooting techniques.
    However, we are interested in solving a family of optimal control problems and thus, we use differential path following (or homotopy) methods.
    We combine multiple shooting and homotopy techniques to study first a practical case ($\imax = 150$) and then to build a synthesis with respect to
    the parameters $\imax$ and $\vmax$.

    The paper is organized as follows. The optimal control problem is defined in section \ref{sec:OCPDefinition}. In section \ref{sec:unconstrained},
    we analyze the state unconstrained problem, while section \ref{sec:constrained} is devoted to the state constrained case. The numerical methods
    are presented in section \ref{sec:numericalMethods} and section \ref{sec:numericalResults} describes the numerical simulations. Section
    \ref{sec:conclusion} concludes the article.

\section{The Mathematical model and the Mayer optimal control problem}
    \label{sec:OCPDefinition}

    \subsection{Preliminaries}


        We first recall some basic facts from symplectic and differential geometries to introduce the Hamiltonian function and the adjoint vector (or covector).
        Most of the notations are taken from \cite{Agrachev2004}.

        Let $M$ be a smooth manifold of dimension $n$. Let $T^*M$ denote
        the cotangent bundle of $M$ and $\sigma$ the canonical symplectic form on $T^*M$. Recall that $(T^*M,\sigma)$ is a
        \emph{symplectic manifold}, that is, $\sigma$ is a smooth exterior 2-form on $T^*M$ which is closed and non-degenerated. For a smooth
        function $h$ on $T^*M$, we write $\vec{h}$ the \emph{Hamiltonian vector field} on $T^*M$ defined by $i_{\vec{h}} \sigma = - \diff h$,
        where $i_{\vec{h}} \sigma$ is the interior product of $\sigma$ by $\vec{h}$: $\forall\, z\in T^*M$, $\forall\, w \in T_z(T^*M)$,
        $i_{\vec{h}} \sigma_z \cdot w  \coloneqq  \sigma_z (\vec{h}(z),w)$. The function $h$ is called the \emph{Hamiltonian} function.

        Let $(q,p)$ denote the Darboux coordinates on $T^*M$. In these canonical coordinates, the Liouville 1-form $s\in \Lambda^1(T^*M)$
        writes $s = p \diff q$, and the canonical symplectic form is $\sigma = \diff s = \diff p \wedge \diff q$. If $\vphi \colon M \to N$
        is a smooth mapping between smooth manifolds, then its differential $\vphi_* \colon T_q M \to T_{x} N$, $q\in M$, $x=\vphi(q)$,
        has an adjoint mapping $\vphi^*  \coloneqq  (\vphi_*)^* \colon T^*_x N \to T^*_{q} M$ defined as follows:
        $\forall\, p_x \in T_{x}^* N$, $\vphi^* p_x = p_x \circ \vphi_*$, and $\forall\, v \in T_q M$,
        $ \prodscal{\vphi^* p_x}{v} = \prodscal{p_x}{\vphi_* v}$.
        If $\vphi$ is a diffeomorphism on $M$, then it induces the lifted diffeomorphism $\Phi$ on $T^*M$,
        $\Phi(p,q)  \coloneqq  (\vphi(q), (\vphi^*)^{-1}p)$, which is a Mathieu symplectic transformation.

    \subsubsection*{Notation}


        Here, we adopt the following notation used throughout the paper. Let $F_0$, $F_1$ be two smooth vector fields on $M$, $c$ a smooth function on M.
        We use $\ad F_0$ to denote the operator acting on vector fields $F_1 \mapsto [F_0,F_1]  \coloneqq  F_0 \cdot F_1 - F_1 \cdot F_0$,  with
        $(F_0 \cdot F_1)(x) = \diff F_1(x) \, F_0(x)$, which gives the \emph{Lie bracket}.
        The \emph{Lie derivative} $\Lcal_{F_0}c$ of $c$ along $F_0$
        is simply written $F_0\cdot c$. Denoting $H_0$, $H_1$ the Hamiltonian lifts of $F_0$, $F_1$, then the \emph{Poisson bracket} of $H_0$ and $H_1$ is
        $\{H_0,H_1\}  \coloneqq  \vvec{H_0}\cdot H_1$. We also use the notation $H_{01}$ (resp. $F_{01}$) to write the bracket $\{H_0,H_1\}$ (resp. $[F_0,F_1]$)
        and so forth. Besides, since $H_0$, $H_1$ are Hamiltonian lifts, we have $\{H_0,H_1\}= \prodscal{ p }{ [F_0,F_1] }$.


    \subsection{Electric car model}
        \label{sec:Model}

        The dynamics of the electric car is given by the smooth control vector field: $(q,u) \in M \times U \to f(q,u) \in T_q M$.
        We make a simple normalization introducing the diffeomorphism on $M$,
        \[
            x  \coloneqq  \vphi(q) = \left( \frac{q_1}{\imax}, \frac{q_2}{\af}, \frac{q_3}{\omega_\mathrm{max}} \right), \quad q = (q_1, q_2, q_3),
        \]
        where $\imax>0$ is the maximal current, $\af>0$ is the wished final position (\ie the position to reach at final time) and
        $w_\mathrm{max}>0$ is the maximal angular velocity. We define also $\vmax  \coloneqq  w_\mathrm{max} \times 3.6\,r / K_r$ the maximal linear velocity.
        The dynamics becomes:
        \[
            \dot{x}(t) = \diff \vphi({\vphi^{-1}(x(t))}) \cdot f(\vphi^{-1}(x(t)),u(t))  \coloneqq  F(x(t),u(t)),
        \]
        and any covector $p \in T^*_q M$ is transformed to the covector $p_x$, with the relation $p = {}^t \diff \vphi({\vphi^{-1}(x)}) \, p_x$.
        \begin{rmrk}
            From now on, we denote by $p$ the covector $p_x$ to simplify notations.
        \end{rmrk}

        We introduce the vector of parameters
        \[
            w  \coloneqq  \left(\frac{1}{L_m}, R_m, K_m, V_\mathrm{alim}, \frac{r}{K_r}, g\, K_f, \frac{1}{M},
            \frac{1}{2}\rho\, S\, C_x, R_\mathrm{bat}, \imax, \af, \omega_\mathrm{max}\right),
        \]
        and we have $\forall i \in \llbracket 1, 12 \rrbracket$, $w_i > 0$.
        In the $x$-coordinates, the dynamics writes
        \begin{equation}
            \dot{x}(t) = F_0(x(t)) + u(t) \, F_1(x(t)), \quad x(t) \in \R^3,
            \label{eq:system}
        \end{equation}
        where the smooth vector fields are given by
        \begin{equation*}
            \begin{aligned}
                F_0(x)  =  (\cste_1 x_1 + \cste_2 x_3)             \fracpartial{}{x_1} +
                            \cste_3 x_3                            \fracpartial{}{x_2} +
                            (\cste_4 + \cste_5 x_1 + \cste_6 x_3^2) \fracpartial{}{x_3}, \quad
                F_1(x)  = \cste_7 \fracpartial{}{x_1},
            \end{aligned}
        \end{equation*}
        with
        \[
        \begin{array}[ht!]{llll}
            \displaystyle \cste_1 = -w_1 w_2,                           &
            \displaystyle \cste_2 = -\frac{w_1 w_3 w_{12}}{w_{10}},     &
            \displaystyle \cste_3 = \frac{w_5 w_{12}}{w_{11}},          & \\[1em]
            \displaystyle \cste_4 = -\frac{w_6}{w_5 w_{12}},            &
            \displaystyle \cste_5 = \frac{w_3 w_7 w_{10}}{w_5^2 w_{12}},&
            \displaystyle \cste_6 = -w_5 w_7 w_8 w_{12},                &
            \displaystyle \cste_7 = \frac{w_1 w_4}{w_{10}}.               \\[1em]
        \end{array}
        \]
        \begin{dfntn}
            \emph{The minimum time control problem of the electric vehicle} is the following optimal
            control problem \eqref{ref:OCP_tf_min}: starting from $x_0=(0,0,0)$, reach in
            the minimum time $t_f$ the fixed normalized position $x_2(t_f)=1$ (corresponding to
            $q_2(t_f) = \af$),
            while satisfying the following control and path constraints:
            \begin{equation*}
                \begin{aligned}
                    u(t) \in U   & \coloneqq  \intervalleff{-1}{1},                                                          \quad t \in \intervalleff{0}{t_f}~ a.e., \\
                    x(t) \in M_c & \coloneqq  \intervalleff{-1}{1} \times \R \times \intervalleff{-1}{1} \subset M  \coloneqq  \R^3,  \quad t \in \intervalleff{0}{t_f}. \\
                \end{aligned}
            \end{equation*}
        \end{dfntn}
        In this article we are interested in solving \eqref{ref:OCP_tf_min} for different values of $w_{10} = \imax$
        and $w_{12} = \omega_\mathrm{max}$.
        The problem \eqref{ref:OCP_tf_min} can be stated as a \emph{Mayer problem} summerized this way:
%
        \begin{equation}
            \tag{$P_\mathrm{tmin}$}
            \left\{
            \begin{array}{l l}
                \displaystyle g(t_f,x(t_f)) &  \coloneqq  t_f \longrightarrow \displaystyle \min_{u(\cdot),\, t_f},                                          \\[1.0em]
                \dot{x}(t)              & = F_0(x(t)) + u(t) \, F_1(x(t)), \quad u(t) \in U,  \quad t\in \intervalleff{0}{t_f}~a.e.,\quad t_f > 0,  \\[0.5em]
                x(0)                    & = x_0,                                                                                                    \\[0.5em]
                x(t)                    & \in M_c, \quad t \in \intervalleff{0}{t_f},                                                               \\[0.5em]
                b(x(t_f))               & \coloneqq x_2(t_f) - 1 = 0.
            \end{array}
            \right.
            \label{ref:OCP_tf_min}
        \end{equation}

        \begin{rmrk}
            In \cite{Messine2014,Messine2015}, the authors investigate the problem of the minimization of the energy consumption,
            the transfer time $t_f$ being fixed.
            It is then reasonable to study the minimum time problem. The energy minimum problem can be stated as a Mayer problem of the following form:
            \begin{equation}
                \tag{$P_\mathrm{emin}$}
                \left\{
                \begin{array}{l l}
                    \displaystyle g(t_f,\xt(t_f)) &  \coloneqq  x_4(t_f) \longrightarrow \displaystyle \min_{u(\cdot)},
                    \quad \xt  \coloneqq  (x,x_4) \in \R^4 ,\quad t_f \text{ fixed},                                                                         \\[1.0em]
                    \dot{\xt}(t)            & = \Ft_0(\xt(t)) + u(t) \, \Ft_1(\xt(t)), \quad u(t) \in U,  \quad t\in \intervalleff{0}{t_f}~a.e.,  \\[0.5em]
                    x(0)                    & = x_0, \quad                                                                                                
                    x_4(0)                   = 0,                                                                                                  \\[0.5em]
                    x(t)                    & \in M_c, \quad t \in \intervalleff{0}{t_f},                                                           \\[0.5em]
                    b(x(t_f))               & = 0,
                \end{array}
                \right.
                \label{ref:OCP_E_min}
            \end{equation}
            where $\displaystyle \Ft_0(\xt)  \coloneqq  F_0(x) + w_9 w_{10}^2 x_1^2 \fracpartial{}{x_4}$ and
            $\displaystyle \Ft_1(\xt)  \coloneqq  F_1(x) + w_4 w_{10} x_1 \fracpartial{}{x_4}$.
        \end{rmrk}

\section{Geometric analysis of problem \eqref{ref:OCP_tf_min} without state constraints}
\label{sec:unconstrained}

    \subsection{Necessary optimality conditions}
        \label{sec:PMP}

        Let consider the problem \eqref{ref:OCP_tf_min} and assume the state constraints are relaxed: $M_c = M$.
        This optimal control problem can be written as the infinite dimensional minimization problem:
        \begin{equation*}
            \min \enstq{  t_f }{ t_f > 0,\, u(\cdot) \in \Ucal,\, b \circ E_{x_0}(t_f,u(\cdot)) = 0 },
        \end{equation*}
        where $x_0 \in M$ is fixed, $\Ucal \coloneqq \xLn{\infty}(\intervalleff{0}{t_f},U)$ is the set of \emph{admissible controls}\footnote{
        The set of admissible controls is the set of $\xLn{\infty}$-mappings on $\intervalleff{0}{t_f}$ taking their values in $U$
        such that the associated trajectory $x(\cdot)$ is globally defined on $\intervalleff{0}{t_f}$.},
        and where $E_{x_0}$ is the \emph{end-point mapping} defined by: $E_{x_0} \colon \R^+ \times \Ucal \to M$,
        $E_{x_0}(t,u(\cdot))  \coloneqq  x(t,x_0,u(\cdot))$, where $t\mapsto x(t,x_0,u(\cdot))$ is the trajectory solution of \eqref{eq:system},
        corresponding to the control $u(\cdot)$ such that $x(0,x_0,u(\cdot))=x_0$.

        This leads to the following necessary optimality conditions according to the Pontryagin Maximum Principle (PMP), see
        \cite{Pontryagin1986}.
        Define the pseudo-Hamiltonian:
        \begin{equation*}
            \begin{array}{rcl}
                H \colon T^*M \times U   & \longrightarrow   & \R \\
                (x,p,u)     & \longmapsto       & H(x,p,u)  \coloneqq  \prodscal{ p}{ F_0(x) + u\, F_1(x) }.
            \end{array}
        \end{equation*}
        By virtue of the \emph{maximum principle},
        if $(\usol(\cdot),\tfsol)$ is optimal then the associated trajectory $\xsol(\cdot)$ is the projection of an
        absolutely continuous integral curve $(\xsol(\cdot),\psol(\cdot)) \colon \intervalleff{0}{\tfsol} \to T^*M$ of
        $\vvec{H}  \coloneqq  (\frac{\partial H}{\partial p},-\frac{\partial H}{\partial x})$
        such that the following \emph{maximization condition} holds for almost every $t \in \intervalleff{0}{\tfsol}$:
        \begin{equation}
                \displaystyle H(\xsol(t),\psol(t),\usol(t)) = \max_{u\in U}H(\xsol(t),\psol(t),u).
                \label{eq:max}
        \end{equation}
        Note that the function $h(\xsol(t),\psol(t)) \coloneqq H(\xsol(t),\psol(t),\usol(t))$ is constant.
        The boundary conditions must be satisfied and we have the following \emph{transversality conditions}:
        \begin{equation*}
            \psol(\tfsol) = \mu \, b'(\xsol(\tfsol)) = (0, \mu, 0),
            \quad \mu \in \R.
        \end{equation*}
        Since $t_f$ is free, if $\usol(\cdot)$ is continuous at time $\tfsol$, then
        \begin{equation}
            H(\xsol(\tfsol),\psol(\tfsol),\usol(\tfsol)) = - p^0 \, \frac{\partial g}{\partial t}(\tfsol,\xsol(\tfsol)) = -p^0,
            \quad p^0 \in \R.
            \label{eq:H}
        \end{equation}
        Moreover, the constant $p^0$ is nonpositive and $(\psol(\cdot),p^0) \ne (0,0)$. Either $p^0 = 0$ (\emph{abnormal case}), or
        $p^0$ can be set to $-1$ by homogeneity (\emph{normal case}).
        Since $x_2$ is cyclic, $p_2$ defines a first integral and the translation $x_2 \to x_2 + c$ defines a one-parameter group of symmetries. 
        This is due to Noether theorem in Hamiltonian form.
        Hence, by symmetry, we can fix $x_2(0) = 0$. The adjoint equation is
        \[
            \dot{p}(t) = - p(t)\, F'_0(x(t)) =  - \left( \cste_1 p_1(t) + \cste_5 p_3(t) \right) \fracpartial{}{p_1}
                                                - \left( \cste_2 p_1(t) + \cste_3 p_2 + 2\cste_6 x_3(t) p_3(t) \right) \fracpartial{}{p_3}.
        \]
        Note that $h$ defines a second first integral and that $\psol_2 = \mu \ne 0$, otherwise we would have $\psol(\tfsol)=0$ which would 
        imply $h(\xsol(t),\psol(t)) = -p^0 = 0$, and $(\psol(\tfsol),p^0) = (0,0)$ is not possible.

        \begin{prpstn}
            $p^0 = 0$ (abnormal case) if and only if $\xsol_3(\tfsol) = 0$.
        \end{prpstn}
        \begin{proof}
            We have $H(\xsol(\tfsol),\psol(\tfsol),\usol(\tfsol)) = \psol_2\cste_3 \xsol_3(\tfsol) = -p^0$, with $\psol_2 \ne 0$.
            The result follows.
        \end{proof}

        \begin{dfntn}
            A triple $(x(\cdot),p(\cdot),u(\cdot))$ where $(x(\cdot),p(\cdot))$ is an integral curve of $\vec{H}$ and satisfying \eqref{eq:max}
            is called an \emph{extremal}. Any extremal satisfying the boundary conditions,
            the transversality conditions and condition \eqref{eq:H} is called a BC-extremal.
            Let $(\zsol(\cdot),\usol(\cdot))$ be an extremal, with $\zsol(\cdot) \coloneqq (\xsol(\cdot),\psol(\cdot))$.
            If $h$ is defined and smooth in a neighborhood of $\zsol$ then
            $h$ defines a \emph{true Hamiltonian}, and $\zsol$ is also an integral curve of $\vec{h}$.
            We define the \emph{Hamiltonian lifts} $H_0(x,p)  \coloneqq  \prodscal{ p }{ F_0(x) }$, $H_1(x,p)  \coloneqq  \prodscal{ p }{ F_1(x) }$ 
            and the \emph{switching function}
            $\Phi(t)  \coloneqq  H_1(x(t),p(t))$.
        \end{dfntn}

        \subsubsection*{\textbf{Regular extremals}}
        It follows from \eqref{eq:max} that if $\Phi(t) \ne 0$, $u(t) = \sign(\Phi(t))$. We say that a trajectory $x(\cdot)$
        restricted to a subinterval $I \subset \intervalleff{0}{t_f}$, not reduced to a singleton, is a \emph{bang arc}
        if $u(\cdot)$ is constant on $I$, taking values in $\{-1,1\}$.
        The trajectory is called \emph{bang-bang} if it is the concatenation of a finite number of bang arcs.

        \subsubsection*{\textbf{Singular extremals}}
        We say that a trajectory $x(\cdot)$ restricted to a subinterval $I \subset \intervalleff{0}{t_f}$, not reduced to a singleton,
        is a \emph{singular arc} if the associated
        extremal lift satisfies $\Phi(t) = 0$, $\forall t \in I$. A trajectory $x(\cdot)$ defined on $\intervalleff{0}{t_f}$, with $x(0)=x_0$,
        is singular if the associated control $u(\cdot)$ is a critical point of $u(\cdot) \mapsto E_{x_0}(u(\cdot),t_f)$. In this case,
        $u(\cdot)$ is said to be singular. From the geometric point of view, a control $\usol(\cdot) \in \Ucal$ such that
        $E_{x_0}(\usol(\cdot),t_f) \in \partial \Acal(x_0,t_f)$, where $\Acal(x_0,t_f) \coloneqq  \enstq{  E_{x_0}(u(\cdot),t_f) }{ u(\cdot)\in\Ucal }$
        is the accessibility set at time $t_f$, is singular.


    \subsection{Lie bracket configuration}
        \label{sec:LieBracketConf}

        The local behavior of the affine control system \eqref{eq:system} is determined by the values of the drift $F_0$, the control vector field $F_1$, and all
        their Lie brackets at a reference point $x \in M$.
        Loosely speaking, the values and dependencies of these vector fields at $x$ is called the \emph{Lie bracket configuration} of the system at $x$.
        According to \cite[Chapter 7]{Schattler2012}, the concept of \emph{codimension} is crucial to organize the Lie bracket conditions into groups of 
        increasing degrees of degeneracy.
        The codimension is given by the number of linearly independent ``relevant'' equality relations that hold between these vector fields at $x$.
        As both the dimension of the state space and the codimension of the Lie bracket configuration at the point $x$ increase, the local optimal
        synthesis becomes increasingly more complex.
        According to the following lemma, the 3-dimensional system \eqref{eq:system} is a part of the codimension-0 case\footnote{The points such that $x_3 = 0$ have
        different linearly independent equality relations than the others.}.
        \begin{lmm}
            \label{lemma:lieBrackets}
            We have:
            \medskip
            \begin{enumerate}[i)]
                \item 
                    $
                    \left\{
                    \begin{array}{ll}
                        F_{01}  &= - \cste_7 \big( \cste_1 \fracpartial{}{x_1} + \cste_5 \fracpartial{}{x_3} \big),                      \\[0.5em]
                        F_{001} &= \cste_7  \Big(
                                                        \big( \cste_1^2 + \cste_2 \cste_5 \big)              \fracpartial{}{x_1} +
                                                            \cste_3 \cste_5                                  \fracpartial{}{x_2} +
                                                        \cste_5 \big( \cste_1 + 2 \cste_6 x_3 \big)   \fracpartial{}{x_3} \Big),         \\[0.5em]
                        F_{101} & = F_{1001} = 0,                                                                                        \\[0.5em]
                        F_{10001} & = 2 \cste_5^2 \cste_6 \cste_7^2 \fracpartial{}{x_3},                                                \\[0.5em]
                    \end{array}
                    \right.
                    $
                    \medskip
                \item $ \dim \vect(F_1(x),F_{01}(x),F_{001}(x)) = 3$ and $\lie_x (\{ F_0, F_1\}) = T_x M$, $x\in M$ 
                    (\footnote{$\lie_x (\{ F_0, F_1\})$ is the Lie algebra generated by the family $\{ F_0, F_1\}$ at $x$.}).
                    \medskip
                \item $ \dim \vect(F_0(x),F_{1}(x),F_{01}(x)) = 3$, $x\in M$, $x_3 \ne 0$.
            \end{enumerate}
        \end{lmm}
        \begin{proof}
            The proof is straightforward.
            \begin{enumerate}[i)]
                \item    Computing, $F_{01} = \ad F_0 \cdot F_1 = - F_1 \cdot F_0$ and $F_{001} = \ad F_0 \cdot F_{01} = - F_{01} \cdot F_0$.
                        Computing again, $F_{101} = F_1 \cdot F_{01} - F_{01} \cdot F_{1} = 0$, since $F_1$ and $F_{01}$ does not depend on $x$.
                        Using the \emph{Jacobi identity}, $F_{1001} = [F_1,[F_0,F_{01}]] = -[F_0,[F_{01},F_1]] - [F_{01},[F_{1},F_0]] =
                        [F_0,F_{101}] + [F_{01},F_{01}] = 0$.
                        Computing, $F_{10001} = 2 \cste_5^2 \cste_6 \cste_7^2 \fracpartial{}{x_3}$.
                \item  Let $x\in M$, then $\det (F_1(x),F_{01}(x),F_{001}(x)) = \cste_3 \cste_5^2 \cste_7^3 \ne 0$, so $\lie_x (\{ F_0, F_1\}) = T_x M$.
                \item   Let $x\in M$, then $\det (F_0(x),F_{1}(x),F_{01}(x)) = \cste_3 \cste_5 \cste_7^2 x_3 \ne 0$, when $x_3\ne 0$.
            \end{enumerate}
        \end{proof}



    \subsection{Regular extremals}
    \label{sec:SingularRegularExtremals}

        \begin{lmm}
            The switching function $\Phi$ is $\xCn{3}$.
        \end{lmm}
        \begin{proof}
            $\Phi = H_1 = \cste_7 p_1$ is absolutely continuous, so differentiable almost everywhere. $\dot{\Phi}=H_{01}=-\cste_7(\cste_1 p_1 + \cste_5 p_3)$
            is absolutely continuous, therefore continuous, and so $\Phi$ is $\xCn{1}$.
            Likewise $\dot{\Phi}=H_{01}$, $\ddot{\Phi}=H_{001} + u\, H_{101} = H_{001}$ (since $H_{101}=0$ by lemma \ref{lemma:lieBrackets}) and
            $\dddot{\Phi}= H_{0001} + u\, H_{1001}=H_{0001}$ (since $H_{1001}=0$) are absolutely continuous and we can conclude that
            $\Phi$ is $\xCn{3}$.
        \end{proof}

        \begin{prpstn}
            \label{prop:bang-bang}
            Every extremal is bang-bang.
        \end{prpstn}
        \begin{proof}
            Let $(z(\cdot),u(\cdot))$ be an extremal defined on $\intervalleff{0}{t_f}$.
            Let assume there exists a singular part $z_s(\cdot)=(x_s(\cdot),p_s(\cdot))$ defined on $I\subset \intervalleff{0}{t_f}$.
            Then for all $t\in I$, $0 = \Phi(t) = H_1(z_s (t))$. Differentiating twice with respect to time,
            we have $0 =  H_{1}(z_s(t)) = H_{01}(z_s(t)) = H_{001}(z_s(t))$
            which implies $p_s(\cdot) = 0$ since $F_1(x_s(t))$, $F_{01}(x_s(t))$ and $F_{001}(x_s(t))$ are independent, which is not possible.
        \end{proof}

        \begin{prpstn}
            The switching function has finitely many zeros.
        \end{prpstn}
        \begin{proof}
            Following the proof of \cite[proposition 2.1]{GeHa2006}, if not, then along $(z(\cdot),u(\cdot))$ defined on $\intervalleff{0}{t_f}$,
            there exists a sequence $(t_k)_k$, $t_k \in \intervalleff{0}{t_f}$, all disctinct, such that $\Phi(t_k) = 0$. But $t_f$ is fixed, so there exists a
            subsequence, always noted $(t_k)_k$, which converges to $\tsol$. As $\Phi$ is $\xCn{0}$, $\Phi(\tsol) = 0$. But $\Phi$ is also $\xCn{1}$, hence
            \[
                \frac{\Phi(t_k) - \Phi(\tsol)}{t_k - \tsol} = 0 \longrightarrow \dot{\Phi}(\tsol) = 0.
            \]
            Besides, by Rolle's theorem, for every $k$, there exists at least one $\tau_k \in \intervalleoo{t_k}{t_{k+1}}$ such that $\dot{\Phi}(\tau_k) = 0$,
            and by squeeze theorem the sequence $(\tau_k)_k$ converges to $\tsol$. But $\Phi$ is also $\xCn{2}$ and we have
            \[
                \frac{\dot{\Phi}(\tau_k) - \dot{\Phi}(\tsol)}{\tau_k - \tsol} = 0 \longrightarrow \ddot{\Phi}(\tsol) = 0.
            \]
            In conclusion, at time $\tsol \in \intervalleff{0}{t_f}$, $\Phi(\tsol) = \dot{\Phi}(\tsol) = \ddot{\Phi}(\tsol) = 0$, whence $p(\tsol) = 0$
            which is impossible.
        \end{proof}

        Hence, any BC-extremal is a concatenation of only bang arcs. It is a difficult task to give upper bounds on the number of switchings
        for global time-optimal trajectories. However, any global time-optimal trajectory is locally time-optimal
        and next section is dedicated to the classification of such local bang-bang time-optimal trajectories.

    \subsection{Classification of bang-bang extremals}
        \label{sec:Classification}

        We define first the colinear set $C  \coloneqq  \enstq{ x \in M }{ F_0(x) \parallel F_1(x) }$.
        Then $x\in C$ if and only if $x_3 = 0$ and $\cste_4 + \cste_5 x_1 = 0$.
        The point $(x,u) = (-\frac{\cste_4}{\cste_5}, x_2, 0, \frac{\cste_1 \cste_4}{\cste_5 \cste_7})$, $x_2 \in \R$, is an equilibrium point of the
        control system \eqref{eq:system}.
        Let then define the following sets:
        \begin{equation*}
            \displaystyle   \Sigma^+_{1}  \coloneqq  \enstq{  z \in T^*M }{ H_{1} (z) > 0 }, \quad
                            \Sigma^0_{1}  \coloneqq  \enstq{  z \in T^*M }{ H_{1} (z) = 0 }, \quad
                            \Sigma^-_{1}  \coloneqq  \enstq{  z \in T^*M }{ H_{1} (z) < 0 },
        \end{equation*}
        and
        \begin{equation*}
            \displaystyle   \Sigma^+_{01} \coloneqq  \enstq{  z \in T^*M }{ H_{01}(z) > 0 }, \quad
                            \Sigma^0_{01} \coloneqq  \enstq{  z \in T^*M }{ H_{01}(z) = 0 }, \quad
                            \Sigma^-_{01} \coloneqq  \enstq{  z \in T^*M }{ H_{01}(z) < 0 }.
        \end{equation*}
        Singular extremals (if any) are entirely contained in $\Sigma_{s} \coloneqq \Sigma^0_{1}\cap \Sigma^0_{01}$.
        A crucial point is to apply the results obtained by Kupka
        \cite{kupka87} (see also \cite{BoChSing}) to classify extremal curves near the \emph{switching surface} $\Sigma^0_{1}$.
        We label $\gamma_b$ a bang extremal, $\gamma_+$, $\gamma_-$, a bang extremal such that $u(\cdot) =  +1$, $u(\cdot) = -1$ and $\gamma_s$ a singular extremal.
        We denote by $\gamma_1\gamma_2$, the concatenation of an arc $\gamma_1$ followed by an arc $\gamma_2$.
        Since for every $z\in T^*M$, $\ker \diff H_1(z)$ is transverse to $\ker \diff H_{01}(z)$, then $\Sigma^0_{1}$ and $\Sigma^0_{01}$ are both
        smooth submanifold of codimension 1 transverse at each point.

        \subsubsection{{Normal switching points}} 
\begin{figure}[ht!]
        \centering
        \begin{tikzpicture}[scale=1.4]

            \coordinate (O) at (0,0);
            \def\xshift{5em}
            \def\xshiftExt{4em}
            \def\yshiftExt{3em}
            \draw[thin] ([xshift=-\xshift] O) -- ([xshift=\xshift] O);
            \draw ([                yshift= 0em]    O) node{ {$\bullet$}};
            \draw ([                yshift=-0.5em]  O) node[below]  {$\zsol$};
            \draw ([xshift=-\xshift]                O) node[left]   {$\Sigma^0_{1}$};
            \draw ([xshift=-\xshift,yshift= 1em]    O) node[above]  {$\Sigma^+_{1}$};
            \draw ([xshift=-\xshift,yshift=-1em]    O) node[below]  {$\Sigma^-_{1}$};

            \draw[very thick] ([xshift=\xshiftExt,yshift= \yshiftExt] O) -- (O) node[midway, sloped, rotate=180] {\footnotesize \ding{228}};
            \draw[very thick] ([xshift=\xshiftExt,yshift=-\yshiftExt] O) -- (O) node[midway, sloped, rotate=  0] {\footnotesize \ding{228}};

            \def\eps{3em}
            \draw[ thin, gray] ([xshift=\xshiftExt-\eps,yshift= \yshiftExt] O) -- ([xshift=-\eps] O) node[midway, sloped, rotate=180] {\footnotesize \ding{228}};
            \draw[ thin, gray] ([xshift=\xshiftExt-\eps,yshift=-\yshiftExt] O) -- ([xshift=-\eps] O) node[midway, sloped, rotate=  0] {\footnotesize \ding{228}};

            \def\eps{-3em}
            \draw[ thin, gray] ([xshift=\xshiftExt-\eps,yshift= \yshiftExt] O) -- ([xshift=-\eps] O) node[midway, sloped, rotate=180] {\footnotesize \ding{228}};
            \draw[ thin, gray] ([xshift=\xshiftExt-\eps,yshift=-\yshiftExt] O) -- ([xshift=-\eps] O) node[midway, sloped, rotate=  0] {\footnotesize \ding{228}};

        \end{tikzpicture}
        \hspace{7em}
        \begin{tikzpicture}[scale=1.4]

            \coordinate (O) at (0,0);
            \def\xshift{5em}
            \def\xshiftExt{4em}
            \def\yshiftExt{3em}
            \draw[thin] ([xshift=-\xshift] O) -- ([xshift=\xshift] O);
            \draw ([                yshift= 0em]    O) node{ {$\bullet$}};
            \draw ([                yshift=-0.5em]  O) node[below]  {$\zsol$};
            \draw ([xshift=-\xshift]                O) node[left]   {$\Sigma^0_{1}$};
            \draw ([xshift=-\xshift,yshift= 1em]    O) node[above]  {$\Sigma^+_{1}$};
            \draw ([xshift=-\xshift,yshift=-1em]    O) node[below]  {$\Sigma^-_{1}$};

            \draw[very thick] ([xshift=\xshiftExt,yshift= \yshiftExt] O) -- (O) node[midway, sloped, rotate=  0] {\footnotesize \ding{228}};
            \draw[very thick] ([xshift=\xshiftExt,yshift=-\yshiftExt] O) -- (O) node[midway, sloped, rotate=180] {\footnotesize \ding{228}};

            \def\eps{3em}
            \draw[ thin, gray] ([xshift=\xshiftExt-\eps,yshift= \yshiftExt] O) -- ([xshift=-\eps] O) node[midway, sloped, rotate=  0] {\footnotesize \ding{228}};
            \draw[ thin, gray] ([xshift=\xshiftExt-\eps,yshift=-\yshiftExt] O) -- ([xshift=-\eps] O) node[midway, sloped, rotate=180] {\footnotesize \ding{228}};

            \def\eps{-3em}
            \draw[ thin, gray] ([xshift=\xshiftExt-\eps,yshift= \yshiftExt] O) -- ([xshift=-\eps] O) node[midway, sloped, rotate=  0] {\footnotesize \ding{228}};
            \draw[ thin, gray] ([xshift=\xshiftExt-\eps,yshift=-\yshiftExt] O) -- ([xshift=-\eps] O) node[midway, sloped, rotate=180] {\footnotesize \ding{228}};

        \end{tikzpicture}
        \caption{
            Bang-Bang extremal with contact of order 1. (Left) $\dot{\Phi}(\tsol) < 0$. (Right) $\dot{\Phi}(\tsol) > 0$.
        }
        \label{fig:contact1}
        \end{figure}
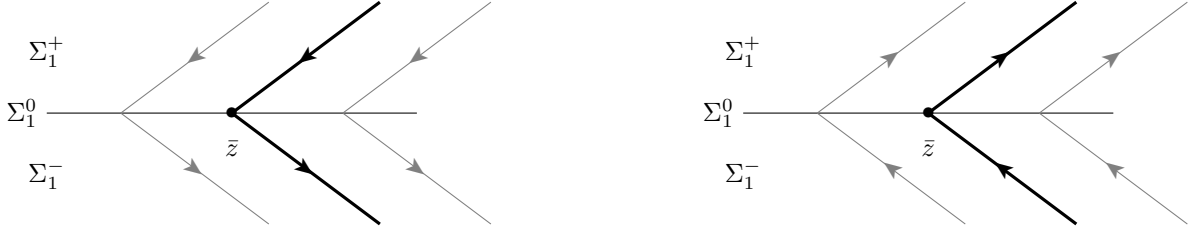

It is the case when a bang arc has a contact of order 1 with the switching surface.
        Let $\zsol \coloneqq (\xsol,\psol) \in \Sigma^0_{1} \setminus \Sigma^0_{01}$ and assume
        $\xsol \notin C$. The point $\zsol$ is called \emph{normal}. Let $(z(\cdot),u(\cdot))$ be a regular extremal passing through $\zsol$ at time $\tsol$,
        then $\Phi(\tsol) = H_1(\zsol) = 0$, $\dot{\Phi}(\tsol) = H_{01}(\zsol)\ne 0$ and $\Phi$ changes its sign at $\tsol$.
        Since $\Phi(\tsol+s) = \dot{\Phi}(\tsol) \, s + \petito{s}$ near $\tsol$, then locally
        \[
            u(\tsol+s) = \frac{\Phi(\tsol+s)}{\abs{\Phi(\tsol+s)}} = \frac{\dot{\Phi}(\tsol)}{\abs{\dot{\Phi}(\tsol)}} \cdot \frac{s}{\abs{s}} + \petito{1},
        \]
        so near $\zsol$, every extremal is
        of the form $\gamma_+\gamma_-$ if $\dot{\Phi}(\tsol) < 0$ and $\gamma_-\gamma_+$ if $\dot{\Phi}(\tsol) > 0$, see Fig.~\ref{fig:contact1}.

        \begin{prpstn}
            Let $\zsol \coloneqq (\xsol,\psol)\in T^*M$, $\xsol \not\in C$, then $\zsol$ is a normal switching point if and only if $\psol_1 = 0$ and $\psol_3 \ne 0$.
            If $\zsol$ is a normal switching point, then locally, any extremal is of the form $\gamma_+\gamma_-$ if $\psol_3>0$
            and $\gamma_-\gamma_+$ if $\psol_3 < 0$.
        \end{prpstn}

        \subsubsection{{The fold case}}
        \label{sec:foldCase}
\begin{figure}[ht!]
        \centering
        \begin{tikzpicture}[scale=1.2]

            \coordinate (O) at (0,0);
            \def\xshift{8em}
            \def\yshift{5em}
            \draw[thin] ([xshift=-\xshift] O) -- ([xshift=\xshift] O);
            \draw[thin] ([yshift=-\yshift] O) -- ([yshift=\yshift] O);

            \draw ([                yshift= 0em]    O) node{ {$\bullet$}};
            \draw ([                yshift=-0.5em]  O) node[below left]  {$\zsol$};

            \draw ([xshift=-\xshift]                O) node[left]   {$\Sigma^0_{1}$};
            \draw ([xshift=-\xshift,yshift= 1em]    O) node[above]  {$\Sigma^+_{1}$};
            \draw ([xshift=-\xshift,yshift=-1em]    O) node[below]  {$\Sigma^-_{1}$};

            \draw ([yshift=-\yshift]                O) node[below]   {$\Sigma^0_{01}$};
            \draw ([xshift=-1em,yshift=-\yshift]    O) node[left]  {$\Sigma^-_{01}$};
            \draw ([xshift= 1em,yshift=-\yshift]    O) node[right]  {$\Sigma^+_{01}$};

            \def\a{0.3}
            \draw [thin, gray] ([xshift=-1*\xshift/8,yshift=4*\yshift/5] O) .. controls +(0,-\a) and +(-\a,0) .. ([xshift=0em,yshift=2*\yshift/5] O)
            node[midway, sloped, rotate=0] {\tiny \ding{228}};
            \draw [thin, gray] ([xshift=0em,yshift=2*\yshift/5]        O) .. controls +(\a,0) and +(0,-\a) .. ([xshift=1*\xshift/8,yshift=4*\yshift/5] O);

            \def\a{0.8}
            \draw [very thick] ([xshift=-3*\xshift/8,yshift=4*\yshift/5] O) .. controls +(0,-\a) and +(-\a,0) .. ([xshift=0em,yshift=0em] O)
            node[midway, sloped, rotate=0] {\tiny \ding{228}};
            \draw [very thick] ([xshift=0em,yshift=0em]          O) .. controls +(\a,0) and +(0,-\a) .. ([xshift=3*\xshift/8,yshift=4*\yshift/5] O);

            \def\a{0.6}
            \def\b{0.3}
            \draw [thin, gray] ([xshift=-5*\xshift/8,yshift=4*\yshift/5] O) .. controls +(0,-\a) and +(-\a,\b) .. ([xshift=-2*\xshift/5,yshift=0em] O)
            node[midway, sloped, rotate=0] {\tiny \ding{228}};
            \draw [thin, gray] ([xshift=-2*\xshift/5,yshift=0em]   O) .. controls +(\a/2,-2*\b) and +(-\b,0) .. ([xshift=0em,yshift=-2.5*\yshift/5] O);
            \draw [thin, gray] ([xshift=0em,yshift=-2.5*\yshift/5] O) .. controls +(\b,0) and +(-\a/2,-2*\b)   .. ([xshift=2*\xshift/5,yshift=0em] O)
            node[near start, sloped, rotate=0] {\tiny \ding{228}};
            \draw [thin, gray] ([xshift=2*\xshift/5,yshift=0em] O) .. controls +(\a,\b) and +(0,-\a) .. ([xshift=5*\xshift/8,yshift=4*\yshift/5] O)
            node[midway, sloped, rotate=0] {\tiny \ding{228}};
            \draw ([xshift=-2*\xshift/5]  O) node[below left]  {$z_0$};
            \draw ([xshift=-2*\xshift/5]  O) node[gray] { {$\bullet$}};

            \def\a{0.6}
            \def\b{0.3}
            \draw [thin, gray] ([xshift=-7*\xshift/8,yshift=4*\yshift/5] O) .. controls +(0,-\a) and +(-\a,\b) .. ([xshift=-3.5*\xshift/5,yshift=0em] O)
            node[midway, sloped, rotate=0] {\tiny \ding{228}};

            \draw [thin, gray] ([xshift=-3.5*\xshift/5,yshift=0em]   O) .. controls +(-\a,-\b) and +(0,\a) ..
            ([xshift=-7*\xshift/8,yshift=-4*\yshift/5] O) node[midway, sloped, rotate=180] {\tiny \ding{228}};

            \def\a{0.6}
            \def\b{0.3}
            \draw [thin, gray] ([xshift=7*\xshift/8,yshift=4*\yshift/5] O) .. controls +(0,-\a) and +(\a,\b) .. ([xshift=3.5*\xshift/5,yshift=0em] O)
            node[midway, sloped, rotate=0] {\tiny \ding{228}};

            \draw [thin, gray] ([xshift=3.5*\xshift/5,yshift=0em]   O) .. controls +(\a,-\b) and +(0,\a) ..
            ([xshift=7*\xshift/8,yshift=-4*\yshift/5] O) node[midway, sloped, rotate=180] {\tiny \ding{228}};

        \end{tikzpicture}
        \caption{
            Local optimal synthesis in a neighborhood of $\zsol$, for $\ddot{\Phi}(\tsol) > 0$.
        }
        \label{fig:parabolic}
        \end{figure}
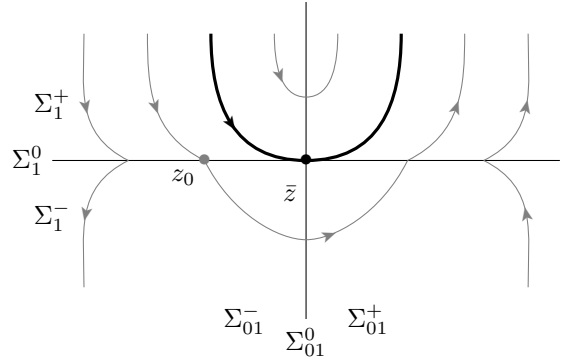

        It is the case when a bang arc has a contact of order 2 with the switching surface.
        Let $\zsol \coloneqq (\xsol,\psol) \in \Sigma_s$
        and assume $\xsol \notin C$. We define $H_\pm  \coloneqq  H_0 \pm H_1$
        and then we get $\Sigma^0_{1} = \enstq{ z \in T^*M }{ H_+(z) - H_-(z) = 0 }$ and
        $\Sigma^0_{01} = \enstq{ z \in T^*M }{ \{H_+,H_-\}(z) = 0 }$.
        \begin{lmm}
            Let $(z(\cdot),u(\cdot))$ be a regular extremal passing through $\zsol \in \Sigma_s$ at time $\tsol$, then both Hamiltonian vector fields
            $\vvec{H_+}$ and $\vvec{H_-}$ are tangent to $\Sigma^0_{1}$ at $z(\tsol) = \zsol$.
        \end{lmm}
        \begin{proof}
            $\vvec{H_+}(\zsol) \in T_{\zsol} \Sigma^0_{1}$ since
            $
                (\diff H_+ - \diff H_-) \cdot \vvec{H_+} = \{H_+, H_+-H_-\} = -\{H_+,H_-\}
            $
            which is zero at $\zsol$ since $\zsol \in \Sigma^0_{01}$. The same result goes for $\vvec{H_-}$.
        \end{proof}
        If we set $\ddot{\Phi}_\pm  \coloneqq  H_{001} \pm H_{101}$ then if both $\ddot{\Phi}_\pm \ne 0$, the contact of the trajectories of
        $\vvec{H_+}$ and $\vvec{H_-}$ with $\Sigma^0_{1}$ is of order 2. Such a point is called a \emph{fold}. According to \cite{kupka87},
        we have three cases:
        \begin{enumerate}[i)]
            \item $\ddot{\Phi}_+ \ddot{\Phi}_- > 0$: parabolic case,
            \item $\ddot{\Phi}_+ > 0$ and $\ddot{\Phi}_- < 0$: hyperbolic case,
            \item $\ddot{\Phi}_+ < 0$ and $\ddot{\Phi}_- > 0$: elliptic case.
        \end{enumerate}
        \begin{prpstn}
            The contacts with the switching surface are at most of order 2. Let $\zsol \coloneqq (\xsol,\psol)\in T^*M$, $\xsol \not\in C$, then $\zsol$ is a fold
            point if and only if $\psol_1=\psol_3=0$ and $\psol_2\ne 0$, and if $\zsol$ is a fold point, it is parabolic.
        \end{prpstn}
        \begin{proof}
            With exactly the same argument as in Proposition~\ref{prop:bang-bang}, contacts of order 3 or more are not possible.
            Let $\zsol \coloneqq (\xsol,\psol)\in T^*M$, then $\zsol$ is a fold point iff $H_1(\zsol)=H_{01}(\zsol)=0$ and $H_{001}(\zsol)\ne 0$, since
            $\ddot{\Phi} = H_{001}$, which gives the result. It is then a parabolic point since $\ddot{\Phi}_+ = \ddot{\Phi}_- = H_{001}$.
        \end{proof}
        We give the generic classification of extremals near a fold point only in the parabolic case. In this case, if there exists a singular
        extremal of minimal order passing through $\zsol$ (\ie $H_{101}(\zsol)\ne 0$), this extremal is not admissible and every extremal curve
        near $\zsol$ is bang-bang with at most two switchings, \ie is of the form $\gamma_+\gamma_-\gamma_+$ or
        $\gamma_-\gamma_+\gamma_-$ (by convention each arc of the sequence can be empty). One extremal is time minimizing while the other is
        time maximizing, depending on the sign of $\ddot{\Phi}(\tsol)$, see Fig.~\ref{fig:parabolic}.

        \subsubsection*{Application to problem \eqref{ref:OCP_tf_min}}

        Let $(z(\cdot),u(\cdot))$ be a regular extremal passing through $\zsol$ at time $\tsol$, then
        $\Phi(\tsol+s) = \frac{1}{2} \ddot{\Phi}(\tsol) \, s^2 + \petito{s^2}$ near $\tsol$, and locally,
        \[
            u(\tsol+s) = \frac{\Phi(\tsol+s)}{\abs{\Phi(\tsol+s)}} = \frac{\ddot{\Phi}(\tsol)}{\abs{\ddot{\Phi}(\tsol)}} + \petito{1}.
        \]
        Since $\zsol$ is a fold point, then $\ddot{\Phi}(\tsol) = \cste_3 \cste_5 \cste_7 p_2$, with $p_2$ constant along any extremal
        and $\cste_3 \cste_5 \cste_7>0$.
        Besides, if $(z(\cdot),u(\cdot))$ is a normal BC-extremal, then according to section \ref{sec:PMP},
        $1 = h(z(t_f)) = H_0(z(t_f)) = \cste_3 x_3(t_f) p_2$, $\cste_3>0$,
        and because we want to steer the vehicle from a zero position to a positive position
        in minimum time, $x_3(t_f)$ must be positive. As a consequence, $p_2>0$.

        Let consider now a normal switching point $z_0$, close to $\zsol$. Let $(z(\cdot),u(\cdot))$ be an extremal of the form $\gamma_+\gamma_-\gamma_+$,
        $\gamma_-$ non-empty, passing through $z_0$ at time $t_0$.
        The switching function may be approximated by
        \[
            \Phi(t_0 + s) \approx \dot{\Phi}(t_0)\, s + \frac{1}{2} \ddot{\Phi}(t_0) \, s^2,
        \]
        with $\dot{\Phi}(t_0) < 0$ and $\ddot{\Phi}(t_0) > 0$. Hence, the switching function vanishes at $t_0$ and near $t_0 + s^*$,
        where
        \[
            s^* = -2 \frac{\dot{\Phi}(t_0)}{\ddot{\Phi}(t_0)} > 0.
        \]
        Putting all together, we have the following proposition.
        \begin{prpstn}
            \label{prop:structureLocalFold}
            Let $\zsol\in T^*M$ be a fold point, then $H_{001}(\zsol)>0$ for any normal BC-extremal passing trough $\zsol$
            and the optimal policy near $\zsol$ is $\gamma_+\gamma_-\gamma_+$. Besides, the length $s$ of $\gamma_-$ is given by
            \[
                s \approx 2\, \frac{p_3(t_0)}{\cste_3 p_2 + p_3(t_0) \left( \cste_1 + 2 \cste_6 x_3(t_0) \right)},
            \]
            where $t_0$ is the first switching time, between $\gamma_+$ and $\gamma_-$.
        \end{prpstn}

    \subsection{Optimality of the $\gamma_+$ trajectory}
    \label{sec:bangOptimal}

        We conclude the section \ref{sec:unconstrained} showing that the strategy $\gamma_+$ is optimal.
        The transversality conditions imply that any BC-extremal ends at a fold point.
        By proposition \ref{prop:structureLocalFold}, the last bang arc must be a positive bang. 
        Let $\tfsol \approx 5.6156$ denote the time when the $\gamma_+$ BC-extremal reaches the target $x_2 = 1$.
        Since the final submanifold $M_f$ is of codimension 1 and by homogeneity, the final adjoint vector is uniquely determined for each $x_f \in M_f$.
        For all $x_f \in M_f$ the final adjoint vector is $p_f \coloneqq (0,p_2,0)$, with $p_2 > 0$, and we can fix $p_2 = 1$.
        Let $Z_f \coloneqq M_f \times \{p_f\} $ and write $z(\cdot,z_f)$, $z_f \in Z_f$, the solution of ${H_+}$ starting from $z_f$ and computed with
        backward integration.
        If for all $z_f \in Z_f$ and for all $t \in \intervalleof{0}{\tfsol}$, $\phi(t) = H_1(z(t,z_f)) \ne 0$, then the strategy $\gamma_+$ is optimal.

\begin{figure}[ht!]
            \def\sizeFig{0.4}
            \centering
            \includegraphics[width=\sizeFig\textwidth]{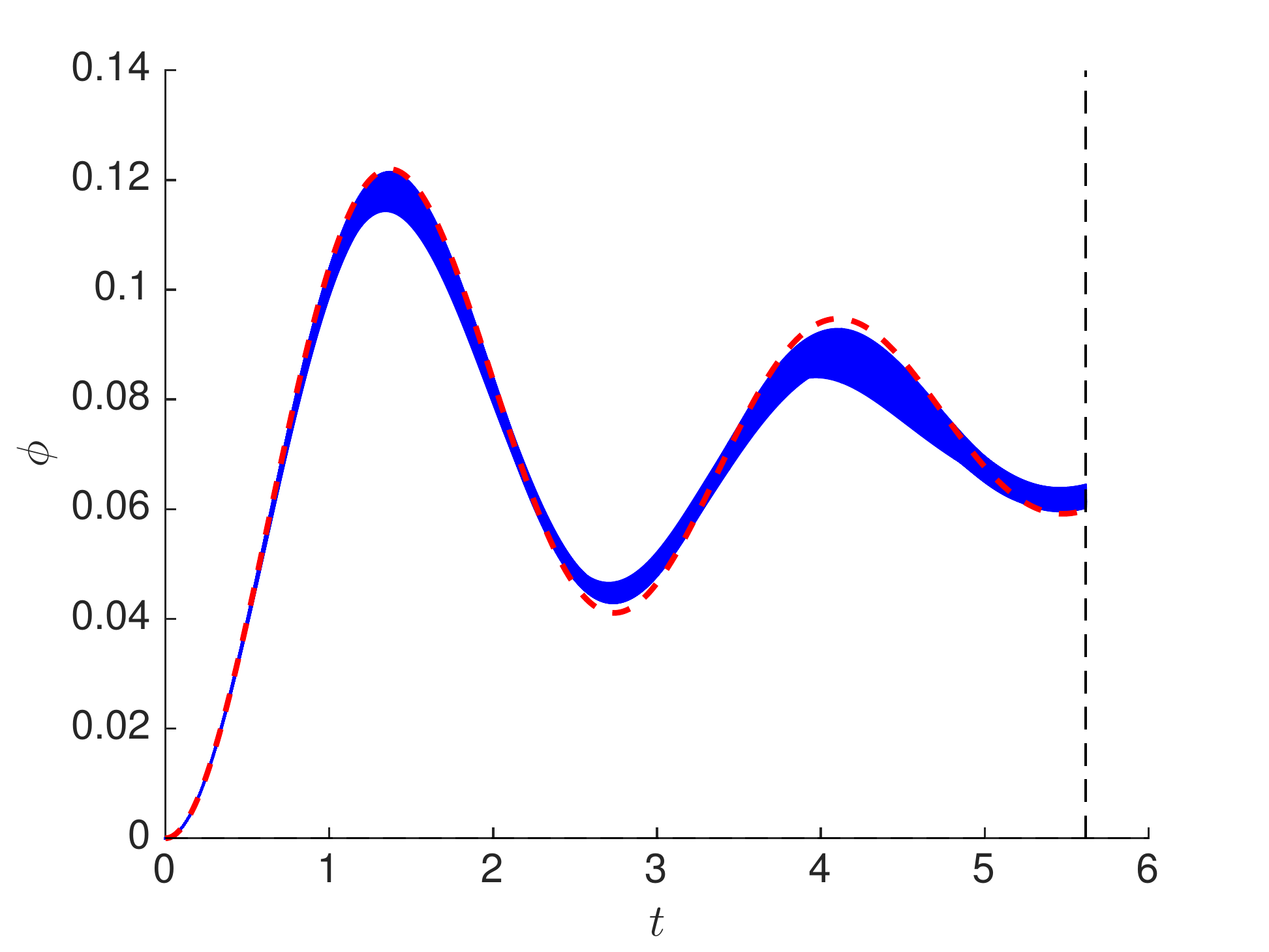}
            \hspace{1em}%
            \includegraphics[width=\sizeFig\textwidth]{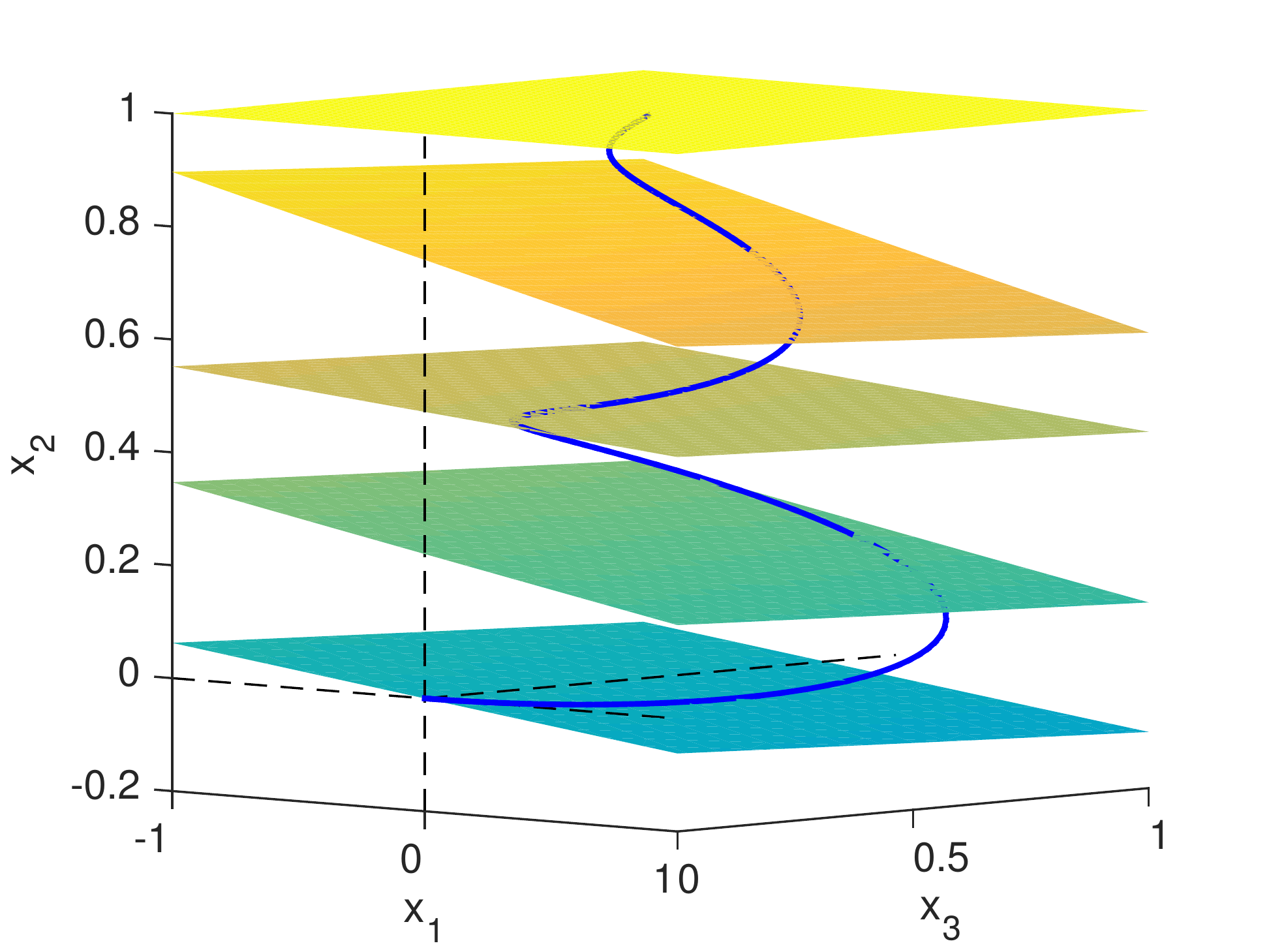}%
            \caption{(Left) The blue plain lines represent the switching functions $\phi(t) = H_1(z(t,z_f))$, $t\in\intervalleff{0}{\tfsol}$,
            for a set of points $z_f \coloneqq (x_f, p_f)$ with $x_f \in M_f$, $p_f=(0,1,0)$, where $M_f$ is sampled using a fine grid,
            and where $z(\cdot,z_f)$ is the solution of ${H_+}$ starting from $z_f$ computed with backward integration. The red dashed line is the switching function
            when the dynamics eq.~\eqref{eq:system} is linear, \ie $\cste_6 = 0$.
            (Right) The evolution of $M_f$ at times $0$, $0.25\, \tfsol$, $0.5\, \tfsol$, $0.75\, \tfsol$ and $\tfsol$, by backward integration.
            The blue curve is the optimal $\gamma_+$ trajectory.}
            \label{fig:gamma+Optimal}
        \end{figure}

        \begin{prpstn}
            For the case study of problem \eqref{ref:OCP_tf_min} described in section \ref{sec:caseStudy}, the strategy $\gamma_+$ is optimal.
            \label{prop:bangOptimal}
        \end{prpstn}
        \begin{proof}
            The final submanifold is $M_f = \intervalleff{-1}{1} \times \{1\} \times \intervalleff{0}{1}$. 
            According to the left subgraph of Figure~\ref{fig:gamma+Optimal}, for all $z_f \in Z_f$ and
            for all $t \in \intervalleof{0}{\tfsol}$, $\phi(t) = H_1(z(t,z_f)) \ne 0$.
            The values of the parameters are given in Table~\ref{table:parametersAndValues} with $\alpha_f=100$, $\imax=1200$ and $\vmax=120$.
        \end{proof}

\section{Geometric analysis of problem \eqref{ref:OCP_tf_min} (with state constraints)}
    \label{sec:constrained}

    \subsection{Abstract formulation}

        In this section, we consider the problem \eqref{ref:OCP_tf_min} with state constraints on $x_1$ and $x_3$, \ie on the electric current
        and the speed of the vehicle. For both state variables $x_i$, $i=1,3$,
        there are two constraints: $x_i\le 1$ and $-1\le x_i$ but these two constraints cannot be active at the same time, that is we cannot
        have $x_i(t)=1$ and $x_i(t)=-1$ at a time $t$.
        Moreover, since $x_1$ and $x_3$ may saturate one of their constraint at the same time only in very particular cases,
        we present the necessary conditions of optimality considering only a scalar constraint of the form $c(x(t)) \le 0$, $\forall\, t\in\intervalleff{0}{t_f}$.
        The optimal control problem \eqref{ref:OCP_tf_min} with a scalar state constraint can be written as the infinite dimensional minimization problem:
        \begin{equation*}
            \min \enstq{t_f}{t_f >  0,\, u(\cdot) \in \Ucal,\, b \circ E_{x_0}(t_f,u(\cdot)) = 0,\, c \circ E_{x_0}(t,u(\cdot)) \le 0,\, t\in\intervalleff{0}{t_f}}.
        \end{equation*}
        It is still possible to express necessary conditions in terms of Lagrange multipliers as before but this has to be done in
        distributions spaces and the Lagrange multipliers must be expressed as Radon measures (see e.~g. \cite{HartlSethiVickson}).

        To emphasize the main difference with the state unconstrained case, we write the following optimization problem with only inequality constraints
        (we omit ``$(\cdot)$'' for readability):
        \begin{equation*}
            \min_{u\in\Ucal} J(u), \quad \text{subject to} \quad S(u) \le 0.
        \end{equation*}
        Here $S$ maps $\Ucal$ into the Banach space of scalar-valued continuous function on $\intervalleff{0}{t_f}$ denoted by $\xCn{0}(\intervalleff{0}{t_f})$
        and supplied with the uniform norm. According to \cite[theorem 1]{Jacobson1971}, if $J$ is a real-valued Fr\'echet differentiable
        function on $\Ucal$ and $S \colon \Ucal \to \xCn{0}(\intervalleff{0}{t_f})$ a Fr\'echet differentiable mapping, then if $\usol \in \Ucal$
        minimizes $J$ subject to $S(\usol) \le 0$, then there exists a scalar $p^0 \le 0$,
        a linear form $\Lambda \in (\xCn{0}(\intervalleff{0}{t_f}))^*$,  $\Lambda \le 0$
        and non-increasing such that the Lagrangian $p^0 J(u) + \prodscal{ \Lambda }{ S(u) }$ is stationary at $\usol$.
        Besides, $\prodscal{ \Lambda }{ S(\usol) } = 0$ and from Riesz's theorem there exists a
        measure\footnote{In our case, $\musol \le 0$ because $p^0 \le 0$. In \cite[theorem 1]{Jacobson1971} $p^0 \ge 0$ and $\musol \ge 0$.} $\musol$ such that
 $           \prodscal{ \Lambda }{ S(\usol) } = \int_0^{t_f} S(\usol) \diff \musol, $
        where the integral is in the Stieltjes sense.

    \subsection{Necessary optimality conditions}

        We recall the necessary conditions due to \cite{Jacobson1971} and \cite{Maurer1977} and follow the presentation of \cite{BFLT2003},
        which exhibits the role of the Lie bracket configuration.

        \subsubsection{Definitions}
        We call a \emph{boundary arc}, labeled $\gamma_c$, an arc defined on an interval $I=\intervalleff{a}{b}$ (not reduced to a singleton),
        such that $c(\gamma_c(t)) = 0$, for every $t\in I$. The times $a$ and $b$ are called the \emph{entry-} and \emph{exit-time} of the
        boundary arc; $a$ and $b$ are also termed \emph{junction times}. An arc $\gamma$ is said to have a \emph{contact point} with the boundary
        at $\tsol \in \intervalleff{0}{t_f}$ if $c(\gamma(\tsol))=0$ and $c(\gamma(t))<0$ for $t\ne\tsol$ in a neighborhood of $\tsol$.
        A subarc $\gamma$ on which $c(\gamma(t)) < 0$ is called an \emph{interior arc}.

        The \emph{generic order} of the constraint is the integer $m$ such that
        $
            F_1 \cdot c = F_1 \cdot (F_0 \cdot c) = \dots = F_1 \cdot (F_0^{m-2} \cdot c) = 0$
            and $F_1 \cdot (F_0^{m-1} \cdot c) \ne 0.
        $
        If the order of a boundary arc $\gamma_c$ is $m$, then its associated feedback control can be generically computed by differentiating $m$
        times the mapping $t \mapsto c(\gamma_c(t))$ and solving with respect to $u$ the linear equation:
        \[
            c^{(m)} = F_0^m \cdot c + u \, F_1 \cdot (F_0^{m-1} \cdot c) = 0.
        \]
        A boundary arc is contained in $c = \dot{c} = \dots = c^{(m)} = 0$, and the constraint $c=0$ is called \emph{primary} while the constraints
        $\dot{c} = \dots = c^{(m)} = 0$ are called \emph{secondary}. The boundary feedback control is denoted by
        \[
            u_c  \coloneqq  - \frac{F_0^m \cdot c }{ F_1 \cdot (F_0^{m-1} \cdot c)}.
        \]

        \subsubsection{Assumptions}
        Let $t\mapsto \gamma_c(t)$, $t\in\intervalleff{t_1}{t_2} \subset \intervalleff{0}{t_f}$,
        be a boundary arc associated with $u_c(\cdot)$. We introduce the assumptions:
        \leqnomode
        \begin{align*}
            \label{Hyp1} \tag{$\mathbf{A_1}$} \quad\quad & \text{$(F_1 \cdot (F_0^{m-1} \cdot c))({\gamma_c(t)}) \ne 0$ 
            for $t\in\intervalleff{t_1}{t_2}$, with $m$ the order of the constraint.}     & \\
            \label{Hyp2} \tag{$\mathbf{A_2}$} \quad\quad & \text{$\abs{u_c(t)}\le1$ for $t\in\intervalleff{t_1}{t_2}$, \ie the boundary control is admissible.}        & \\
            \label{Hyp3} \tag{$\mathbf{A_3}$} \quad\quad & \text{$\abs{u_c(t)} < 1$ for $t\in\intervalleoo{t_1}{t_2}$,
            \ie the control is not saturating on the interior of the boundary arc.} & \\
            \label{Hyp12} \tag{$\mathbf{A^2_1}$} \quad\quad & \text{$(F_{01} \cdot c)({\gamma_c(t)}) \ne 0$ for $t\in\intervalleff{t_1}{t_2}$.}     &
        \end{align*}
        \reqnomode

        \begin{rmrk}
            Assumptions \ref{Hyp1} and \ref{Hyp12} are always satisfied for problem \eqref{ref:OCP_tf_min}, see the proof of proposition \ref{prop:ControlEta}.
            Assumption \ref{Hyp2} is satisfied by any BC-extremals by definition.
            Assumption \ref{Hyp3} is always satisfied for the constraint $c_1$ (see section \ref{sec:ApplicationControleContraint} 
            for the definition of $c_1$) of order 1 since $u_{c_1}(\cdot)$ is strictly increasing (see the second item of section \ref{sec:influenceBounds}
            page \pageref{sec:influenceBounds}). About the constraint $c_3$ of order 2, for some limit cases, $u_{c_3}(\cdot) \equiv 1$.
            This can be seen on Figure \ref{fig:limiteCaseHom3Hom4} page \pageref{fig:limiteCaseHom3Hom4}.
        \end{rmrk}

        \subsubsection{A maximum principle with a single state constraint}
        \label{ref:PMPConstraint}
        Define the pseudo-Hamiltonian:
        \begin{equation*}
            \begin{array}{rcl}
                H \colon T^*M \times U \times \R & \longrightarrow   & \R \\
                (x,p,u,\eta)        & \longmapsto       & H(x,p,u,\eta)  \coloneqq  \prodscal{ p }{ F_0(x) + u\, F_1(x) } + \eta\, c(x),
            \end{array}
        \end{equation*}
        where $\eta$ is the \emph{Lagrange multiplier of the constraint}.
        Consider $(\usol(\cdot), \tfsol) \in \Ucal \times \Rsp$ an optimal solution with associated trajectory $\xsol(\cdot)$.
        Assume that the optimal control is piecewise smooth and
        that along each boundary arc, assumptions \ref{Hyp1} and \ref{Hyp2} are satisfied.
        Then we have the following necessary optimality conditions:
        \begin{enumerate}[i)]
            \item There exists a function $\etasol(\cdot) \le 0$, a real number $p^0\le 0$ and a function 
                $p(\cdot) \in BV(\intervalleff{0}{\tfsol},(\R^n)^*)$ such that:
                \[
                    \dot{\xsol}(t) =  \frac{\partial H}{\partial p}(\xsol(t),\psol(t),\usol(t),\etasol(t)), \quad
                    \dot{\psol}(t) = -\frac{\partial H}{\partial x}(\xsol(t),\psol(t),\usol(t),\etasol(t)), \quad
                    t \in \intervalleff{0}{\tfsol} \text{ a.e.}
                \]
            \item The maximization condition holds for almost every $t \in \intervalleff{0}{\tfsol}$:
                \begin{equation}
                    \displaystyle H(\xsol(t),\psol(t),\usol(t),\etasol(t)) = \max_{u\in U}H(\xsol(t),\psol(t),u,\etasol(t)).
                    \label{eq:maxPbConstraint}
                \end{equation}
            \item The boundary conditions are satisfied and we have the following {transversality conditions}:
                    $
                    \psol(\tfsol) = (0, \mu, 0)$, $\mu \in \R.
                    $
                Since $t_f$ is free, if $\usol(\cdot)$ is continuous at time $\tfsol$, then
                $
                 \displaystyle   H(\xsol(\tfsol),\psol(\tfsol),\usol(\tfsol),\etasol(\tfsol)) = - p^0.
                $
            \item The function $\etasol(\cdot)$ is continuous on the interior of the boundary arcs and
                $$
                 \displaystyle   \etasol(t) \, c(\xsol(t)) = 0, \quad \forall t \in \intervalleff{0}{\tfsol}.
                $$
            \item Let $\Tcal$ denote the set of contact and junction times with the boundary. Then at $\tau \in \Tcal$ we have
                \begin{align}
                    H[\tau^+] = H[\tau^-],   &    \quad \text{ where $[\tau]$ stands for $(\xsol(\tau),\psol(\tau),\usol(\tau),\etasol(\tau))$}, \\
                    \label{eq:adjointJump} \psol(\tau^+)    = \psol(\tau^-) - \nu_\tau \, c'(\xsol(\tau)),  &   \quad \text{ where }
                    \nu_\tau  \coloneqq  \musol(\tau^+) - \musol(\tau^-) \le 0 \text{ (since $\musol$ is non-increasing)}.
                \end{align}
        \end{enumerate}

        \vspace{1em}
        \begin{rmrk} ~
            \begin{itemize}
                \item In this context, an extremal is a quadruple $(x(\cdot),p(\cdot),u(\cdot),\eta(\cdot))$ satisfying i), ii), iv) and v).
                    It is a BC-extremal if iii) also holds.
                \item On a boundary arc, the maximisation condition \eqref{eq:maxPbConstraint} with assumption \ref{Hyp3} imply that $\Phi = 0$
                    on the interior of the boundary arc.
                \item The adjoint vector may be discontinuous at $\tau \in \mathcal{T}$, and $\etasol(\cdot) = {\diff \musol(\cdot)}/{\diff t}$
                on $\intervalleff{0}{\tfsol} \setminus \mathcal{T}$.
            \end{itemize}
        \end{rmrk}

    \subsection{Computations of the multiplier and the jump and the junction conditions}
        \label{sec:BoundaryCase}

        We may find in \cite{BFLT2003,Maurer1977} the determination of the multiplier $\eta$ and the jump $\nu_\tau$ together with the analysis of the junction
        conditions, which is based on the concept of order and related to the classification of extremals. We give next some results for
        $m=1$ and $m=2$ since higher order constraints are not present in the problem \eqref{ref:OCP_tf_min}.

        \subsubsection{Case $m=1$}
        For a first-order constraint, assuming \ref{Hyp1} and \ref{Hyp3}, we have the following lemma from \cite{BFLT2003}.
        \begin{lmm}
            \label{lemma:m1}
            Let $m=1$. Then:
            \begin{enumerate}[i)]
                \item along the boundary, $\displaystyle \eta = \frac{H_{01}}{F_1 \cdot c}$;
                \item if the control is discontinuous at a contact or junction time $\tau$ of a bang arc with the boundary then $\nu_\tau = 0$.
                \item we have
                    \[
                        \nu_\tau = \frac{{\Phi}(\tau^-)}{(F_{1} \cdot c) (x(\tau))} \text{ at an entry point and }
                        \nu_\tau = -\frac{{\Phi}(\tau^+)}{(F_{1} \cdot c) (x(\tau))} \text{ at an exit point.}
                    \]
            \end{enumerate}
            \label{lmm:order1}
        \end{lmm}
        \begin{proof}
            \begin{enumerate}[i)]
                \item Along the boundary, $\Phi = 0$ since \ref{Hyp3} is satisfied; differentiating, we obtain (with a slight abuse of notation)
                        \[
                            0 = \dot{\Phi} = \{H,H_1\} = H_{01} + \eta \, \{c, H_1\} + c \, \{\eta, H_1\} = H_{01} - \eta\, F_1 \cdot c
                        \]
                        and $F_1 \cdot c \ne 0$ along the boundary arc since \ref{Hyp1} holds. Item i) is proved.
                \item See \cite[lemma 2.4]{BFLT2003}.
                \item According to \eqref{eq:adjointJump}, $p(\tau^+) = p(\tau^-) - \nu_\tau \, c'(x(\tau))$
                    at a junction point $x(\tau)$. If $\tau$ is an entry-time, then $\Phi(\tau^+) = H_1(x(\tau^+),p(\tau^+)) = \prodscal{ p(\tau^+) }{ F_1(x(\tau^+)) } = 0$
                    and we have $0 = \Phi(\tau^-) - \nu_\tau (F_1 \cdot c) (x(\tau))$. The proof is similar at an exit-point.
            \end{enumerate}
        \end{proof}

        \begin{rmrk}
            The second item of lemma \ref{lmm:order1} is a particular case of \cite[corollary 5.2]{Maurer1977} which states
            that if $\tau$ is a junction time between a bang arc and a boundary arc and if $m+r$ is odd, where $r$ is the lowest
            order derivative of the optimal control which is discontinuous at $\tau$, then $\nu_\tau = 0$. For this particular case where $m=1$
            and $r=0$, it is still true even for a contact point (see lemma \ref{lmm:order1}) but this is not proved in general in
            \cite[corollary 5.2]{Maurer1977}.
        \end{rmrk}

        We profit from this lemma to get new junction conditions between bang and boundary arcs. We first need to define the sign of an arc.
        To do that, we first consider a non-empty arc $\gamma_b$ defined on the interval $\intervalleff{t_0}{t_1}$. We define the sign of the bang arc by
        $s_b  \coloneqq  \sign (\Phi|_{\gamma_b})$,
        \ie the sign of $\Phi$ on $\intervalleoo{t_0}{t_1}$, which is constant. In a same way, for a first-order constraint $c$, we define the sign of a
        non-empty arc $\gamma_c$ by $s_c  \coloneqq  \sign ((F_1 \cdot c)|_{\gamma_c})$, assuming \ref{Hyp1}.
        Using these definitions, we have the following lemma which gives more insight into the junction conditions for the first-order case and which is
        necessary to define the multiple shooting function, see section \ref{sec:multipleShooting}.
        \begin{lmm}
            \label{lemma:junction1}
            Let consider two non-empty arcs $\gamma_b$ and $\gamma_c$, with $c$ a constraint of order 1, assume \ref{Hyp1} and \ref{Hyp3} along $\gamma_c$
            and note $s_b = \sign (\Phi|_{\gamma_b})$, $s_c = \sign ((F_1 \cdot c)|_{\gamma_c})$.
            We note $\tau$ the junction time between $\gamma_b$ and $\gamma_c$ and $z_\tau \in T^*M$ its associated point.
            \begin{enumerate}[i)]
                \item If the sequence is $\gamma_b \gamma_c$ and
                    \begin{enumerate}[a)]
                        \item if $s_b\, s_c > 0$ or $z_{\tau^-} \in \Sigma^0_{1} $, then $\nu_\tau = 0$, else
                        \item if $s_b\, s_c < 0$ and $z_{\tau^-} \not\in \Sigma^0_{1} $, then $\nu_\tau <  0$ (so the control is continuous at $\tau$).
                    \end{enumerate}
                \item If the sequence is $\gamma_c \gamma_b$ and
                    \begin{enumerate}[a)]
                        \item if $s_b\, s_c < 0$ or $z_{\tau^+} \in \Sigma^0_{1} $, then $\nu_\tau = 0$, else
                        \item if $s_b\, s_c > 0$ and $z_{\tau^+} \not\in \Sigma^0_{1}$, then $\nu_\tau < 0$ (so the control is continuous at $\tau$).
                    \end{enumerate}
            \end{enumerate}
        \end{lmm}
        \begin{proof}
            Let prove first item \textit{i-a)}.
            From lemma \ref{lemma:m1} and equation \eqref{eq:adjointJump}, at time $\tau$ we have
            \[
                 \frac{{\Phi}(\tau^-)}{(F_{1} \cdot c) (x(\tau))} = \nu_\tau \le 0,
            \]
            and if $s_b\, s_c > 0$ then
            \[
                \frac{{\Phi}(\tau^-)}{(F_{1} \cdot c) (x(\tau))} \ge 0,
            \]
            thus $\nu_\tau = 0$ and $\Phi$ is continuous at $\tau$ and $\Phi(\tau) = 0$.
            Item \textit{i-a)} is proved. Now, if $s_b\, s_c < 0$ and $z_{\tau^-} \not\in \Sigma^0_{1} $, then ${\Phi}(\tau^-) \ne 0$ and then $\nu_\tau < 0$.
            Item \textit{i-b)} is proved. Items \textit{ii-a)} and \textit{ii-b)} are proved with the same argumentation.
        \end{proof}

        \begin{rmrk}
            In the numerical results, see section \ref{sec:numericalResults}, we have trajectories of the form $\gamma_+ \gamma_{c_1} \gamma_+$, with $F_1 \cdot c_1 > 0$.
            By lemma \ref{lemma:junction1}, the jump is zero at the entry-time of $\gamma_{c_1}$ and the control is continuous at
            the exit-time $\tau$ of $\gamma_{c_1}$ if at $t = \tau^+$, the extremal has no contact with the switching surface.
        \end{rmrk}

        \subsubsection{Case $m=2$}
        For a second-order constraint, assuming \ref{Hyp1}, \ref{Hyp12} and \ref{Hyp3}, we have the following
        lemma.
        \begin{lmm}
            \label{lemma:m2}
            Let $m=2$. Then:
            \begin{enumerate}[i)]
                \item along a boundary arc, we have
                    $\displaystyle
                        \eta = \frac{H_{001} + u_c \, H_{101}}{F_{01} \cdot c};
                    $
                \item at a contact or junction point at time $\tau$, we have
                    $
                        \Phi(\tau^+) = \Phi(\tau^-);
                    $
                \item we have
                    \[
                        \nu_\tau = \frac{\dot{\Phi}(\tau^-)}{(F_{01} \cdot c) (x(\tau))}  \text{ at an entry point and }
                        \nu_\tau = -\frac{\dot{\Phi}(\tau^+)}{(F_{01} \cdot c) (x(\tau))} \text{ at an exit point.}
                    \]
            \end{enumerate}
        \end{lmm}

        For an order 2 constraint $c$, we define the sign of a non-empty arc $\gamma_c$ by $s_c  \coloneqq  \sign ((F_{01} \cdot c)|_{\gamma_c})$, assuming \ref{Hyp12}.
        \begin{lmm}
            \label{lemma:junction2}
            Let consider two non-empty arcs $\gamma_b$ and $\gamma_c$, with $c$ a constraint of order 2, assume \ref{Hyp1},
            \ref{Hyp12} and \ref{Hyp3} along $\gamma_c$ and note $s_b = \sign (\Phi|_{\gamma_b})$, $s_c = \sign ((F_{01} \cdot c)|_{\gamma_c})$.
            We note $\tau$ the junction time between $\gamma_b$ and $\gamma_c$ and $z_\tau \in T^*M$ its associated point.
            For both sequences $\gamma_b \gamma_c$ and $\gamma_c \gamma_b$,
            if $s_b\, s_c < 0$ or $z_{\tau^-} \in \Sigma^0_{01} $, then $\nu_\tau = 0$, else
            if $s_b\, s_c > 0$ and $z_{\tau^-} \not\in \Sigma^0_{01} $, then $\nu_\tau <  0$.
        \end{lmm}
        \begin{proof}
            Just notice that the signs of $\Phi$ and $\dot{\Phi}$ are equal after an exit point and are opposite before an entry point. 
            Then, the proof is similar to the proof of lemma \ref{lemma:junction1}.
        \end{proof}

        \begin{rmrk}
            In the numerical results, see section \ref{sec:numericalResults}, we have trajectories with sequences of the form 
            $\gamma_- \gamma_{c_3}$, with $F_{01} \cdot c_3 < 0$.
            By lemma \ref{lemma:junction2}, the jump is non zero at the junction if at this time, the extremal has no contact of order 2 with the switching surface.
        \end{rmrk}

        \subsubsection{Application to problem \eqref{ref:OCP_tf_min}}
        \label{sec:ApplicationControleContraint}

        Let define the space $C  \coloneqq  C_1 \cup C_{3} \cup C_{13}$ where the boundary sets are
        $C_{1}  \coloneqq   \enstq{ x \in M} {\abs{x_1}   = 1,~ \abs{x_3} \ne 1  }$, $C_{3}  \coloneqq   \enstq{ x \in M} {\abs{x_1} \ne 1,~ \abs{x_3}   = 1  }$ and $C_{13}   \coloneqq   \enstq{ x \in M} {\abs{x_1} = \abs{x_3} = 1}$.
        For each state constraint, we write $c_i(x)  \coloneqq  s_i \, x_i - 1$, $i \in \{1,3\}$, $s_i \in \{-1,1\}$,
        and we label $\gamma_{c_\alpha}$, $\alpha \in \{1,3,13\}$, a boundary arc defined on $I$ such that for every $t \in I$, $\gamma_{c_\alpha}(t) \in C_\alpha$.
        Assuming \ref{Hyp3}, we have the following proposition.
        \begin{prpstn}
            \label{prop:ControlEta} We have:
            \begin{enumerate}[i)]
                \item the constraint on the electric current, \ie on $x_1$, is of order 1. The feedback boundary control and the multiplier associated
                    to a boundary arc $\gamma_{c_1}$ are given by:
                    \begin{align*}
                        u_{c_1}(x(t))      \coloneqq  -\frac{1}{\cste_7} \left( \cste_1 s_1 + \cste_2 x_3(t) \right), \quad 
                        \eta_{c_1}(z(t))   \coloneqq  - s_1 \cste_5 \, p_3(t).
                    \end{align*}
                \item the constraint on the angular velocity, \ie on $x_3$, is of order 2. The feedback boundary control and the multiplier associated
                    to a boundary arc $\gamma_{c_3}$ are constant and given by:
                    \begin{align*}
                        u_{c_3}(x(t))      \coloneqq  -\frac{1}{\cste_7} \left( - \frac{\cste_1 \left( \cste_4 + \cste_6 \right)}{\cste_5} + \cste_2 s_3 \right), \quad 
                        \eta_{c_3}(z(t))   \coloneqq  - s_3 \cste_3 \, p_2.
                    \end{align*}
                    Besides, on the boundary arc $\gamma_{c_3}$ the electric current is constant and strictly positive and we have
                    \[
                        x_1(\cdot) = - \frac{\cste_4 + \cste_6}{\cste_5} > 0.
                    \]
                \item There exists boundary arcs $\gamma_{c_{13}}$ only if $\cste_4 + \cste_5 + \cste_6 = 0$ and if $\gamma_{c_{13}}$ is a boundary
                    arc, then $s_1 = 1$ along $\gamma_{c_{13}}$.
            \end{enumerate}
        \end{prpstn}
        \begin{proof}
            \begin{enumerate}[i)]
                \item Let consider a boundary arc $\gamma_{c_1}$ defined on $I$. Then we have for every $t\in I$,
                    \begin{align*}
                        0 = \frac{\diff }{\diff t} c_1(\gamma_{c_1}(t)) & = (F \cdot c_1)({\gamma_{c_1}(t)})
                            = s_1 \left( \cste_1 x_1(t) + \cste_2 x_3(t) \right) + u(t)\, s_1\, \cste_7, 
                    \end{align*}
                    with $\left( F_1 \cdot c_1 \right)({\gamma_{c_1}(t)}) = s_1 \cste_7 \ne 0$ (so \ref{Hyp1} is satisfied) and $x_1(t) = s_1$.
                    Besides, $\eta_{c_1} = H_{01} / F_1 \cdot c_1$, \cf\ lemma \ref{lemma:m1} with $H_{01}$ given at lemma \ref{lemma:lieBrackets}
                    and $p_1 = 0$ along $\gamma_{c_1}$ since \ref{Hyp3} is assumed.
                    Whence the results.
                \item Let consider a boundary arc $\gamma_{c_3}$ defined on $I$. Then we have for every $t\in I$, $x_3(t) = s_3$ and
                    \begin{align*}
                        0 = \frac{\diff }{\diff t} c_3(\gamma_{c_3}(t)) & = (F \cdot c_3)({\gamma_{c_3}(t)})
                            = s_3 \left( \cste_4 + \cste_5 x_1(t) + \cste_6 x_3^2(t)  \right)
                            = s_3 \left( \cste_4 + \cste_5 x_1(t) + \cste_6            \right).
                    \end{align*}
                    The electric current is constant and strictly positive since
                    $\cste_4 + \cste_5 x_1(t) + \cste_6 = 0$, with $\cste_4 < 0$, $\cste_5 > 0$ and $\cste_6 < 0$.
                    Differentiating a second time we have
                    \begin{align*}
                        0= \frac{\diff^2}{\diff t^2} c_3(\gamma_{c_3}(t)) & = s_3 \cste_5 \, \dot{x}_1(t)
                            = s_3 \cste_5 \left( \cste_1 x_1(t) + \cste_2 x_3(t) + u(t) \cste_7 \right), \quad x_3(t) = s_3,
                    \end{align*}
                    with $\left( F_1 \cdot (F_0 \cdot c_3) \right)({\gamma_{c_3}(t)}) = s_3 \cste_5 \cste_7 \ne 0$, so \ref{Hyp1} and \ref{Hyp12} are satisfied since here
                    $F_{01} \cdot c = - F_{1} \cdot (F_{0} \cdot c)$.
                    The result follows with $\eta_{c_3}$ given at lemma \ref{lemma:m2}, $H_{001}$, $H_{101}$ and $F_{01}$ given at lemma \ref{lemma:lieBrackets},
                    and with $p_1=p_3=0$ along $\gamma_{c_3}$ since \ref{Hyp3} holds.
                \item According to $i)$ and $ii)$, on a boundary arc $\gamma_{c_{13}}$, $x_1 = s_1$ and $\cste_4 + \cste_5 x_1 + \cste_6 = 0$,
                    $x_1 > 0$. Whence the result.
            \end{enumerate}
        \end{proof}

    \subsection{Local time minimal synthesis}
        \label{sec:timeMinimalSynthesis}

        We end this section with an application of theorem 4.4 from \cite{BFLT2003}.
        This theorem is local and valid only in the parabolic case.
        It asserts the following.
        One consider the time minimization problem for a 3-dimensional system of the form \eqref{eq:system}, $\abs{u}\le1$,
        with the scalar constraint $c(x)\le 0$. Let $\xsol\in c^{-1}(\{0\})$ and assume the following:
        \begin{enumerate}
            \item at the point $\xsol$, the vectors $F_0$, $F_1$ and $F_{01}$ form a frame and
                \[
                    [F_0 \pm F_1, F_{01} ](\xsol) = a\, F_0(\xsol) + b\, F_1(\xsol) + c\, F_{01}(\xsol),
                \]
                where $a>0$;
            \item the constraint is of order 2 and assumptions \ref{Hyp1} and \ref{Hyp3} hold at the point $\xsol$.
        \end{enumerate}
        Then the boundary arc passing through $\xsol$ is small time optimal if and only if the arc $\gamma_+$ is contained in the nonadmissible domain $c \ge 0$.
        In this case, the local time minimal synthesis with a boundary arc is of the form $\gamma_+\gamma_-^T\gamma_c\gamma_-^T\gamma_+$, where $\gamma_-^T$
        are arcs tangent to the boundary arc. Therefore, we have the following proposition.
        \begin{prpstn}
            Let $\xsol \in c_3^{-1}(\{0\})$ with $s_3 = 1$, and assume \ref{Hyp3} at the point $\xsol$.
            Assume also that the arc $\gamma_+$ passing through $\xsol$ is contained in the nonadmissible domain $c_3 \ge 0$.
            Then the local time minimal synthesis is of the form $\gamma_+\gamma_-^T\gamma_{c_3}\gamma_-^T\gamma_+$.
        \end{prpstn}
        \begin{proof}
            The constraint $c_3$ is of order 2, assumption \ref{Hyp1} holds since $\left( F_1 \cdot (F_0 \cdot c_3) \right)(\xsol) = \cste_5 \cste_7 \ne 0$
            and $[F_0 \pm F_1, F_{01} ](\xsol) = F_{001}(\xsol) = a\, F_0(\xsol) + b\, F_1(\xsol) + c\, F_{01}(\xsol)$, with $a=\cste_5 \cste_7 / \xsol_3 > 0$.
        \end{proof}

\section{Numerical methods}
\label{sec:numericalMethods}

    In this section, we present the numerical methods used to find BC-extremals of problem \eqref{ref:OCP_tf_min}. All the methods are implemented within the
    \hampath\ code, see \cite{OMS}. The software is based upon indirect methods: simple, multiple shooting,
    differential path following (or homotopy) methods and exponential mappings.
    One focus on the description of the shooting and homotopic functions, the methods being automatically generated by the \hampath\ code.
    One pay a special attention in section \ref{sec:multipleShooting} on the strategy we develop to define the shooting functions
    associated to problems with state constraints.

    \begin{rmrk}
        From now on, the state constraints are $c_1(x) = x_1 - 1$ and $c_3(x) = x_3 - 1$.
    \end{rmrk}

    \begin{rmrk}
        We only present the shooting functions associated to the structures we encounter in section \ref{sec:influenceBounds} during numerical experiments.
    \end{rmrk}

    \subsection{Simple shooting method}

        \subsubsection{Structure $\gamma_+$}
        If $\imax$ and $\omega_\mathrm{max}$ (or $\vmax = \omega_\mathrm{max} \times 3.6 \, r / K_r$) are big enough the solution of
        problem \eqref{ref:OCP_tf_min} has only interior arcs and then the optimal trajectory is of the form $\gamma_+$, see proposition \ref{prop:bangOptimal}.
        We need the following definitions.
        \begin{dfntn}
            For fixed $\zsol_0 \in T^*M$ and $\tsol \ge 0$, we define in a neighborhood of $(\zsol_0,\tsol)$ (if possible),
            the following \emph{exponential mapping}
            $(z_0,t) \mapsto \expmap{z_0}{t}{\vvec{H}}$ as the trajectory $z(\cdot)$ at time $t$ of
            the Hamiltonian vector $\vvec{H}$, \ie $\dot{z}(s) = \vvec{H}(z(s))$ for every $s\in \intervalleff{0}{t}$, satisfying $z(0) = z_0$.
        \end{dfntn}
        \begin{dfntn}
            Let $\vvec{H}$ be a Hamiltonian vector on $T^*M$, and let $z(\cdot)$ be a trajectory of $\vvec{H}$ defined on $\intervalleff{0}{t_f}$.
            The differential equation on $\intervalleff{0}{t_f}$
            \[
                \dot{\wideparen{\delta z}}(t) = \diff \vvec{H} (z(t)) \cdot \delta z(t)
            \]
            is called a \emph{Jacobi equation}, or \emph{variational system}, along $z(\cdot)$. Let $\delta z(\cdot)$ be a solution of
            the variational system along $z(\cdot)$, we write $\delta z(t) =:  \expmap{\delta z(0)}{t}{\diff\vvec{H}|_{z(\cdot)} }$.
        \end{dfntn}
        Let $(z(\cdot),u(\cdot))$ be a regular extremal defined on $\intervalleff{0}{t_f}$ with only one single positive bang arc.
        We have $z(t) = \expmap{z(0)}{t}{\vvec{H_+}}$, $t\in \intervalleff{0}{t_f}$, where $H_+(z) = H(z,u(z)) = H_0(z) + H_1(z)$.
        Let denote by $\pi$ the following projection: in coordinates, writing $z  \coloneqq  (x,p) \in T^*M$, then $\pi(z)  \coloneqq  (x_2,p_1,p_3)$.
        We define now the \emph{simple shooting function}
        \[
            S_1(p_0,t_f)  \coloneqq
                            \begin{pmatrix}
                                \pi\left( \expmap{z_0}{t_f}{\vvec{H_+}} \right) \\[0.5em]
                                H_+\left( \expmap{z_0}{t_f}{\vvec{H_+}} \right)
                            \end{pmatrix}
                        -
                            \begin{pmatrix}
                                1 \\[0.0em]
                                0 \\[0.0em]
                                0 \\[0.0em]
                                1
                            \end{pmatrix},
        \]
        where $z_0 \coloneqq (x_0,p_0)$, with $x_0 = (0,0,0)$ fixed. If $(\psol_0,\tfsol)$ is a zero of $S_1$, then the constant control $u(\cdot) = +1$
        with the integral curve
        $t \mapsto \expmap{\zsol_0}{t}{\vvec{H_+}}$, with $\zsol_0  \coloneqq  (x_0,\psol_0)$, for $t \in \intervalleff{0}{\tfsol}$, is a BC-extremal,
        \ie the extremal satisfies the necessary optimality conditions, see section \ref{sec:PMP}.
        The \emph{simple shooting method} consists in finding a zero of the simple shooting function $S_1$, \ie in solving $S_1(p_0,t_f) = 0$.

        \subsubsection{\hampath\ code}
        \label{sec:hampathCode}
        The Fortran hybrid Newton method \hybrj{} (from the \minpack{} library \cite{minpack80}) is used
        to solve the nonlinear system $S_1(p_0,t_f)=0$. Providing $H_+$ and $S_1$ to \hampath,
        the code generates automatically the Jacobian of the shooting function given to the solver. To make the implementation
        of $S_1$ easier, \hampath{} supplies the exponential mapping. Automatic Differentiation (\tapenade{} software \cite{tapenade2012})
        is used to produce $\vvec{H_+}$ and is combined with Runge-Kutta integrators (\dopri, \dop, see \cite{HaNoWa93} and \radau, see \cite{HaWa96})
        to assemble the exponential mapping.

        We detail how the Jacobian of the shooting function is computed. If we note $z(t,z_0)  \coloneqq  \expmap{z_0}{t}{\vvec{H_+}}$, then
        the Jacobian is given by
        \begin{equation*}
            \begin{aligned}
                \frac{\partial S_1}{\partial p_0}(p_0,t_f) \, \delta p_0 & =
                    \begin{pmatrix}
                        \displaystyle \diff \pi\left( z(t_f,z_0) \right) \cdot \frac{\partial z}{\partial z_0} (t_f,z_0)    \, \delta z_0 \\[0.5em]
                        \displaystyle \diff H_+\left( z(t_f,z_0) \right) \cdot \frac{\partial z}{\partial z_0} (t_f,z_0) \, \delta z_0
                    \end{pmatrix},
                    \quad \delta z_0  \coloneqq  \begin{pmatrix} 0_{\R^n} \\ \delta p_0 \end{pmatrix},
                    \quad \delta p_0 \in \R^n,
                    \\[1em]
                \frac{\partial S_1}{\partial t_f}(p_0,t_f) \, \delta t_f & =
                           \delta t_f
                    \begin{pmatrix}
                        \displaystyle \diff \pi\left( z(t_f,z_0) \right) \cdot \vvec{H_+} (z(t_f,z_0)) \\[0.5em]
                        \displaystyle \diff H_+\left( z(t_f,z_0) \right) \cdot \vvec{H_+} (z(t_f,z_0))
                    \end{pmatrix},
                    \quad \delta t_f \in \R,
            \end{aligned}
        \end{equation*}
        where
        \[
            \frac{\partial z}{\partial z_0} (t_f,z_0) \, \delta z_0 = \expmap{\delta z_0}{t_f}{\diff\vvec{H_+} |_{z(\cdot,z_0)} },
        \]
        and $\diff H_+\left( z(t_f,z_0) \right) \cdot \vvec{H_+} (z(t_f,z_0)) = \{H_+,H_+\}(z(t_f,z_0)) = 0$.
        To assemble automatically the Jacobian, \hampath\ uses AD to compute $\diff \pi$, $\diff H_+$ and $\diff \vvec{H_+}$, and produces the exponential mapping
        associated to the variational system.

    \subsection{Multiple shooting method}
    \label{sec:multipleShooting}

    We need the following propositions (inspired by \cite[lemma 20.21]{Agrachev2004} and \cite[proposition 1]{Caillau2012})
        to define the shooting functions when the solutions have boundary arcs.
        We introduce for that the canonical $x$-projection: $\pi_x (x,p) \coloneqq  x$, $(x,p) \in T^*M$.
        \begin{prpstn}
            \label{prop:invariant1}
            Let $c(x)\le 0$ be a scalar constraint of order 1, and define for $z \coloneqq (x,p) \in T^*M$ the true Hamiltonian
            \[
                H_{c}(z)  \coloneqq  H(z,u_{c}(x),\eta_{c}(z)) = H_0(z) + u_{c}(x)\, H_1(z) + \eta_{c}(z)\, c(x),
            \]
            with
            \[
                u_c(x) = - \frac{(F_0 \cdot c) (x) }{ (F_1 \cdot c) (x)}, \quad
                \eta_c(z) = \frac{H_{01}(z)}{(F_1 \cdot c) (x)}.
            \]
            Let $\zsol  \coloneqq  (\xsol,\psol) \in \Sigma^0_1$, $c(\xsol) = 0$; there is exactly one extremal $(x(\cdot),p(\cdot),u(\cdot),\eta(\cdot))$
            passing through $\zsol$, such that $c \circ x(\cdot) = 0$, $H_1 \circ z(\cdot) = 0$, $z(\cdot) \coloneqq (x(\cdot), p(\cdot))$,
            and it is defined by the flow of $H_{c}$.
        \end{prpstn}
        \begin{proof}
            First we show that the space $\enstq{ x \in M } {c(x) = 0 } \cap \Sigma^0_1$ is invariant with respect to the flow of $H_{c}$.
            Let $\zsol  \coloneqq  (\xsol,\psol) \in \Sigma^0_1$, $c(\xsol) = 0$, and let $z(\cdot)  \coloneqq  (x(\cdot),p(\cdot))$ be the
            associated integral curve of $H_{c}$ passing through $\zsol$ at time $0$.
            Let $\vphi  \coloneqq  (g, H_1) \circ z(\cdot)$, with $g  \coloneqq  c \circ \pi_x$; then $\vphi$ is differentiable and
            \begin{align*}
                \dot{\vphi}_1(t)  & = \{H_{c}, g\} (z(t)) \\[0.5em]
                & = \Big(
                    \underbrace{ \{ H_0                 , g \} + \{ H_1, g\} \, (u_{c}\circ \pi_x)  }_{ = \, 0 \text{ by definition of } u_{c}}
                +   \underbrace{ \{ u_{c} \circ \pi_x , g \}}_{= \, 0} \, H_1
                +   \underbrace{ \{ c \circ \pi_x     , g \}}_{= \, 0} \, \eta_{c}
                +                \{\eta_{c}           , g \} \, (c\circ \pi_x)            \Big) ({z(t)}). \\[0.5em]
                \dot{\vphi}_2(t)  & = \{H_{c}, H_1\} (z(t)) \\[0.5em]
                & = \Big(
                    \underbrace{ H_{01} + \{ c \circ \pi_x, H_1\} \, \eta_{c}  }_{ = \, 0 \text{ by definition of } \eta_{c}}
                +   \underbrace{ \{ H_1                 , H_1 \}}_{= \, 0} \, (u_{c} \circ \pi_x)
                +   \{ u_{c} \circ \pi_x , H_1 \} \, H_1
                +                \{\eta_{c}           , H_1 \} \, (c\circ \pi_x)            \Big) ({z(t)}). \\
            \end{align*}
            so $\dot{\vphi}(t) = A(t) \, \vphi(t)$, with
            \[
                A(t)  \coloneqq
                \begin{bmatrix*}[l]
                    \{\eta_{c}, c \circ \pi_x   \}  &   0                               \\
                    \{\eta_{c}, H_1 \}                &   \{ u_{c} \circ \pi_x , H_1 \}
                \end{bmatrix*}(z(t)).
            \]
            Since $\vphi(0) = 0$, $\vphi$ is indentically zero and $z(\cdot)$ remains in $\enstq{ x \in M }{ c(x) = 0 } \cap \Sigma^0_1$.
            Now,
            \[
                H'_{c}(z) =
                \frac{\partial H}{\partial z}       (z,u_{c}(x),\eta_{c}(z)) +
                \frac{\partial H}{\partial u}       (z,u_{c}(x),\eta_{c}(z)) \, (u_{c} \circ \pi_x)'(z) +
                \frac{\partial H}{\partial \eta}    (z,u_{c}(x),\eta_{c}(z)) \, \eta'_{c}(z),
            \]
            with
            \[
                \frac{\partial H}{\partial u}(z,u_{c}(x),\eta_{c}(z)) = H_1(z) \quad \text{and} \quad
                \frac{\partial H}{\partial \eta}(z,u_{c}(x),\eta_{c}(z)) = c(x),
            \]
            so $\vvec{H_c}(z(t)) = \vvec{H}(z(t),u_{c}(x(t)),\eta_{c}(z(t)))$ as $H_1$ and $c \circ \pi_x$ vanish along $z(\cdot)$,
            and $(z(\cdot), u_{c} \circ x(\cdot), \eta_{c} \circ z(\cdot))$ is extremal.
        \end{proof}

        For the second-order case, we have the following proposition.
        \begin{prpstn}
            \label{prop:invariant2}
            Let $c(x)\le 0$ be a scalar constraint of order 2, and define for $z \coloneqq (x,p) \in T^*M$ the true Hamiltonian
            \[
                H_{c}(z)  \coloneqq  H(z,u_{c}(x),\eta_{c}(z)) = H_0(z) + u_{c}(x)\, H_1(z) + \eta_{c}(z)\, c(x),
            \]
            with
            \[
                u_c(x) = - \frac{(F_0^2 \cdot c) (x) }{ (F_1 \cdot (F_0 \cdot c)) (x)}, \quad
                \eta_c(z) = \frac{H_{001}(z) + u_c(x) \, H_{101}(z)}{(F_{01} \cdot c)(x)}.
            \]
            \begin{enumerate}[i)]
                \item   Assume that the control is constant, \ie $\forall x\in M$, $u'_c(x) = 0$.
                        Let $\zsol  \coloneqq  (\xsol,\psol) \in T^*M$, $c(\xsol) = 0$ and $(F_0 \cdot c)(\xsol) = 0$;
                        there is exactly one extremal $(x(\cdot),p(\cdot),u(\cdot),\eta(\cdot))$
                        passing through $\zsol$, such that
                        \[
                            c \circ x(\cdot) = 0 \quad \text{and} \quad  (F_0 \cdot c) \circ x(\cdot) = 0,
                        \]
                        and it is defined by the flow of $H_{c}$.
                \item[]
                \item   Let $\zsol  \coloneqq  (\xsol,\psol) \in \Sigma^0_1 \cap \Sigma^0_{01}$, $c(\xsol) = 0$ and $(F_0 \cdot c)(\xsol) = 0$;
                        there is exactly one extremal $(x(\cdot),p(\cdot),u(\cdot),\eta(\cdot))$
                        passing through $\zsol$, such that
                        \[
                            c \circ x(\cdot) = 0, \quad  (F_0 \cdot c) \circ x(\cdot) = 0,
                            \quad H_1 \circ z(\cdot) = 0 \quad \text{and} \quad H_{01} \circ z(\cdot) = 0,
                        \]
                        $z(\cdot)  \coloneqq  (x(\cdot),p(\cdot))$, and it is defined by the flow of $H_{c}$.
            \end{enumerate}
        \end{prpstn}
        \begin{proof}
            First we show that the spaces
            \[
                E_1   \coloneqq  \enstq{ x \in M }{ c(x) = 0 } \cap \enstq{ x \in M }{ (F_0 \cdot c)(x) = 0 } \quad \text{and} \quad 
                E_2   \coloneqq  E_1 \cap \Sigma^0_1 \cap \Sigma^0_{01} = E_1 \cap \Sigma_s
            \]
            are invariant with respect to the flow of $H_c$. Let $\zsol  \coloneqq  (\xsol, \psol) \in T^*M$, and let $z(\cdot)  \coloneqq  (x(\cdot),p(\cdot))$ be the
            associated integral curve of $H_{c}$ passing through $\zsol$ at time $0$.
            Let $\vphi  \coloneqq  (g, f, H_1, H_{01}) \circ z(\cdot)$, with $g  \coloneqq  c \circ \pi_x$ and $f  \coloneqq  (F_0 \cdot c) \circ \pi_x$;
            then $\vphi$ is differentiable and we have the following.
            \begin{align*}
                \dot{\vphi}_1(t)  & = \{H_{c}, g\} (z(t)) \\
                & = \Big(
                                    \{ H_0                  , g \}
                +   \underbrace{    \{ H_1                  , g \}  \, (u_{c}\circ \pi_x)   }_{ = \, 0 \text{ ($c$ of order 2)}}
                +   \underbrace{    \{ u_{c} \circ \pi_x    , g \}                          }_{= \, 0} \, H_1
                +   \underbrace{    \{ c \circ \pi_x        , g \}                          }_{= \, 0} \, \eta_{c}
                +                   \{\eta_{c}              , g \}  \, (c\circ \pi_x)
                    \Big) ({z(t)})  \\
                & = \Big(
                                    (F_0 \cdot c) \circ \pi_x + \{\eta_{c}, g \} \, (c\circ \pi_x)
                \Big) ({z(t)}), \\[0.5em]
                \dot{\vphi}_2(t)  & = \{H_{c}, f\} (z(t)) \\
                & = \Big(
                    \underbrace{    \{ H_0                  , f \}
                +                   \{ H_1                  , f \}  \, (u_{c}\circ \pi_x)   }_{ = \, 0 \text{ by definition of $u_c$}}
                +   \underbrace{    \{ u_{c} \circ \pi_x    , f \}                          }_{= \, 0} \, H_1
                +   \underbrace{    \{ c \circ \pi_x        , f \}                          }_{= \, 0} \, \eta_{c}
                +                   \{\eta_{c}              , f \}  \, (c\circ \pi_x)
                    \Big) ({z(t)}), \\[0.5em]
                \dot{\vphi}_3(t)  & = \{H_{c}, H_1\} (z(t)) \\
                & = \Big(
                                                              H_{01}
                +   \underbrace{    \{ H_1                  , H_1 \}                        }_{= \, 0}  \, (u_{c}\circ \pi_x)
                +                   \{ u_{c} \circ \pi_x    , H_1 \}                                    \, H_1
                +   \underbrace{    \{ c \circ \pi_x        , H_1 \}                        }_{= \, 0}  \, \eta_{c}
                +                   \{\eta_{c}              , H_1 \}  \, (c\circ \pi_x)
                    \Big) ({z(t)}), \\[0.5em]
                \dot{\vphi}_4(t)  & = \{H_{c}, H_{01}\} (z(t)) \\
                & = \Big(
                    \underbrace{                              H_{001}
                +                                             H_{101}   \, (u_{c}\circ \pi_x)
                +                   \{ c \circ \pi_x        , H_{01} \} \, \eta_{c}             }_{= \, 0 \text{ by definition of $\eta_c$}}
                +                   \{ u_{c} \circ \pi_x    , H_{01} \} \, H_1
                +                   \{ \eta_{c}             , H_{01} \} \, (c \circ \pi_x)
                    \Big) ({z(t)}), \\
            \end{align*}
            so $\dot{\vphi}(t) = A(t) \, \vphi(t)$, with
            \[
                A(t)  \coloneqq
                \begin{bmatrix*}[l]
                    \{\eta_{c}, g       \}      &   1   &   0                                   &   0   \\
                    \{\eta_{c}, f       \}      &   0   &   0                                   &   0   \\
                    \{\eta_{c}, H_1     \}      &   0   &   \{ u_{c} \circ \pi_x, H_1    \}     &   1   \\
                    \{\eta_{c}, H_{01}  \}      &   0   &   \{ u_{c} \circ \pi_x, H_{01} \}     &   0
                \end{bmatrix*}(z(t)),
            \]
            and we have also
            \[
                \begin{bmatrix*}[l]
                    \dot{\vphi}_1(t) \\
                    \dot{\vphi}_2(t)
                \end{bmatrix*}
                =
                \begin{bmatrix*}[l]
                    \{\eta_{c}, g       \}(z(t))      &   1   \\
                    \{\eta_{c}, f       \}(z(t))      &   0   \\
                \end{bmatrix*} ~
                \begin{bmatrix*}[l]
                    {\vphi}_1(t) \\
                    {\vphi}_2(t)
                \end{bmatrix*}.
            \]
            So if $\vphi_1(0) = \vphi_2(0) = 0$, \ie $\zsol \in E_1$, then $\vphi_1$ and $\vphi_2$ are indentically zero and $z(\cdot)$ remains in $E_1$.
            Now if $\vphi(0) = 0$, \ie $\zsol \in E_2$, then $\vphi$ is indentically zero and $z(\cdot)$ remains in $E_2$.
            Now,
            \[
                H'_{c}(z) =
                \frac{\partial H}{\partial z}       (z,u_{c}(x),\eta_{c}(z)) +
                \frac{\partial H}{\partial u}       (z,u_{c}(x),\eta_{c}(z)) \, (u_{c} \circ \pi_x)'(z) +
                \frac{\partial H}{\partial \eta}    (z,u_{c}(x),\eta_{c}(z)) \, \eta'_{c}(z),
            \]
            with
            \[
                \frac{\partial H}{\partial u}(z,u_{c}(x),\eta_{c}(z)) = H_1(z) \quad \text{and} \quad
                \frac{\partial H}{\partial \eta}(z,u_{c}(x),\eta_{c}(z)) = c(x),
            \]
            so
            \begin{itemize}
                \item   if $\zsol \in E_1$ and if we assume $(u_{c} \circ \pi_x)' = 0$ then
                        $\vvec{H_c}(z(t)) = \vvec{H}(z(t),u_{c}(x(t)),\eta_{c}(z(t)))$ as $z(\cdot)$ remains in $E_1$,
                        and $(z(\cdot), u_{c} \circ x(\cdot), \eta_{c} \circ z(\cdot))$ is extremal. Item i) is proved.
                \item   if $\zsol \in E_2$ then
                        $\vvec{H_c}(z(t)) = \vvec{H}(z(t),u_{c}(x(t)),\eta_{c}(z(t)))$ as $z(\cdot)$ remains in $E_2$,
                        and $(z(\cdot), u_{c} \circ x(\cdot), \eta_{c} \circ z(\cdot))$ is extremal. Item ii) is proved.
            \end{itemize}
        \end{proof}

        \subsubsection{Structure $\gamma_+ \gamma_{c_1} \gamma_+$ and limit case $\gamma_+^{{c_1}}$}
        \label{sec:Shooting2}

        We introduce for $z \coloneqq (x,p) \in T^*M$ the true Hamiltonian
        \[
            H_{c_1}(z)  \coloneqq  H(z,u_{c_1}(x),\eta_{c_1}(z)) = H_0(z) + u_{c_1}(x)\, H_1(z) + \eta_{c_1}(z)\, c_1(x),
        \]
        where $c_1(x) = x_1 - 1$, and $u_{c_1}$ and $\eta_{c_1}$ are given by proposition \ref{prop:ControlEta}.
        Here we consider a trajectory with a structure of the form $\gamma_+ \gamma_{c_1} \gamma_+$, with non-empty arcs.
        We note $t_1 < t_2$ the junction times. To define the
        \emph{multiple shooting function} associated to this structure, we need the results about the transversality conditions and the level of the Hamiltonian
        from section \ref{sec:PMP}, and we use also lemmas \ref{lemma:m1}, \ref{lemma:junction1} and proposition \ref{prop:invariant1}.
        We assume that \ref{Hyp3} holds along $\gamma_{c_1}$, and that at $t_2$, there is no contact with the switching surface.
        Then, since $F_1 \cdot c_1 \equiv \cste_7 > 0$, then \ref{Hyp1} is also satisfied, and from lemma \ref{lemma:junction1}
        we have at $t_1$, $\nu_1 = 0$, and at $t_2$, $\nu_2 < 0$
        and the control is continuous at $t_2$. Proposition \ref{prop:invariant1} gives two conditions at time $t_1$: $c_1 \left( \pi_x \left(z_1\right) \right) = H_1 \left( z_1 \right) = 0$.
        The multiple shooting function $S_2(p_0,t_f,t_1,t_2,\nu_2,z_1,z_2)$ is then given by the following equations.
        \[
        \begin{array}[ht!]{lll}
            0           \displaystyle = c_1       \left( \pi_x    \left(  z_1       \right) \right),                                    &
            0           \displaystyle = H_1                       \left(  z_1       \right)        ,                                    &
            1           \displaystyle = u_{c_1}   \left( \pi_x    \left(  z_2       \right) \right),                                        \\[0.5em]
            (1,0,0)     \displaystyle = \pi                       \left(  \expmap{z_2^+}{(t_f-t_2)}{\vvec{H_+}}        \right)        , &
            1           \displaystyle = H_+                       \left(  \expmap{z_2^+}{(t_f-t_2)}{\vvec{H_+}}        \right)        , &   \\[0.5em]
            0           \displaystyle =                                   \expmap{z_0  }{t_1}      {\vvec{H_+}}       - z_1           , &
            0           \displaystyle =                                   \expmap{z_1  }{(t_2-t_1)}{\vvec{H_{c_1}}}   - z_2           , &
        \end{array}
        \]
        where $z_0  \coloneqq  (x_0,p_0)$, $x_0=(0,0,0)$ is fixed, $z_2^+  \coloneqq  z_2 - \nu_2\, c'_1 (\pi_x (z_2)) \frac{\partial}{\partial p}$.
        The last two equations are classical matching conditions: see \cite{StBu} for details about multiple shooting methods.
        The \emph{multiple shooting method} in the case of structure $\gamma_+ \gamma_{c_1} \gamma_+$ consists in finding a zero of
        the multiple shooting function $S_2$, \ie in solving
        \[
            S_2(p_0,t_f,t_1,t_2,\nu_2,z_1,z_2) = 0.
        \]
        A zero of the shooting function $S_2 = 0$ gives a BC-extremal of the form $\gamma_+ \gamma_{c_1} \gamma_+$
        which satisfies the necessary conditions of optimality of section \ref{ref:PMPConstraint}.

        \medbreak
        \paragraph{\textbf{Limit case} $\gamma_+^{{c_1}}$}
        In the limit case where $t_1 = t_2$, we have a contact with the boundary $C_1$ instead of a junction. We note $\gamma_+^{{c_1}}$
        an arc $\gamma_+$ with a contact point with $C_1$.

        \subsubsection{Structure $\gamma_+ \gamma_{c_1} \gamma_+ \gamma_- \gamma_+^{{c_3}}$ and limit case
        $\gamma_+ \gamma_{c_1} \gamma_+^{{H_1},{c_3}}$}
        We define the multiple shooting function $S_3$ associated to the structure $\gamma_+ \gamma_{c_1} \gamma_+ \gamma_- \gamma_+^{{c_3}}$.
        We assume that \ref{Hyp3} holds along $\gamma_{c_1}$, and that at $t_2$ (the exit-time of $\gamma_{c_1}$),
        there is no contact with the switching surface.
        The last bang arc has a contact point with the boundary $C_3$ so is labeled $\gamma_+^{{c_3}}$.
        The multiple shooting function $S_3(p_0,t_f,t_1,t_2,\nu_2,t_3,t_4,t_5,\nu_5,z_1,z_2,z_3,z_4,z_5)$ is defined by the following equations.
        \[
        \begin{array}[ht!]{lll}
            0           \displaystyle = c_1       \left( \pi_x    \left(  z_1       \right) \right),                                    &
            0           \displaystyle = H_1                       \left(  z_1       \right)        ,                                    &
            1           \displaystyle = u_{c_1}   \left( \pi_x    \left(  z_2       \right) \right),                                        \\[0.5em]
            0           \displaystyle = H_1                       \left(  z_3         \right),                                          &
            0           \displaystyle = H_1                       \left(  z_4         \right),                                          &   \\[0.5em]
            0           \displaystyle = {c_3}   \left( \pi_x    \left(  z_5         \right) \right),                                    &
            0           \displaystyle = (F_0 \cdot c_3)(\pi_x    \left(  z_5         \right) ),                                            &   \\[0.5em]
            (1,0,0)     \displaystyle = \pi                       \left(  \expmap{z_5^+}{(t_f-t_5)}{\vvec{H_+}}        \right)        , &
            1           \displaystyle = H_+                       \left(  \expmap{z_5^+}{(t_f-t_5)}{\vvec{H_+}}        \right)        , &   \\[0.5em]
            0           \displaystyle =                                   \expmap{z_0  }{t_1}      {\vvec{H_+}}       - z_1           , &
            0           \displaystyle =                                   \expmap{z_1  }{(t_2-t_1)}{\vvec{H_{c_1}}}   - z_2           , &
            0           \displaystyle =                                   \expmap{z_2^+}{(t_3-t_2)}{\vvec{H_+}}       - z_3           , \\[0.5em]
            0           \displaystyle =                                   \expmap{z_3  }{(t_4-t_3)}{\vvec{H_-}}       - z_4           , &
            0           \displaystyle =                                   \expmap{z_4  }{(t_5-t_4)}{\vvec{H_+}}       - z_5           , &
        \end{array}
        \]
        where $z_0  \coloneqq  (x_0,p_0)$, $x_0=(0,0,0)$ is fixed, $z_2^+  \coloneqq  z_2 - \nu_2\, c'_1 (\pi_x (z_2)) \frac{\partial}{\partial p}$,
        and $z_5^+  \coloneqq  z_5 - \nu_5\, c'_3 (\pi_x (z_5)) \frac{\partial}{\partial p}$.
        The \emph{multiple shooting method} in the case of structure $\gamma_+ \gamma_{c_1} \gamma_+ \gamma_- \gamma_+^{{c_3}}$
        consists in finding a zero of the multiple shooting function $S_3$, \ie in solving
        \[
            S_3(p_0,t_f,t_1,t_2,\nu_2,t_3,t_4,t_5,\nu_5,z_1,z_2,z_3,z_4,z_5) = 0,
        \]
        and we get a BC-extremal of the form $\gamma_+ \gamma_{c_1} \gamma_+ \gamma_- \gamma_+^{{c_3}}$.

        \medbreak
        \paragraph{\textbf{Limit case} $\gamma_+ \gamma_{c_1} \gamma_+^{{H_1},{c_3}}$}
        In the limit case where $t_3 = t_4$, we have a contact of order 2 with the switching surface instead of two consecutive contacts of order 1,
        see \ref{sec:foldCase}.
        We note $\gamma_+^{{H_1},{c_3}}$ an arc $\gamma_+$ with a contact point of order 2 with the switching surface defined by $H_1=0$
        followed by a contact point with $C_3$.
        In this case, we have one unknown $\tau$ instead of two ($t_3$ and $t_4$), with associated point $z_\tau$
        (instead of $z_3$ and $z_4$) but we may add as extra unknown the parameter $\vmax$; all others parameters from $w$ remain fixed. Now,
        we may replace $H_1(z_3) = H_1(z_4) = 0$ by $H_1(z_\tau) = H_{01}(z_\tau) = 0$.

        \subsubsection{Structure $\gamma_+ \gamma_{c_1} \gamma_+ \gamma_- \gamma_{c_3}$}
        \label{sec:Shooting4}

        We introduce for $z \coloneqq (x,p) \in T^*M$ the true Hamiltonian
        \[
            H_{c_3}(z)  \coloneqq  H(z,u_{c_3}(x),\eta_{c_3}(z)) = H_0(z) + u_{c_3}(x)\, H_1(z) + \eta_{c_3}(z)\, c_3(x),
        \]
        where $c_3(x) = x_3 - 1$, and $u_{c_3}$ and $\eta_{c_3}$ are given by proposition \ref{prop:ControlEta}.
        We define now the multiple shooting function $S_4$ associated to the strucutre $\gamma_+ \gamma_{c_1} \gamma_+ \gamma_- \gamma_{c_3}$,
        We assume that \ref{Hyp3} holds along $\gamma_{c_1}$ and $\gamma_{c_3}$ (\ref{Hyp1} is satisfied along $\gamma_{c_1}$ and $\gamma_{c_3}$),
        and that at $t_2$ (the exit-time of $\gamma_{c_1}$), there is no contact with the switching surface.
        We use propositions \ref{prop:invariant1} and \ref{prop:invariant2} and lemmas \ref{lemma:m1}, \ref{lemma:junction1}, \ref{lemma:m2} and \ref{lemma:junction2}
        with the results from section \ref{sec:PMP} to get the multiple shooting function $S_4$.
        One can recall that the transversality conditions $p_1(t_f) = p_3(t_f) = 0$ are equivalent to $\Phi(t_f) = \dot{\Phi}(t_f) = 0$,
        and that the control is constant along $\gamma_{c_3}$. Then according to proposition \ref{prop:invariant2}, we only need to check that
        $c_3$ and $F_0 \cdot c_3$ are zero at the entry-time of $\gamma_{c_3}$.
        The following equations describe the multiple shooting function $S_4(p_0,t_f,t_1,t_2,\nu_2,t_3,t_4,\nu_4,z_1,z_2,z_3,z_4)$.
        \[
        \begin{array}[ht!]{lll}
            0           \displaystyle = c_1       \left( \pi_x    \left(  z_1       \right) \right),                                    &
            0           \displaystyle = H_1                       \left(  z_1       \right)        ,                                    &
            1           \displaystyle = u_{c_1}   \left( \pi_x    \left(  z_2       \right) \right),                                        \\[0.5em]
            0           \displaystyle = H_1                       \left(  z_3         \right),                                          &   \\[0.5em]
            0           \displaystyle = {c_3}   \left( \pi_x    \left(  z_4         \right) \right),                                    &
            0           \displaystyle = (F_0 \cdot c_3)(\pi_x    \left(  z_4         \right) ),                                            &   \\[0.5em]
            (1,0,0)     \displaystyle = \pi                       \left(  \expmap{z_4^+}{(t_f-t_4)}{\vvec{H_+}}        \right)        , &
            1           \displaystyle = H_+                       \left(  \expmap{z_4^+}{(t_f-t_4)}{\vvec{H_+}}        \right)        , &   \\[0.5em]
            0           \displaystyle =                                   \expmap{z_0  }{t_1}      {\vvec{H_+}}       - z_1           , &
            0           \displaystyle =                                   \expmap{z_1  }{(t_2-t_1)}{\vvec{H_{c_1}}}   - z_2           , &
            0           \displaystyle =                                   \expmap{z_2^+}{(t_3-t_2)}{\vvec{H_+}}       - z_3           , \\[0.5em]
            0           \displaystyle =                                   \expmap{z_3  }{(t_4-t_3)}{\vvec{H_-}}       - z_4           , &
        \end{array}
        \]
        where $z_0  \coloneqq  (x_0,p_0)$, $x_0=(0,0,0)$ fixed, $z_2^+  \coloneqq  z_2 - \nu_2\, c'_1 (\pi_x (z_2)) \frac{\partial}{\partial p}$,
        and $z_4^+  \coloneqq  z_4 - \nu_4\, c'_3 (\pi_x (z_4)) \frac{\partial}{\partial p}$.
        The \emph{multiple shooting method} in the case of structure $\gamma_+ \gamma_{c_1} \gamma_+ \gamma_- \gamma_{c_3}$
        consists in finding a zero of the multiple shooting function $S_4$, \ie in solving
        \[
            S_4(p_0,t_f,t_1,t_2,\nu_2,t_3,t_4,\nu_4,z_1,z_2,z_3,z_4) = 0,
        \]
        and we get a BC-extremal of the form $\gamma_+ \gamma_{c_1} \gamma_+ \gamma_- \gamma_{c_3}$.

        \bigbreak
        \paragraph{\textbf{Limit case} $\gamma_+ \gamma_{c_1} \gamma_+ \gamma_- \gamma_{c_3}$ with $u_{c_3}(\cdot) \equiv +1$}
        This is a particular case when the boundary and regular arcs are identical. This phenomenon happens at the interface between
        $\gamma_+ \gamma_{c_1} \gamma_+ \gamma_- \gamma_+^{{c_3}}$ and $\gamma_+ \gamma_{c_1} \gamma_+ \gamma_- \gamma_{c_3}$ (with $u_{c_3} < 1$) trajectories,
        see section \ref{sec:timeMinimalSynthesis}.

        \bigbreak
        \paragraph{\textbf{Limit case} $\gamma_+ \gamma_{c_1} \gamma_+ \gamma_- \gamma_{c_3}$ with $t_2 = t_3$} In the limit case when $t_2 = t_3$,
        we have a trajectory of the form $\gamma_+ \gamma_{c_1} \gamma_- \gamma_{c_3}$.
        From $S_4$, it is straightforward to obtain the multiple shooting function $S_5$ associated to the structure $\gamma_+ \gamma_{c_1} \gamma_- \gamma_{c_3}$.

        \subsubsection{\hampath\ code}
        For any case presented in section \ref{sec:multipleShooting}, the user must provide the true Hamiltonians and the shooting function. The \hampath{} code
        supply automatically (by AD and integrating the Hamiltonian and variational systems) the exponential mappings, the Jacobian of the shooting function and
        the nonlinear solver.

    \subsection{Differential path following method}
    \label{sec:homotopy}

        The shooting method solves a single optimal control problem, \ie all the parameters in $w$ are fixed.
        To solve a one-parameter family of optimal control problems, \textit{e.g.} for different values of $w_{10} = \imax$,
        we use differential path following techniques with arclength parametrization (or homotopy method).
        Let $h \colon \R^N\times\R\rightarrow\R^N$, $h(y,\lambda)$, denote an homotopic function. For example, one can consider the homotopic function
        defined by $S_1$, with $y \coloneqq (p_0,t_f)$, $N = 4$ and $\lambda  \coloneqq  w_{10}$ ($\imax$ is here an independent variable).
        We are interested in solving $h = 0$. Under some assumptions, the solutions set forms a one-dimensional manifold.

        The classical difficulties about homotopic methods consist in assuring that a curve in $h^{-1}(\{0\})$ exists, is sufficiently smooth and will intersect a fixed
        target homotopic level in finite length. Suppose $h$ is continuously differentiable and that we know $y_0$ such that $h(y_0,\lambda_0)=0$ and
        \[
            \rank \left(\frac{\partial h}{\partial y}(y_0,\lambda_0)\right) = N.
        \]
        Suppose also that $0$ is a regular value of $h$.
        Then a continuously differentiable curve starting from $(y_0, \lambda_0)$ exists and is either diffeomorphic to a circle or the real line.
        The curves in $h^{-1}(\{0\})$ are disjoints, and we call each branch of $h^{-1}(\{0\})$ a path of zeros.

        Unlike well-known prediction-correction methods, see \cite{Allgower2003p1618}, the \hampath{} code implements an algorithm which
        merely follow the path of zeros by integrating the associated differential system with a high order Runge-Kutta scheme, without any correction.
        The key point is to compute efficiently the Jacobian of the homotopic function with the tools presented in section \ref{sec:hampathCode}.
        See \cite{M3AS,OMS} for more details about the algorithm.

        We group together in Table \ref{table:homotopies}, the different homotopies we need for the numerical results.
        \begin{table}[ht!]
            \centering
            \begin{tabular}{llll}
                \medhrule
                Shooting fun.                                                       & Homotopic fun.        &
                Structure                                                   & Homotopic par.        \\
                \bighrule
                $S_1(p_0,t_f)$ &                                                    $\Hom_1$                &
                $\gamma_+$                                                          &$\imax$       \\
                $S_2(p_0,t_f,t_1,t_2,\nu_2,z_1,z_2)$ &                              $\Hom_2^{(a)}$          &
                $\gamma_+ \gamma_{c_1} \gamma_+$                                    &$\imax$       \\
                $S_2(p_0,t_f,t_1,t_2,\nu_2,z_1,z_2)$ &                              $\Hom_2^{(b)}$          &
                $\gamma_+ \gamma_{c_1} \gamma_+$                                    &$\vmax$       \\
                $S_3(p_0,t_f,t_1,t_2,\nu_2,t_3,t_4,t_5,\nu_5,z_1,z_2,z_3,z_4,z_5)$  &$\Hom_3$               &
                $\gamma_+ \gamma_{c_1} \gamma_+ \gamma_- \gamma_+^{{c_3}}$   &$\vmax$       \\
                $S_4(p_0,t_f,t_1,t_2,\nu_2,t_3,t_4,\nu_4,z_1,z_2,z_3,z_4)$          &$\Hom_4$               &
                $\gamma_+ \gamma_{c_1} \gamma_+ \gamma_- \gamma_{c_3}$              &$\vmax$       \\
                $S_5(p_0,t_f,t_1,t_2,t_3,\nu_3,z_1,z_2,z_3)$ &                      $\Hom_5$                &
                $\gamma_+ \gamma_{c_1} \gamma_- \gamma_{c_3}$                       &$\vmax$       \\
                \medhrule
                \\
            \end{tabular}
            \caption{The homotopic function name with the associated shooting function, the associated strategy and the homotopic parameter.}
            \label{table:homotopies}
        \end{table}

\section{Numerical results}
    \label{sec:numericalResults}

        \begin{table}[ht!]
            \centering
            \begin{tabular}{*{13}{l}}
                \medhrule
                $L_m$               &
                $R_m$               &
                $K_m$               &
                $V_\mathrm{alim}$   &
                $r$                 &
                $K_r$               &
                $g$                 &
                $K_f$               &
                $M$                 &
                $\rho$              &
                $S$                 &
                $C_x$               &
                $R_\mathrm{bat}$    \\
                \bighrule
                0.05 &
                0.03 &
                0.27 &
                150.0&
                0.33 &
                10.0 &
                9.81 &
                0.03 &
                250.0&
                1.293&
                2.0  &
                0.4  &
                0.05 \\
                \medhrule
                \\
            \end{tabular}
            \caption{Names and values of the electric solar car parameters.}
            \label{table:parametersAndValues}
        \end{table}

    \subsection{Case study}
    \label{sec:caseStudy}

    In the following subsections, we present the numerical results we obtained for different instances of problem \eqref{ref:OCP_tf_min}.
    From the experimental point of view we are interested in parameter values which are real world data, and we choose to model an electrical solar car.
    These values are taken from \cite{Messine2015} and given in Table~\ref{table:parametersAndValues}.
    The different instances of problem \eqref{ref:OCP_tf_min} depend on the values of $\imax$ and $\vmax$. The parameter $\af$ will be fixed to $100$.
    In this setting, the state unconstrained problem is quite simple since the optimal trajectory is $\gamma_+$, see proposition \ref{prop:bangOptimal}.
    The complexity comes
    from the two state constraints of orders 1 and 2. Different phenomenons are expected depending on which constraint is active.
    Let consider first an homotopy on the bound of an order 1 state constraint, let's say $\imax$. For big enough values of $\imax$, the optimal trajectory
    has no contact point nor boundary arc. It is known that if we start to decrease the value of $\imax$ until we get a trajectory with a contact point,
    then for values of $\imax$ just smaller, the contact point is turned into a boundary arc of small length.
    This is quite different for an order 2 state constraint.
    We start again with a big value of $\vmax$ and then we decrease $\vmax$ until we have a trajectory with a contact point.
    Then, for just smaller values of $\vmax$, we still have trajectories
    with only a contact point and not a boundary arc. This can be seen in subsection \ref{sec:influenceBounds}. Some examples with state constraints
    of different orders are solved (analytically) in Bryson \textit{et al.} \cite{Bryson63} and Jacobson \textit{et al.} \cite{Jacobson1971}.

    \subsection{Procedure to study the influence of $\imax$ and $\vmax$ on BC-extremals}
    \label{sec:method}

    We start this section with a remark which is crucial for the numerical experiments of section \ref{sec:numericalResults}.
    \begin{rmrk}
        Actually, from any zero of any shooting function is associated a unique 4-tuple $(x(\cdot),p(\cdot),u(\cdot),\eta(\cdot))$ which is not
        necessarily a BC-extremal. The shooting method does not guarantee that the state $x(\cdot)$ satisfies the path constraints and
        that the times are well ordered. Indeed, solving for example $S_2=0$ could give a time $t_2$ smaller than the time $t_1$. Hence, when we 
        solve the shooting equations, we have to check a posteriori that the associated 4-tuple is a BC-extremal. To guarantee that the 4-tuple
        is a BC-extremal we check that $c_1(x(t)) \le 0$ and $c_3(x(t)) \le 0$ for all $t \in \intervalleff{t_0}{t_f}$, 
        and we check that $t_i \le t_{i+1}$, where $t_i$ is either the initial time or the final time or a switching time.
    \end{rmrk}

    \begin{dfntn}
        We say that the 4-tuple associated to a zero of a shooting function is \emph{admissible} if it is a BC-extremal and \emph{non admissible}
        if it is not a BC-extremal. We say that a structure is \emph{admissible} if it exists a zero of the associated shooting function for which the associated
        4-tuple is admissible, and we say that a structure is \emph{non admissible} in the other case.
    \end{dfntn}

    The procedure to study the influence of $\imax$ and $\vmax$ on the structure of the BC-extremals is the following:
    \begin{enumerate}
        \item The maximum principle combined with geometric analysis is applied first to reduce the set of possible types of extremals,
            then to compute the parameterization of each kind of possible
            extremals (\ie it gives the analytical expressions of the control and the Lagrange multiplier $\eta$), and finally to give junction conditions.
        \item For a given admissible structure, in order to apply indirect numerical methods, we need first to define the associated shooting function, and
            then the shooting function is computed, and its regularity and invertibility are checked.
        \item All the study is about how the structures evolve with respect to the parameters $\imax$ and $\vmax$.
            Hence, we fix starting values for $\imax$ and $\vmax$, and we fix an initial admissible structure.
            In our case, we fix $\imax$ and $\vmax$ big enough such that $\gamma_+$ is admissible (it is even the optimal structure,
            see proposition \ref{prop:bangOptimal}).
        \item Numerical simulations based on indirect methods are applied; starting from a BC-extremal given by the shooting methods,
            we use homotopy methods to deform the BC-extremal until we detect\footnote{With the \hampath\ code, we can check if the structure is admissible
            at each integration step of the homotopy method, and we can stop the homotopy if a change in the structure occurs.}
            that the associated 4-tuple becomes non admissible.
        \item When we detect along an homotopy that the 4-tuple is non admissible, we determine the new structure combining the theoretical results
            with the reason for the change in the structure and then we look for a new adimissible 4-tuple.
            Let note $\Lambda \coloneqq (\imax,\vmax)$ the value when the structure
            has to change. Before starting a new homotopy we valid the new structure by checking first
            that this new structure is admissible and then by checking that the limit cases for both structures (the old and the new) at $\Lambda$ give the same
            trajectory.
    \end{enumerate}

    \begin{rmrk}
        In this procedure, we do not check the global optimality of the BC-extremals but only the fact that the structures are valid, \ie admissible and
        have compatible limit cases.
    \end{rmrk}

    This procedure is exactly the method we use in section \ref{sec:influenceBounds}. To explain in detail the construction of the synthesis in section 
    \ref{sec:Numericalsynthesis}, we need the following definition about the degree of an admissible structure at a given point $\Lambda \coloneqq (\imax,\vmax)$.

    \begin{dfntn}
        We say that a structure $\gamma$ is of \emph{degree 1 at $\bar{\Lambda}$} if it exists a neighborhood
        $\mathcal{V} \subset \R^2$ of $\bar{\Lambda}$ such that for every $\Lambda \in \mathcal{V}$, $\gamma$ is admissible at $\Lambda$ (\ie when the values of the 
        parameters $\imax$ and $\vmax$ are given by $\Lambda$).
        We say that a structure $\gamma$ is of \emph{degree k at $\bar{\Lambda}$}, $k \ge 2$,
        if for every sufficiently small neighborhood $\mathcal{V}$ of $\bar{\Lambda}$, there exists $k$ structures $\gamma_i$,
        there exists $k$ sets $D_i$ of non empty interior
        and a set $I$ of empty interior such that $\{I, D_1, \cdots, D_{k} \}$ is a partition of $\mathcal{V}$ and such that
        for every $\Lambda_i \in D_i$, $i=1,\cdots,k$, $\gamma_i$ is of degree 1 at $\Lambda_i$, and the structures are compatible, \ie the limit
        cases at a given point of $\mathcal{V}$ give the same trajectory.
    \end{dfntn}

    \begin{rmrk}
        In the Figure \ref{fig:synthesis} presenting the synthesis, the structures of degree 2 are given by the blue lines,
        while the structures of degree greater than 3
        are represented by red points. The structures of degree 1 are contained in the domains of non-empty interiors.
    \end{rmrk}

    In section \ref{sec:influenceBounds}, to study the practical case when $\imax = 150$, we only perform homotopies on structures of degree 1 and
    we detect a change when we reach a structure\footnote{We can consider in this case that the probability to
    reach a structure of degree greater than 3 is zero.} of degree 2. To construct the synthesis with respect to the parameters $\imax$ and $\vmax$ we
    do the following:
    \begin{enumerate}
        \item We start an homotopy on a structure of degree 1 until we reach a structure of degree 2.
            We start a new homotopy with the new valid structure of degree 1 and continue the exploration.
        \item We compute the branch (a blue line on Figure \ref{fig:synthesis}) corresponding to a structure of degree 2 by homotopy until
            we reach a structure of degree greater than 3.
    \end{enumerate}

    We can see that the construction of the synthesis is heuristic and about exploration. Thus, the synthesis represented by the
    Figures \ref{fig:synthesis} and \ref{fig:synthesisZooms} is not necessarily complete. Moreover, this procedure does not guarantee the optimality
    of the BC-extremals but this synthesis is a first step toward the optimal synthesis.

    \subsection{Influence of the maximal current and the maximal velocity (with $\af = 100$)}
        \label{sec:influenceBounds}

        \begin{figure}[ht!]
            \centering
            \begin{tikzpicture}[scale=0.8]

                \coordinate (O) at (0,0);
                \def\xshift{18em}
                \def\yshift{16em}
                \def\epsx{5em}
                \def\epsy{3em}
                \def\ecartx{4em}
                \def\ecarty{3em}
                \draw[->,thin] (O) -- ([xshift=\xshift] O) node[below] {$\imax$};
                \draw[->,thin] (O) -- ([yshift=\yshift] O) node[left]  {$\vmax$};

                \coordinate (A0) at ([xshift=\xshift-\epsx, yshift=\yshift-\epsy] O);

                \coordinate (A1) at ([xshift=-\ecartx] A0); \draw (A0) node{$\circ$};
                \draw[very thick] (A0) -- (A1) node[midway,above]{$\Hom_1$} node{$\bullet$};

                \coordinate (A2a) at ([xshift=-\ecartx] A1);
                \draw[very thick] (A1) -- (A2a) node[midway,above]{$\Hom_2^{(a)}$};

                \coordinate (A2b) at ([yshift=-\ecarty] A2a);  \draw (A2a) node{$\circ$};
                \draw[very thick] (A2a) -- (A2b) node[midway,left]{$\Hom_2^{(b)}$} node{$\bullet$};

                \coordinate (A3) at ([yshift=-\ecarty] A2b);
                \draw[very thick] (A2b) -- (A3) node[midway,left]{$\Hom_3$} node{$\bullet$};

                \coordinate (A4) at ([yshift=-\ecarty] A3);
                \draw[very thick] (A3) -- (A4) node[midway,left]{$\Hom_4$} node{$\bullet$};

                \coordinate (A5) at ([yshift=-\ecarty] A4);
                \draw[very thick] (A4) -- (A5) node[midway,left]{$\Hom_5$} node{$\circ$};

                \draw[very thin, dashed]
                ([xshift=\xshift-\epsx-2*\ecartx,yshift=\yshift-\epsy] O) -- ([xshift=0,yshift=\yshift-\epsy] O)                        node[left]{$110$};
                \draw[very thin, dashed]
                ([xshift=\xshift-\epsx-2*\ecartx,yshift=\yshift-\epsy-\ecarty] O) -- ([xshift=0,yshift=\yshift-\epsy-\ecarty] O)        node[left]{$\vmax^{c_3}$};
                \draw[very thin, dashed]
                ([xshift=\xshift-\epsx-2*\ecartx,yshift=\yshift-\epsy-2*\ecarty] O) -- ([xshift=0,yshift=\yshift-\epsy-2*\ecarty] O)    node[left]{$\vmax^{\gamma_{c_3}}$};
                \draw[very thin, dashed]
                ([xshift=\xshift-\epsx-2*\ecartx,yshift=\yshift-\epsy-3*\ecarty] O) -- ([xshift=0,yshift=\yshift-\epsy-3*\ecarty] O)    node[left]{$\vmax^{+}$};
                \draw[very thin, dashed]
                ([xshift=\xshift-\epsx-2*\ecartx,yshift=\yshift-\epsy-4*\ecarty] O) -- ([xshift=0,yshift=\yshift-\epsy-4*\ecarty] O)    node[left]{$10$};
                \draw[very thin, dashed]
                ([xshift=\xshift-\epsx-2*\ecartx,yshift=\yshift-\epsy-4*\ecarty] O) -- ([xshift=\xshift-\epsx-2*\ecartx,yshift=0] O)    node[below]{$150$};
                \draw[very thin, dashed]
                ([xshift=\xshift-\epsx-1*\ecartx,yshift=\yshift-\epsy-0*\ecarty] O) -- ([xshift=\xshift-\epsx-1*\ecartx,yshift=0] O)    node[below]{$\imax^{c_1}$};
                \draw[very thin, dashed]
                ([xshift=\xshift-\epsx-0*\ecartx,yshift=\yshift-\epsy-0*\ecarty] O) -- ([xshift=\xshift-\epsx-0*\ecartx,yshift=0] O)    node[below]{$1100$};

            \end{tikzpicture}
            \caption{Schematic view of the method used to compute the sub-optimal synthesis for $\alpha_f = 100$, $\imax = 150$ and $\vmax \in \intervalleff{10}{110}$,
            starting from $(\imax,\vmax)=(1100,110)$. 
            The symbol $\bullet$ represents a change of structure while $\circ$ means a choice to start or stop an homotopy.
                The scale is not respected. See Table~\ref{table:homotopies} for the corresponding structures.}
            \label{fig:synthesis_imax_150}
        \end{figure}
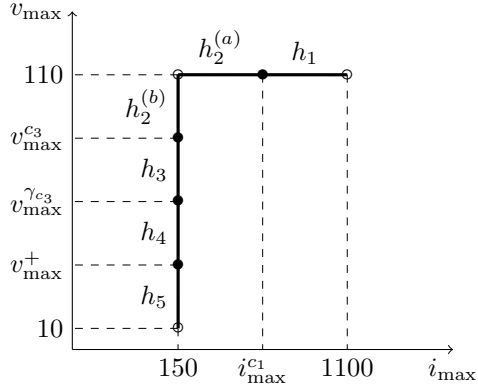

        From now on, the parameter $\af$ is fixed to $100$, so the car has to cover $100$m in minimum time.
        We are interested in the case $\imax = 150$ since it has a practical interest, see \cite{Messine2015}.
        In this case, the optimal trajectories may have boundary arcs.
        Thus, we start with values of $\imax$ and $\vmax$ big enough to deal with trajectories without boundary arcs.
        Then we first decrease the value of $\imax$ to $\imax = 150$ and finally we build a sub-optimal synthesis
        for $\imax = 150$ and $\vmax \in \intervalleff{10}{110}$.
        A schematic view of the method is presented in Figure~\ref{fig:synthesis_imax_150}, that we explain hereinafter.

        \begin{itemize}
            \item Let fix $(\imax,\vmax)=(1100,110)$. 
                By proposition \ref{prop:bangOptimal} the optimal trajectory is in this case $\gamma_+$.
            \item When $(\imax,\vmax)=(\imax^{c_1},110)$, the trajectory has now a contact point with the boundary $C_1$,
                thus the strategy is $\gamma_+^{c_1}$, $c_1(x) = x_1 - 1$.
                For $(\imax,\vmax)=(\imax^{c_1}-\varepsilon,110)$, $\varepsilon > 0$ small, the contact point has turned into a boundary arc since the constraint
                is of order 1.
                The structure becomes $\gamma_+ \gamma_{c_1} \gamma_+$, with $t_1<t_2$ the junction times with the boundary.
                According to lemma \ref{lemma:junction1}, $\nu_1 = 0$ and $\nu_2 < 0$, so the control is continuous at $t_2$. Moreover, the boundary
                control $u_{c_1}(\cdot)$ is strictly increasing since $\cste_4 + \cste_5 + \cste_6 > 0$,
                so it has to be discontinuous at time $t_1$ and continuous at time $t_2$ and \ref{Hyp3} is satisfied.
            \item When $(\imax,\vmax)=(150,\vmax^{c_3})$, the last arc has a contact point with
                the boundary $C_3$ ($c_3(x) = x_3 - 1$).
                For just smaller values of $\vmax$, \ie for
                $\vmax = \vmax^{c_3} - \varepsilon$, $\varepsilon > 0$ small, the structure has to change. Indeed, the last bang arc would cross
                the boundary, so this arc becomes a sequence $\gamma_+ \gamma_{-} \gamma_+$ (there is no singular arc, see proposition \ref{prop:bang-bang}).
                It is quite different from previous item, since here the constraint is of order 2; the trajectory for $\vmax = \vmax^{c_3} - \varepsilon$
                still has a contact point with $C_3$.
                As a consequence, the limit case for $\vmax = \vmax^{c_3}$ is not $\gamma_+ \gamma_{c_1} \gamma_+^{{c_3}}$ but
                $\gamma_+ \gamma_{c_1} \gamma_+^{{H_1,c_3}}$ and the last bang arc has a contact of order 2 with the switching surface.
            \item When $(\imax,\vmax)=(150,\vmax^{\gamma_{c_3}})$, the boundary control $u_{c_3}(\cdot)$ is admissible, \ie $u_{c_3}(\cdot)= 1$, which
                gives $\vmax^{\gamma_{c_3}} \approx 65.6042$. The structure for $\vmax = \vmax^{\gamma_{c_3}} - \varepsilon$ is
                $\gamma_+ \gamma_{c_1} \gamma_+ \gamma_- \gamma_{c_3}$.
            \item When $(\imax,\vmax)=(150,\vmax^{+})$ the second bang arc collapses,
                and the trajectories become $\gamma_+ \gamma_{c_1} \gamma_- \gamma_{c_3}$.
        \end{itemize}

        \subsubsection{Homotopies $\Hom_1$ and $\Hom_2^{(a)}$ and intermediate trajectory of the form $\gamma_+^{{c_1}}$}
            \label{sec:Hom1}

            We first choose $\imax$ and $\vmax$ big enough to get a trajectory with only interior arcs.
            In this case, the optimal trajectory is of the form $\gamma_+$.
            We fix $\imax = 1100$A and $\vmax = 110$km.h$^{-1}$ and we solve $S_1(y_1) = 0$, where $y_1 \coloneqq (p_0, t_f)$.
            The solution is $\ysol_1 \approx (0.3615, 6.4479, 0.2416, 5.6156)$\footnote{The value of $\ysol_1$ with all the digits gives a very accurate solution
            since $\norme{S_1(\ysol_1)} \approx 1.5e^{-15}$.}.
            The trajectory $x(\cdot)$ with its associated costate $p(\cdot)$ are portrayed in Figure~\ref{fig:gamma+}.
            Let
            $
                \imax^{c_1} = \imax \times \max_{t\in \intervalleff{0}{t_f}} x_1(t) \approx 1081.94$, $\imax = 1100,
            $
            denote the maximal current along the trajectory $x(\cdot)$. This maximal current does not depend on $\imax$ for the trajectories of the form
            $\gamma_+$. For $(\imax,\vmax) = (\imax^{c_1}, 110)$, then the structure is $\gamma_+^{c_1}$.
            We use the differential homotopy method presented in section \ref{sec:homotopy} to solve $\Hom_1(y_1,\imax) = 0$,
            for $\imax \in \intervalleff{\imax^{c_1}}{1100}$, starting from $\imax = 1100$.

            When $\imax < \imax^{c_1}$ the structure is $\gamma_+ \gamma_{c_1} \gamma_+$. 
            To initialize the shooting method and solve $S_2(y_2) = 0$, $y_2  \coloneqq  (p_0, t_f, t_1, t_2, \nu_2, z_1, z_2)$, we use the BC-extremal $\gamma_+$ with
            $\imax = \imax^{c_1}$ from the path of zeros of $\Hom_1$. Then we solve $\Hom_2^{(a)}(y_2,\imax) = 0$,
            for $\imax \in \intervalleff{150}{\imax^{c_1}}$, starting from $\imax = \imax^{c_1}$. Figure~\ref{fig:homotopy2} displays the initial,
            the junction and the final times and the jump $\nu_2$ along the path of zeros of $\Hom_2^{(a)}$.
            One may notice that the jump is not zero when $\imax = \imax^{c_1}$. As a consequence,
            the extremals at $\imax = \imax^{c_1}$ from the paths of zeros of $\Hom_1$ and $\Hom_2^{(a)}$ are not equal. The state trajectories are the same but not
            the costate trajectories. However, the differences between the costates are clearly shown in Figure~\ref{fig:limiteCaseHom1Hom2}. As a matter of fact,
            from the solution $(\ysol_1,\imax^{c_1})$ of $\Hom_1(y_1,\imax^{c_1}) = 0$, it is possible to determine completely the solution $(\ysol_2,\imax^{c_1})$
            of $\Hom_2^{(a)}(y_2,\imax^{c_1}) = 0$, and we can easily initialize the homotopy $\Hom_2^{(a)}$.

            \def\sizeFig{0.3}
            \begin{figure}[ht!]
                \centering
                \includegraphics[width=\sizeFig\textwidth]{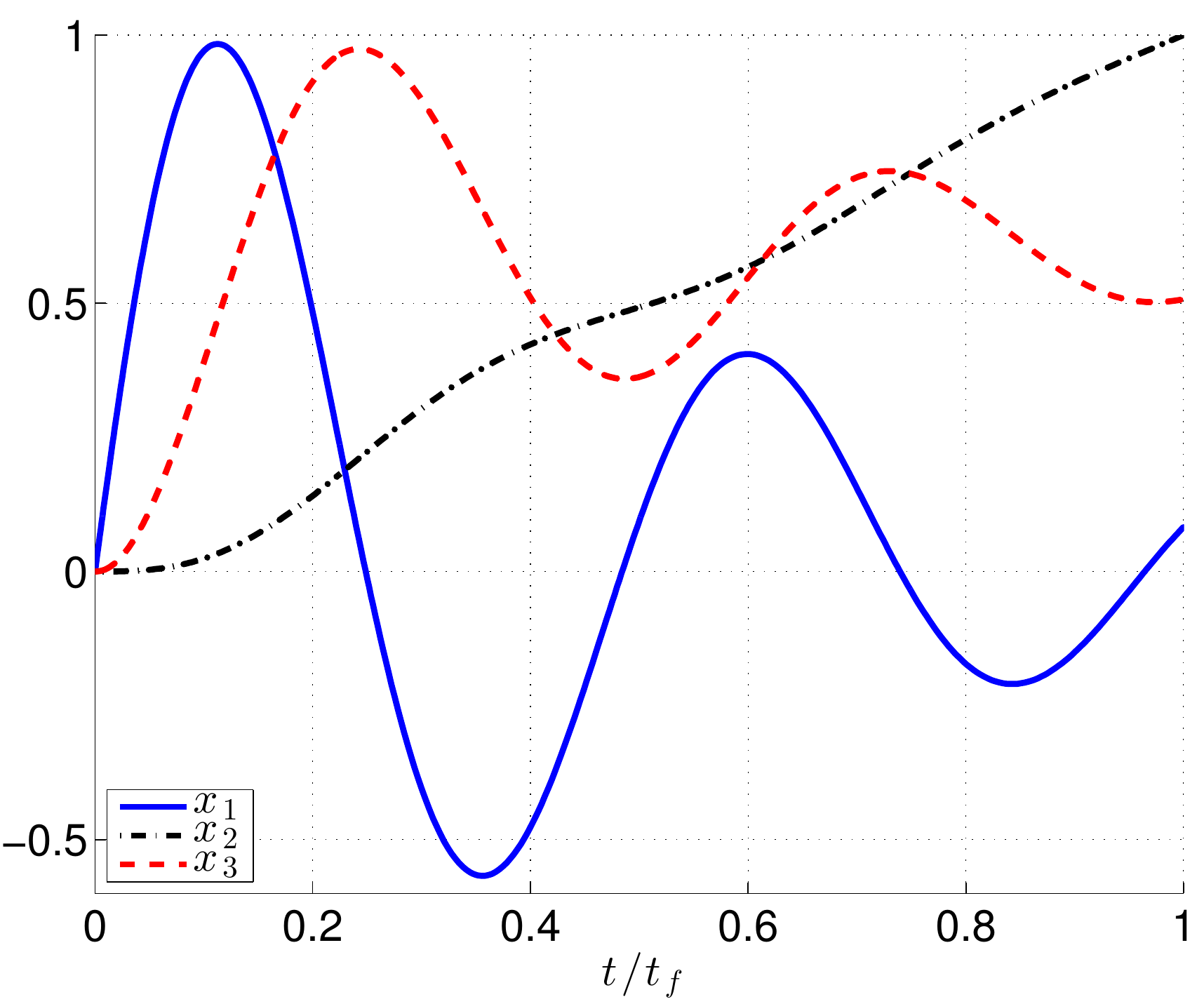}%
                \hspace{3em}
                \includegraphics[width=\sizeFig\textwidth]{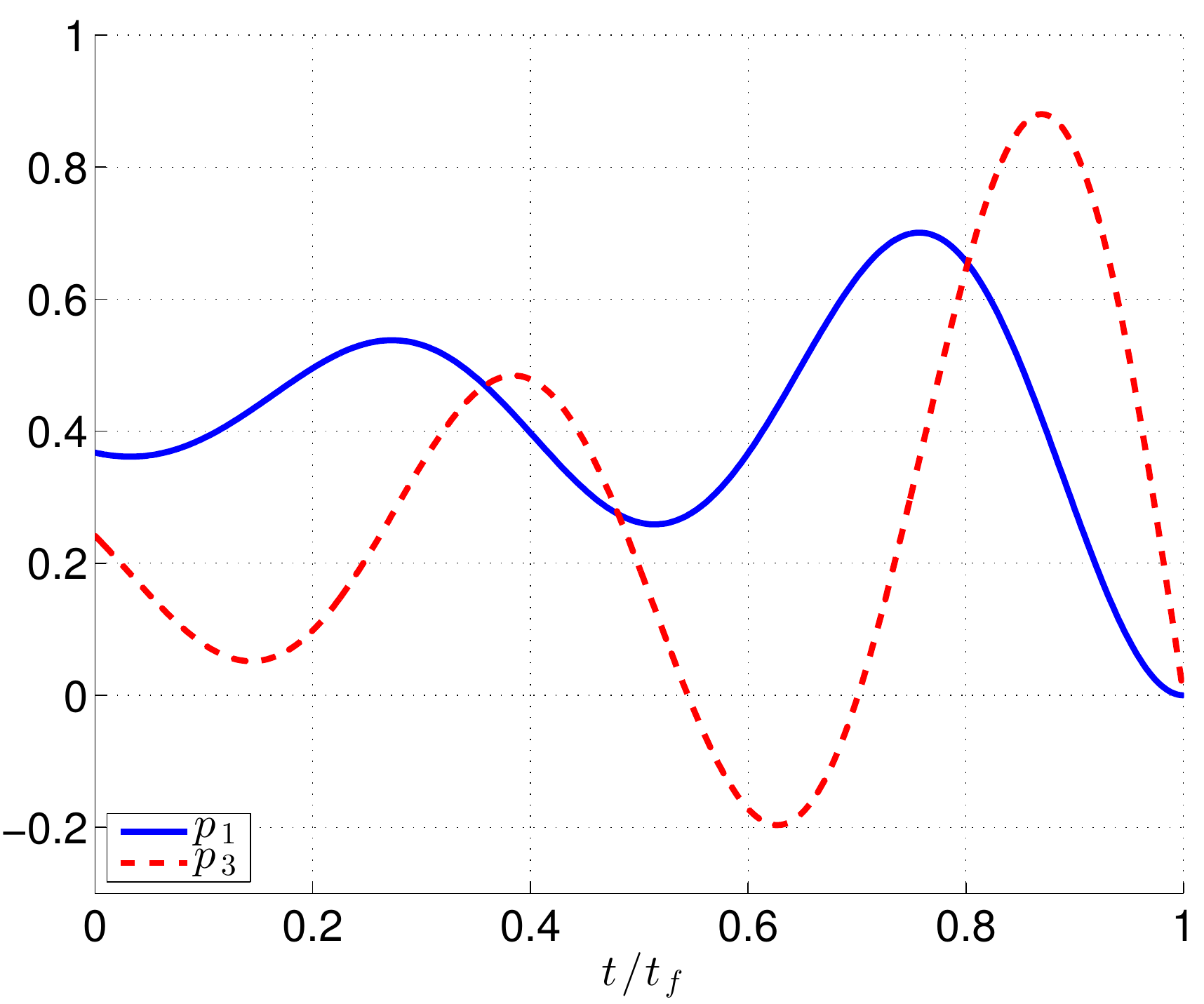}%
                \caption{\textbf{Trajectory $\gamma_+$}. State and costate for $(\imax,\vmax) = (1100, 110)$. The component $p_2 \approx 6.4479$ is constant.
                The terminal constraint $x_2(t_f) = 1$ is satisfied.}
                \label{fig:gamma+}
            \end{figure}

            \def\sizeFig{0.3}
            \begin{figure}[ht!]
                \centering
                \def\x{517}
                \def\y{407}
                \begin{tikzgraphics}{\sizeFig\textwidth}{\x}{\y}{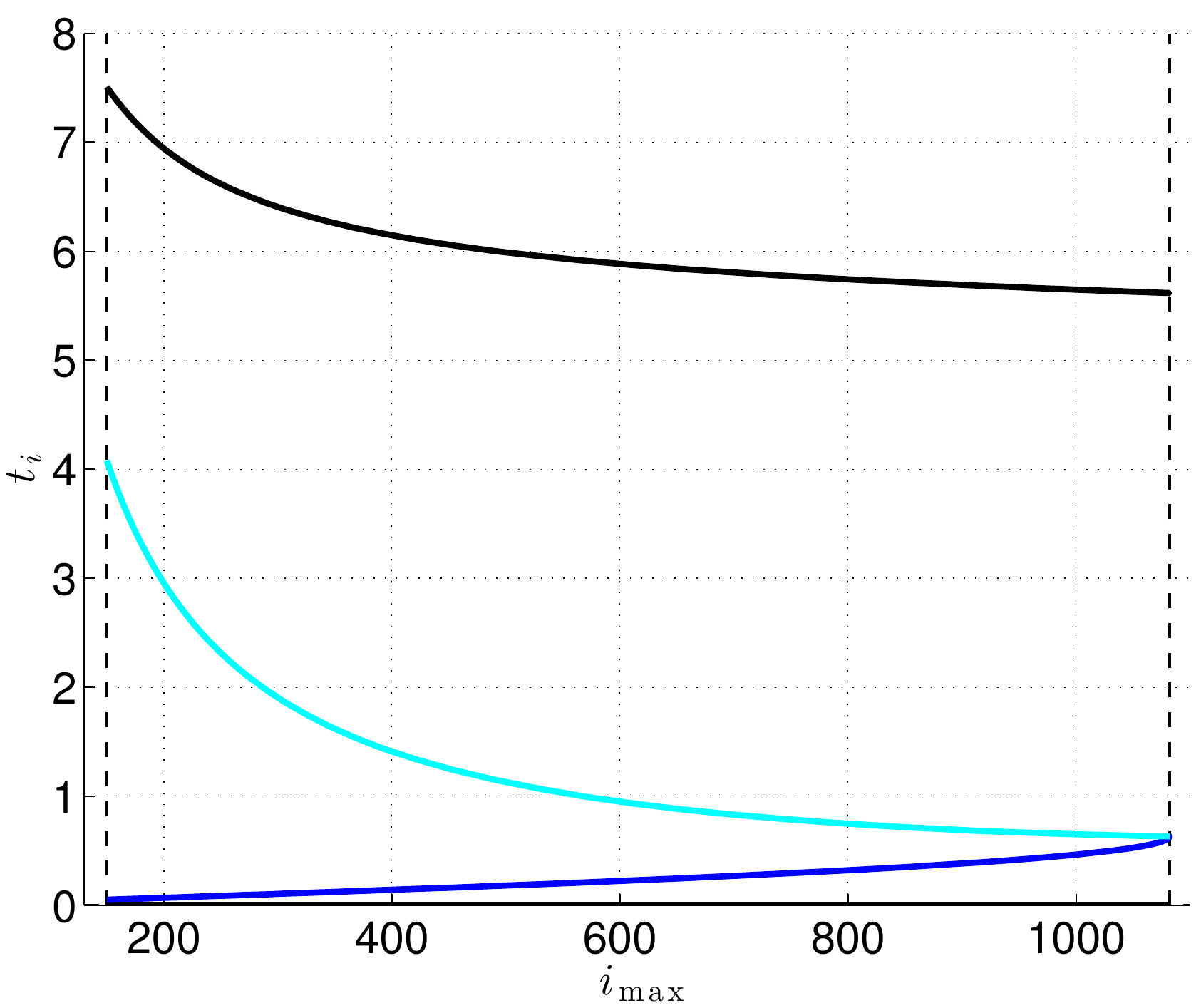}
                    \pxcoordinate{0.09*\x}{-0.05*\y}{A}; \draw (A) node {$150$};
                    \pxcoordinate{0.99*\x}{-0.05*\y}{A}; \draw (A) node {$\imax^{c_1}$};
                \end{tikzgraphics}
                \hspace{0.5em}
                \def\x{485}
                \def\y{407}
                \begin{tikzgraphics}{\sizeFig\textwidth}{\x}{\y}{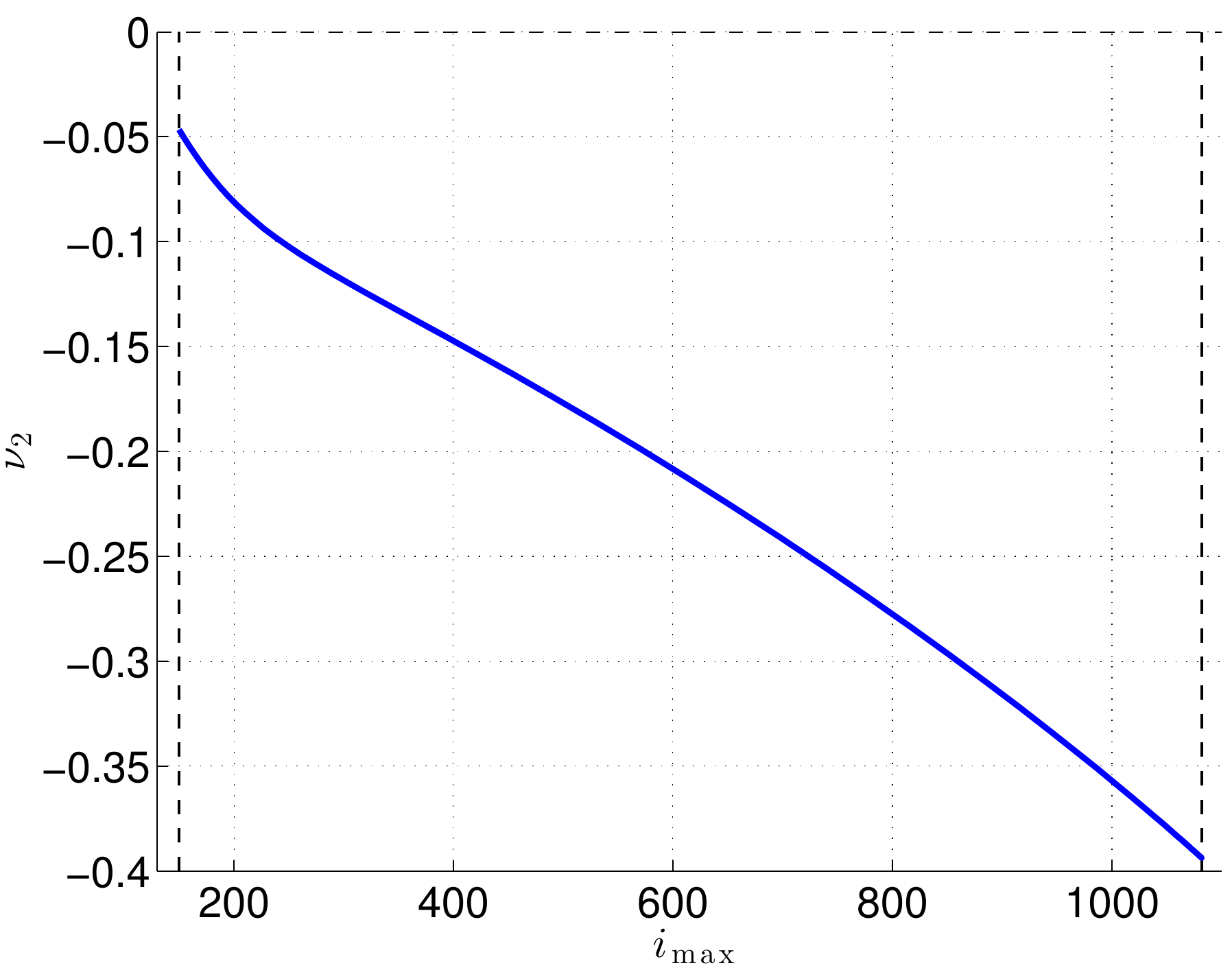}
                    \pxcoordinate{0.14*\x}{-0.01*\y}{A}; \draw (A) node {$150$};
                    \pxcoordinate{0.99*\x}{-0.01*\y}{A}; \draw (A) node {$\imax^{c_1}$};
                \end{tikzgraphics}
                \caption{\textbf{Homotopy $\Hom_2^{(a)}$}. (Left) The initial, the junction and the final times along the path of zeros of $\Hom_2^{(a)}$.
                The junction times collapse at $\imax = \imax^{c_1}$. The length of the boundary arc increases as $\imax$ decreases.
                (Right) The jump $\nu_2 \le 0$ along the path.}
                \label{fig:homotopy2}
            \end{figure}

            \def\sizeFig{0.3}
            \begin{figure}[ht!]
                \centering
                \def\x{495}
                \def\y{386}
                \begin{tikzgraphics}{\sizeFig\textwidth}{\x}{\y}{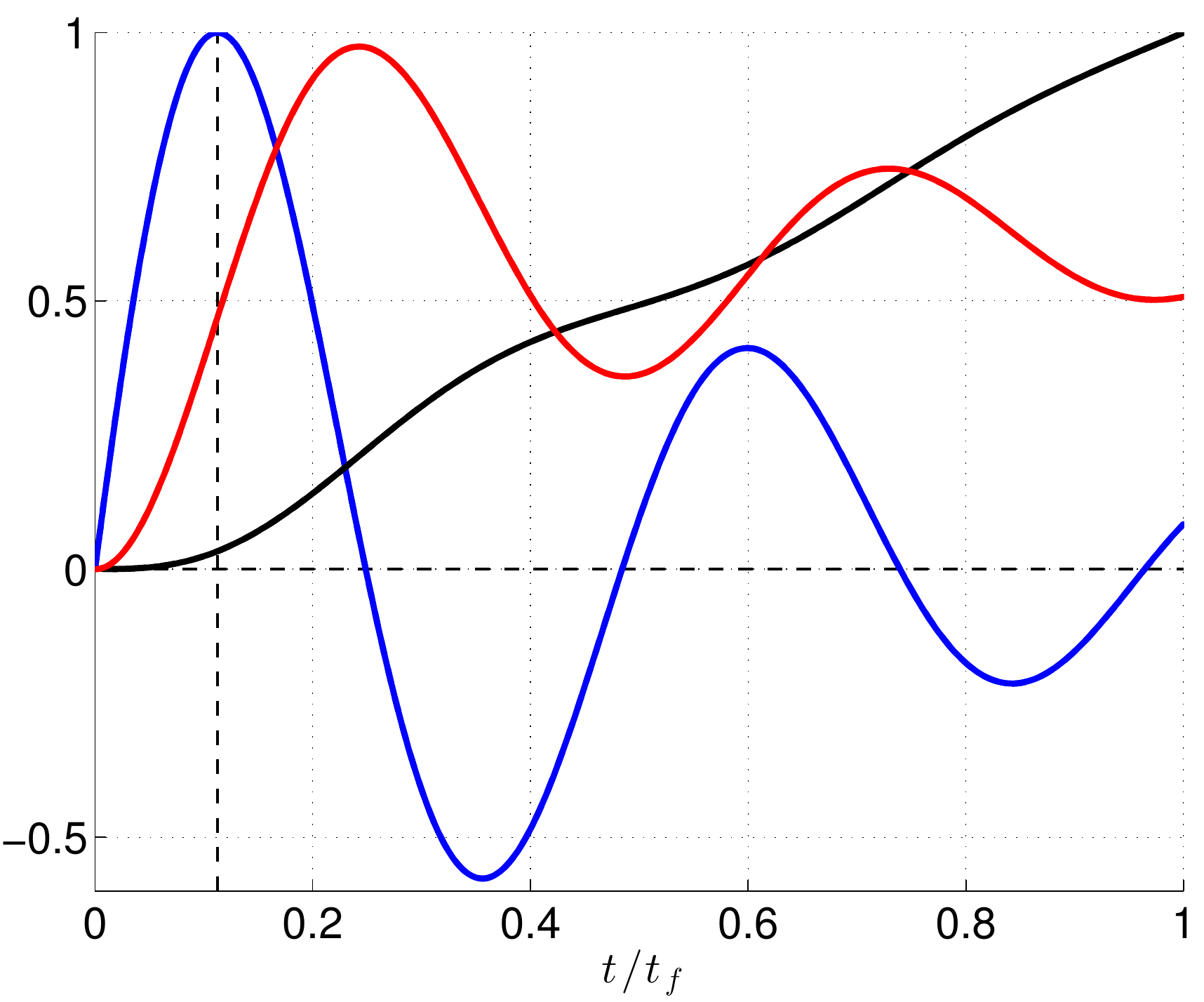}
                    \pxcoordinate{0.18*\x}{0.985*\y}{A};
                    \draw (A) node {$\tau$};
                \end{tikzgraphics}
%
                \begin{tikzgraphics}{\sizeFig\textwidth}{\x}{\y}{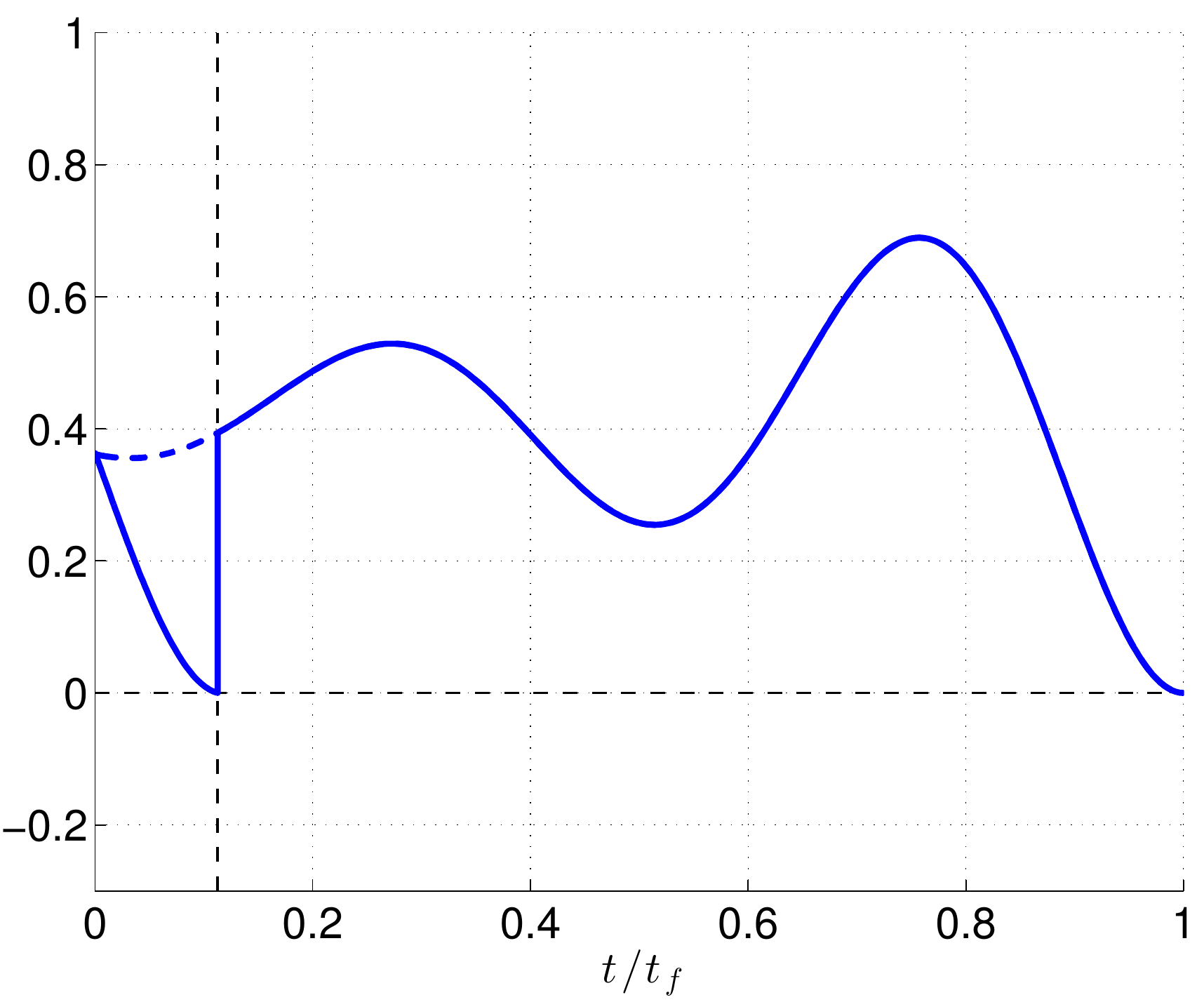}
                    \pxcoordinate{0.18*\x}{0.985*\y}{A};
                    \draw (A) node {$\tau$};
                \end{tikzgraphics}
                \begin{tikzgraphics}{\sizeFig\textwidth}{\x}{\y}{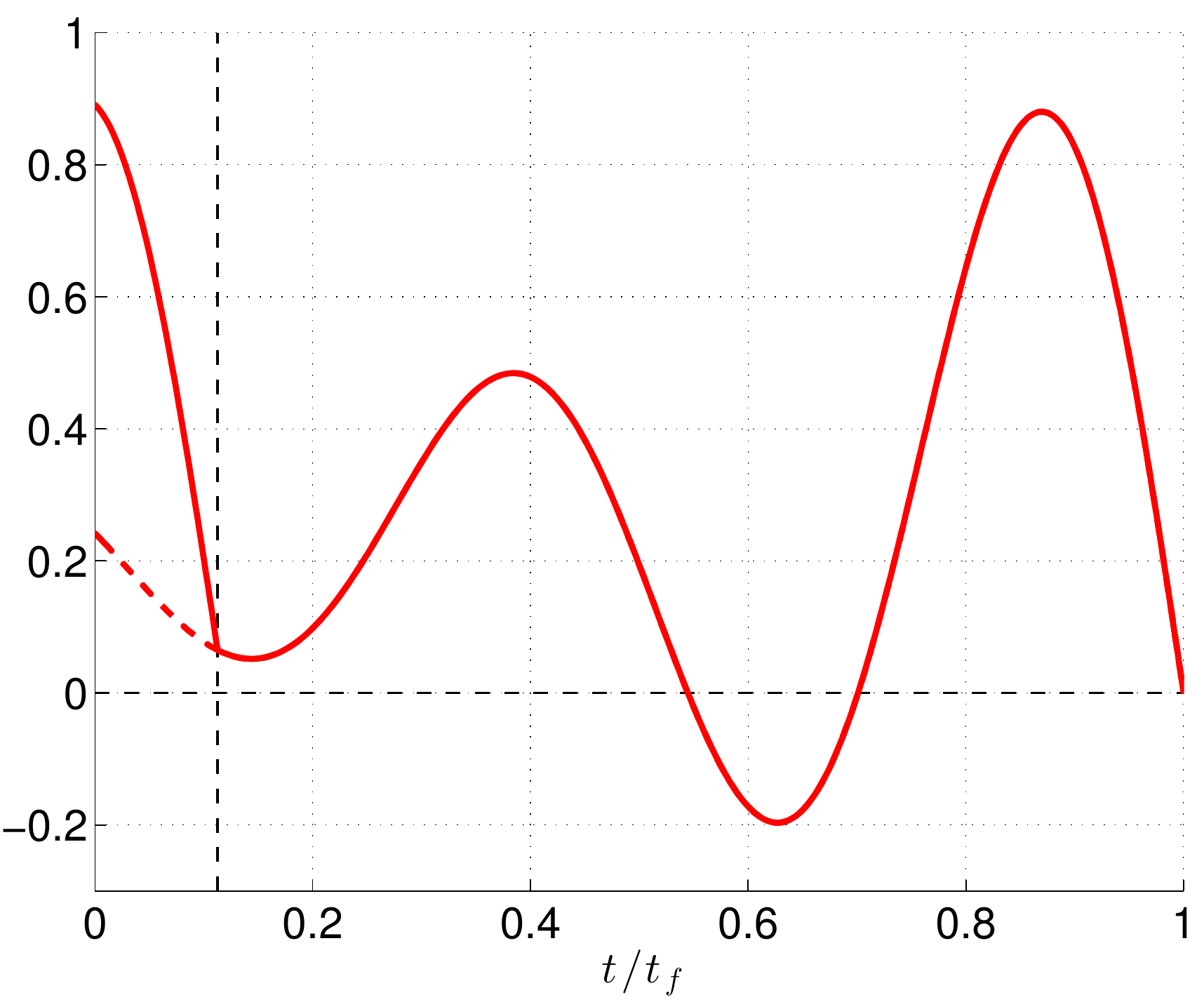}
                    \pxcoordinate{0.18*\x}{0.985*\y}{A};
                    \draw (A) node {$\tau$};
                \end{tikzgraphics}
                \caption{\textbf{Trajectories $\gamma_+^{c_1}$}.
                The BC-extremal of the form $\gamma_+^{c_1}$ at $\imax = \imax^{c_1}$ from $\Hom_1$ is represented by dashed lines
                while the one from $\Hom_2^{(a)}$ is represented by solid lines.
                The contact with the boundary is at time $\tau = t_1 = t_2$, where $t_1$, $t_2$ are the junction times
                of trajectories $\gamma_+\gamma_{c_1}\gamma_+$.
                (Left) The state trajectories (compare with Figure~\ref{fig:gamma+}).
                (Middle) $p_1(\cdot)$.
                (Right)  $p_3(\cdot)$.
                The trajectories are identical while the adjoint vectors are similar only after $\tau$.
                }
                \label{fig:limiteCaseHom1Hom2}
            \end{figure}

        \subsubsection{Homotopies $\Hom_2^{(b)}$ and $\Hom_3$ and intermediate trajectory of the form $\gamma_+ \gamma_{c_1} \gamma_+^{{H_1},{c_3}}$}
        \label{sec:homotopieC3}

            Up to this point, we have a trajectory $\gamma_+ \gamma_{c_1} \gamma_+$ for $(\imax,\vmax) = (150, 110)$. Now we start to decrease the value
            of $\vmax$ to obtain all the trajectories for $\imax = 150$. We start by computing the path of zeros of $\Hom_2^{(b)}$ until we reach
            $\vmax = \vmax^{c_3}$, for which the associated trajectory has a contact point with the boundary $C_3$.
            Let $x(\cdot)$ denote the trajectory of the form $\gamma_+ \gamma_{c_1} \gamma_+$ for $(\imax,\vmax) = (150, 110)$. Then the maximal speed
            is given by
            $
                \vmax^{c_3} = \vmax \times \max_{t\in \intervalleff{0}{t_f}} x_3(t) \approx 70.3716$, $\vmax = 110,
            $
            and it does not depend on $\vmax$ for the trajectories of the form $\gamma_+ \gamma_{c_1} \gamma_+$, for which $\imax = 150$ is fixed.
            For $(\imax,\vmax) = (150, \vmax^{c_3})$, then the structure is $\gamma_+ \gamma_{c_1} \gamma_+^{H_1,c_3}$.
            We solve $\Hom_2^{(b)}(y_2,\vmax)$, for $\vmax \in \intervalleff{\vmax^{c_3}}{110}$, starting from $\vmax = 110$.
            Along the path of zeros of $\Hom_2^{(b)}$ the initial, the junction and the final times, and the jump $\nu_2$ are constant.
            When $\vmax < \vmax^{c_3}$ the structure is $\gamma_+ \gamma_{c_1} \gamma_+ \gamma_- \gamma_+^{c_3}$.
            We fix now $\vmax = \vmax^{c_3}$, and $\imax = 150$.
            To initialize the shooting method and solve $S_3(y_3) = 0$, $y_3  \coloneqq  (p_0,t_f,t_1,t_2,\nu_2,t_3,t_4,t_5,\nu_5,z_1,z_2,z_3,z_4,z_5),$
            we use the BC-extremal $\gamma_+ \gamma_{c_1} \gamma_+$ at $\vmax = \vmax^{c_3}$ from the path of zeros of $\Hom_2^{(b)}$.
            Then we solve $\Hom_3(y_3,\vmax) = 0$ for $\vmax \in \intervalleff{\vmax^{\gamma_{c_3}}}{\vmax^{c_3}}$, starting from $\vmax = \vmax^{c_3}$ and with
            $\vmax^{\gamma_{c_3}} \approx 65.6042$.
            Figure~\ref{fig:homotopy3_times} shows the initial, the switching, the junction and the final times along the path of zeros of $\Hom_3$,
            while the jumps $\nu_2$ and $\nu_5$ are portrayed in Figure~\ref{fig:homotopy3_jumps}. Figure~\ref{fig:limiteCaseHom2Hom3} gives details about
            the limit case when $\vmax = \vmax^{c_3}$.

            \def\sizeFig{0.3}
            \begin{figure}[ht!]
                \def\posg{0.08}
                \def\posc{0.96}
                \centering
                \def\x{489}
                \def\y{403}
                \begin{tikzgraphics}{\sizeFig\textwidth}{\x}{\y}{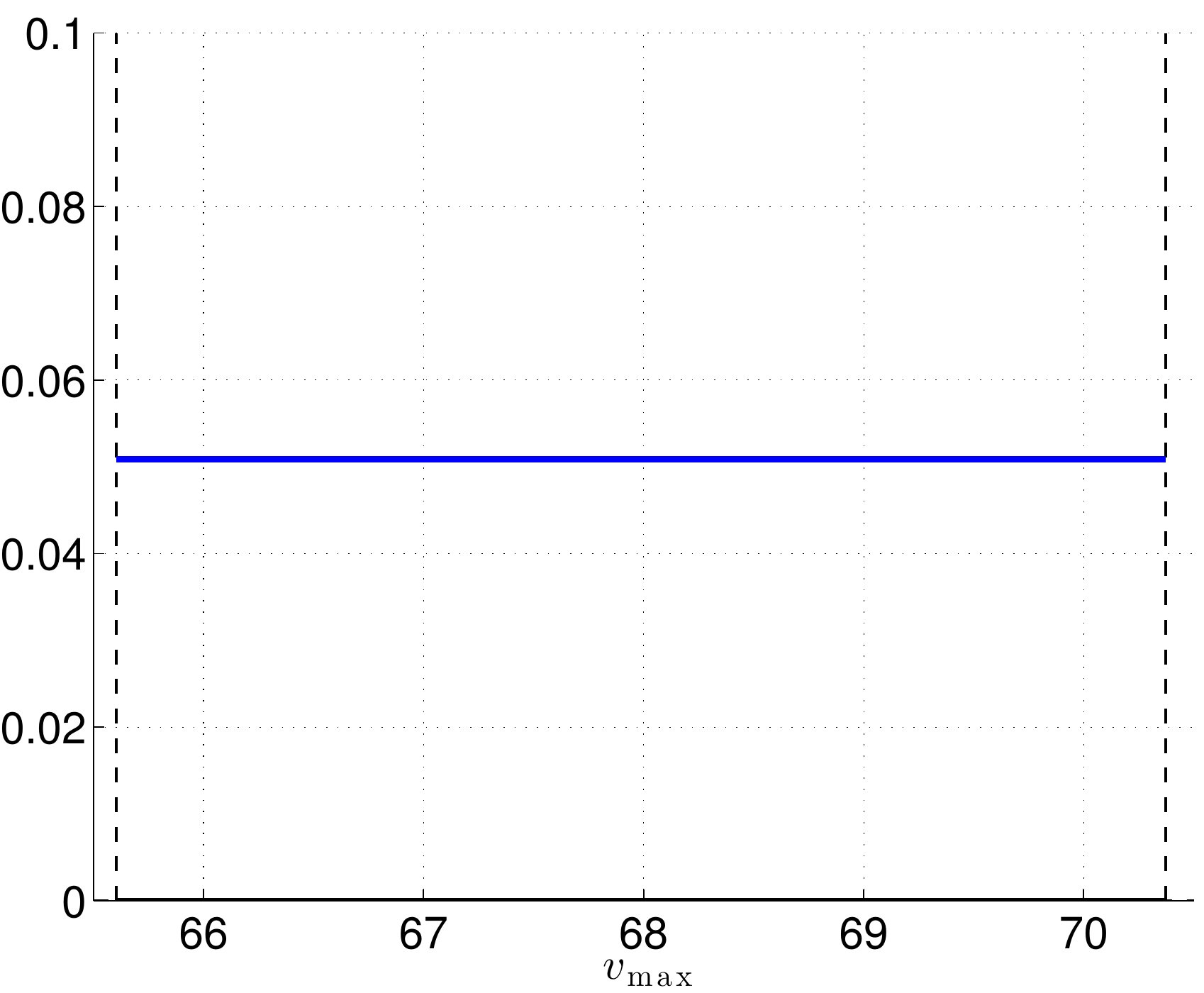}
                    \pxcoordinate{\posg*\x}{-0.05*\y}{A}; \draw (A) node {$\vmax^{\gamma_{c_3}}$};
                    \pxcoordinate{\posc*\x}{-0.05*\y}{A}; \draw (A) node {$\vmax^{c_3}$};
                \end{tikzgraphics}
                \def\x{479}
                \def\y{403}
                \begin{tikzgraphics}{\sizeFig\textwidth}{\x}{\y}{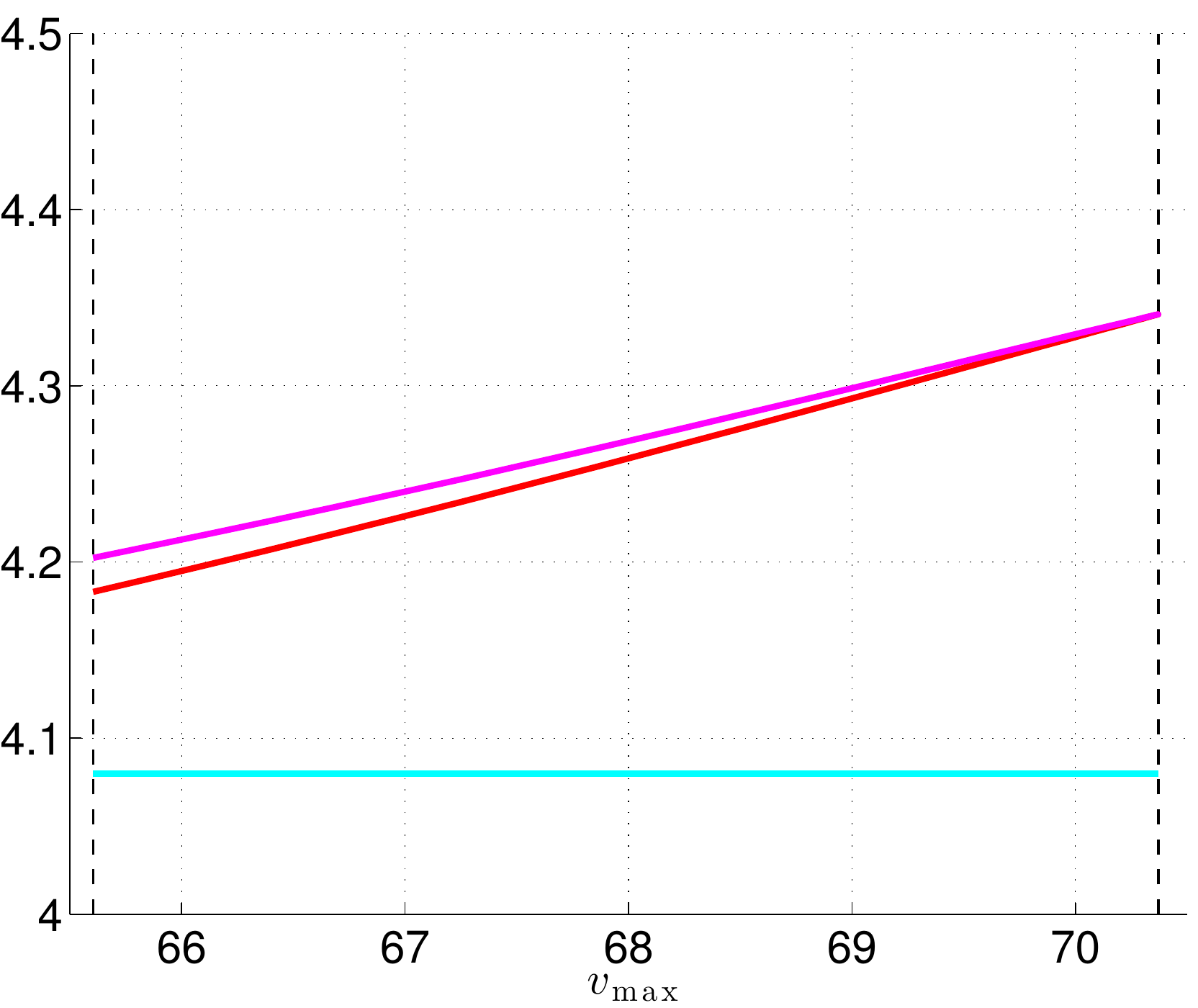}
                    \pxcoordinate{\posg*\x}{-0.05*\y}{A}; \draw (A) node {$\vmax^{\gamma_{c_3}}$};
                    \pxcoordinate{\posc*\x}{-0.05*\y}{A}; \draw (A) node {$\vmax^{c_3}$};
                \end{tikzgraphics}
                \def\x{479}
                \def\y{403}
                \begin{tikzgraphics}{\sizeFig\textwidth}{\x}{\y}{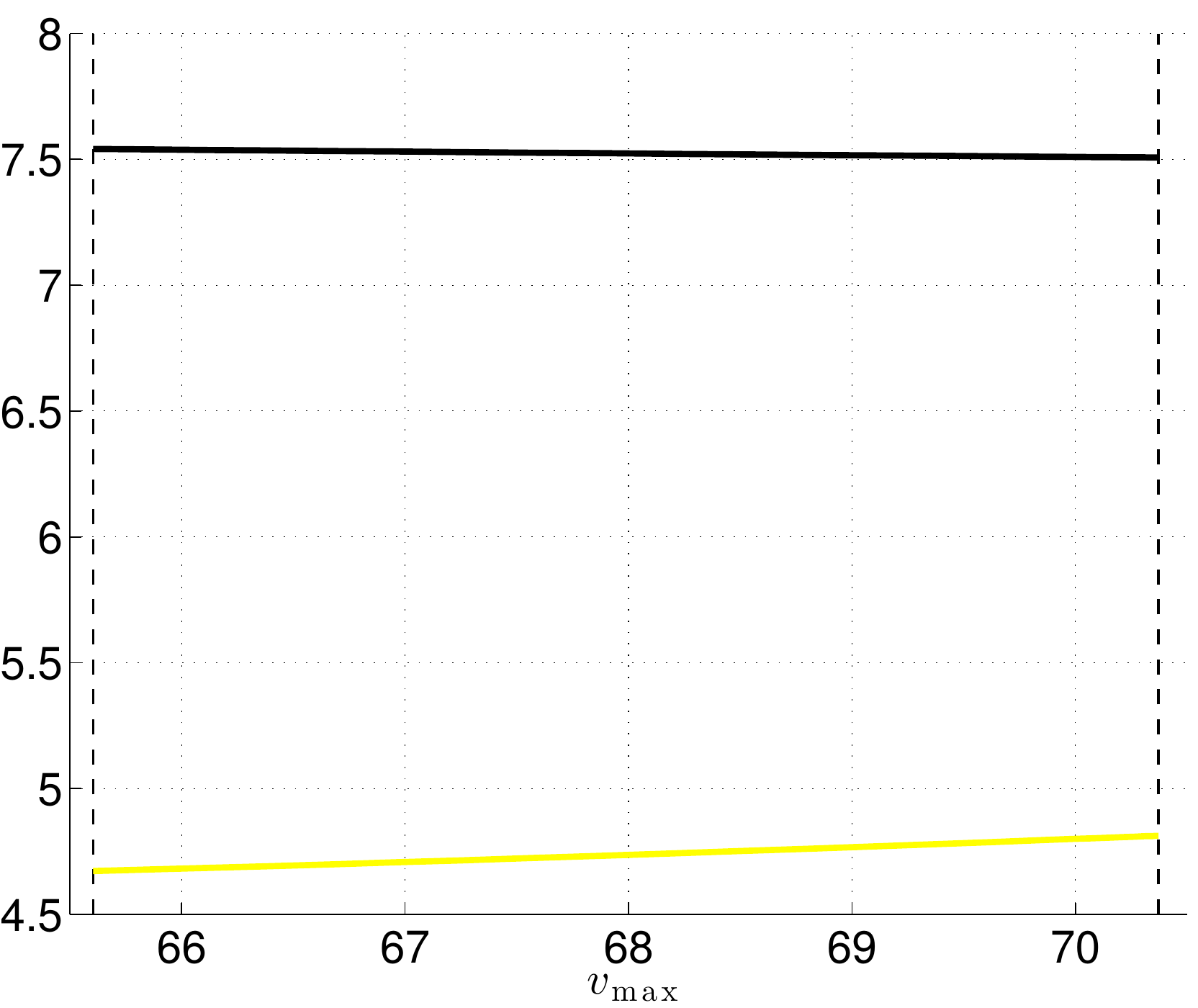}
                    \pxcoordinate{\posg*\x}{-0.05*\y}{A}; \draw (A) node {$\vmax^{\gamma_{c_3}}$};
                    \pxcoordinate{\posc*\x}{-0.05*\y}{A}; \draw (A) node {$\vmax^{c_3}$};
                \end{tikzgraphics}
                \caption{\textbf{Homotopy $\Hom_3$}. The initial, the switching, the junction and the final times along the path of zeros of $\Hom_3$.
                (Left) The times $t_0 \equiv 0$ and $t_1$.
                (Middle) The times $t_2$, $t_3$ and $t_4$.
                (Right) The times $t_5$ and $t_f$.
                The length $t_2-t_1$ of $\gamma_{c_1}$ is constant while the length of $\gamma_-$ vanishes when $\vmax = \vmax^{c_3}$, \ie when the trajectory is
                of the form
                $\gamma_+ \gamma_{c_1} \gamma_+^{{H_1},{c_3}}$.
                }
                \label{fig:homotopy3_times}
            \end{figure}
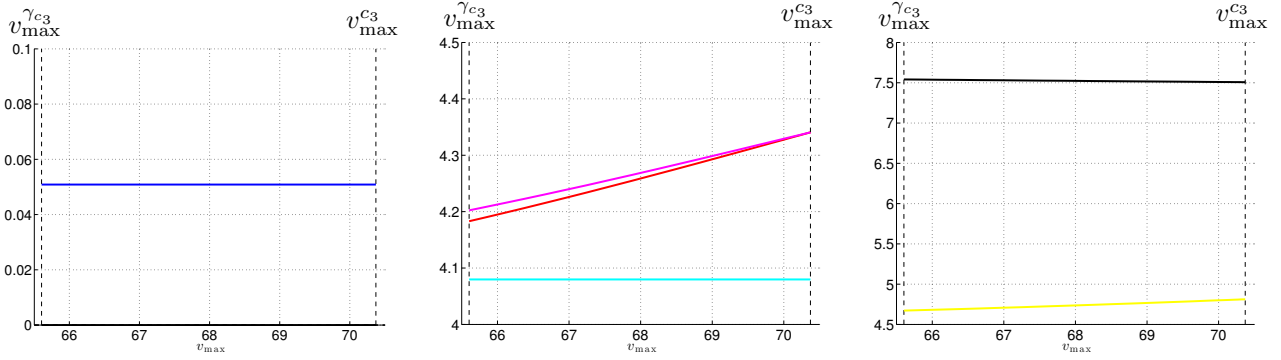

            \vspace{-0.5em}
            \def\sizeFig{0.3}
            \begin{figure}[ht!]
                \def\posg{0.15}
                \def\posc{0.96}
                \centering
                \def\x{526}
                \def\y{403}
                \begin{tikzgraphics}{\sizeFig\textwidth}{\x}{\y}{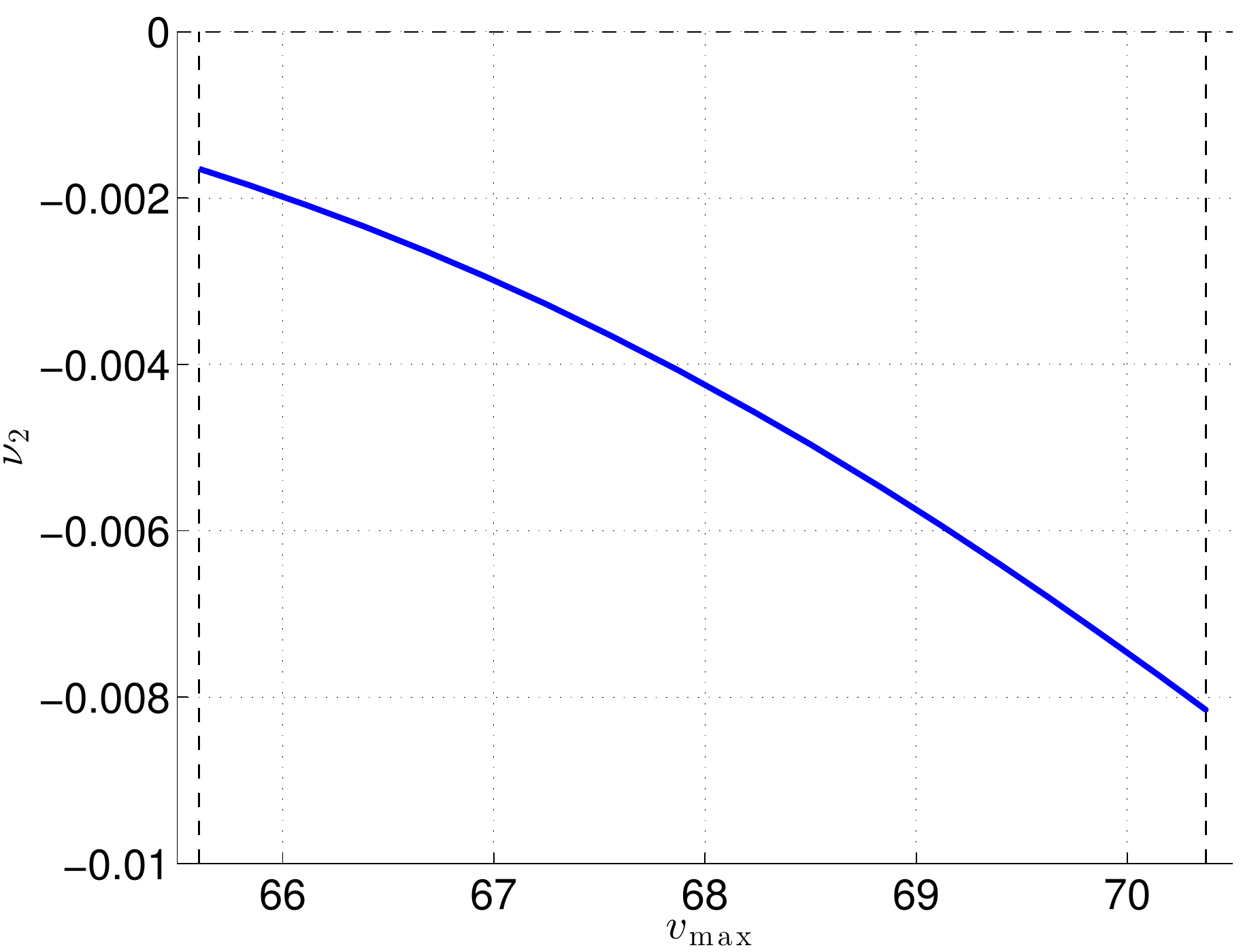}
                    \pxcoordinate{\posg*\x}{-0.05*\y}{A}; \draw (A) node {$\vmax^{\gamma_{c_3}}$};
                    \pxcoordinate{\posc*\x}{-0.05*\y}{A}; \draw (A) node {$\vmax^{c_3}$};
                \end{tikzgraphics}
                \hspace{3em}
                \def\x{517}
                \def\y{403}
                \begin{tikzgraphics}{\sizeFig\textwidth}{\x}{\y}{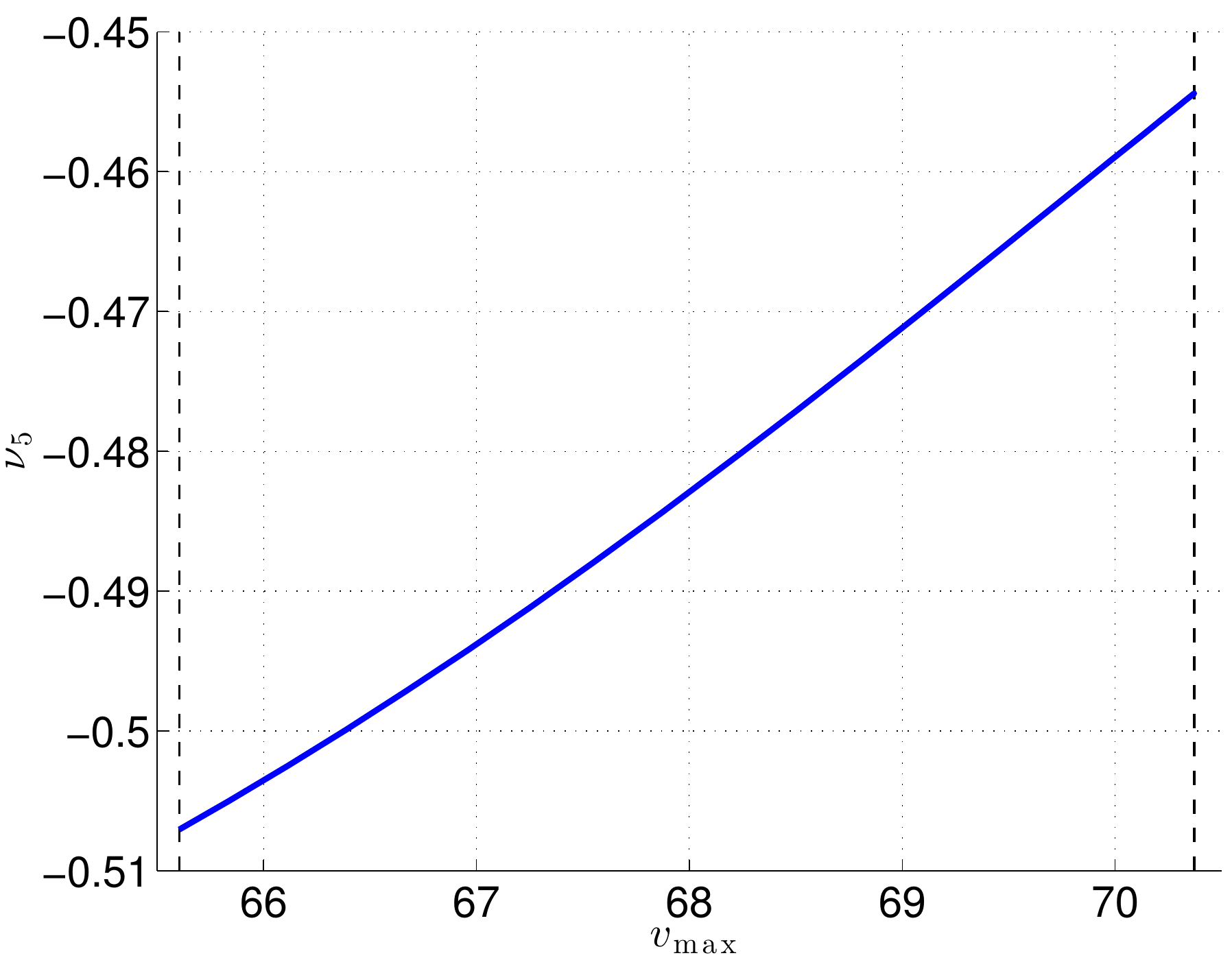}
                    \pxcoordinate{0.14*\x}{-0.04*\y}{A}; \draw (A) node {$\vmax^{\gamma_{c_3}}$};
                    \pxcoordinate{\posc*\x}{-0.04*\y}{A}; \draw (A) node {$\vmax^{c_3}$};
                \end{tikzgraphics}
                \caption{\textbf{Homotopy $\Hom_3$}. The jumps $\nu_2 \le 0$ (left) and $\nu_5 \le 0$ (right) along the path of zeros of $\Hom_3$.}
                \label{fig:homotopy3_jumps}
            \end{figure}

            \def\sizeFig{0.26}
            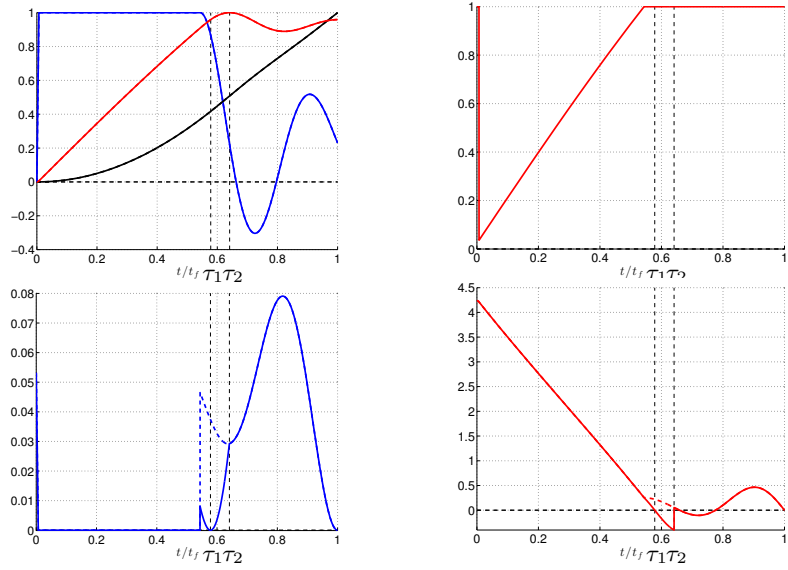
\begin{figure}[ht!]
                \centering
                \def\x{494}
                \def\y{411}
                \begin{tikzgraphics}{\sizeFig\textwidth}{\x}{\y}{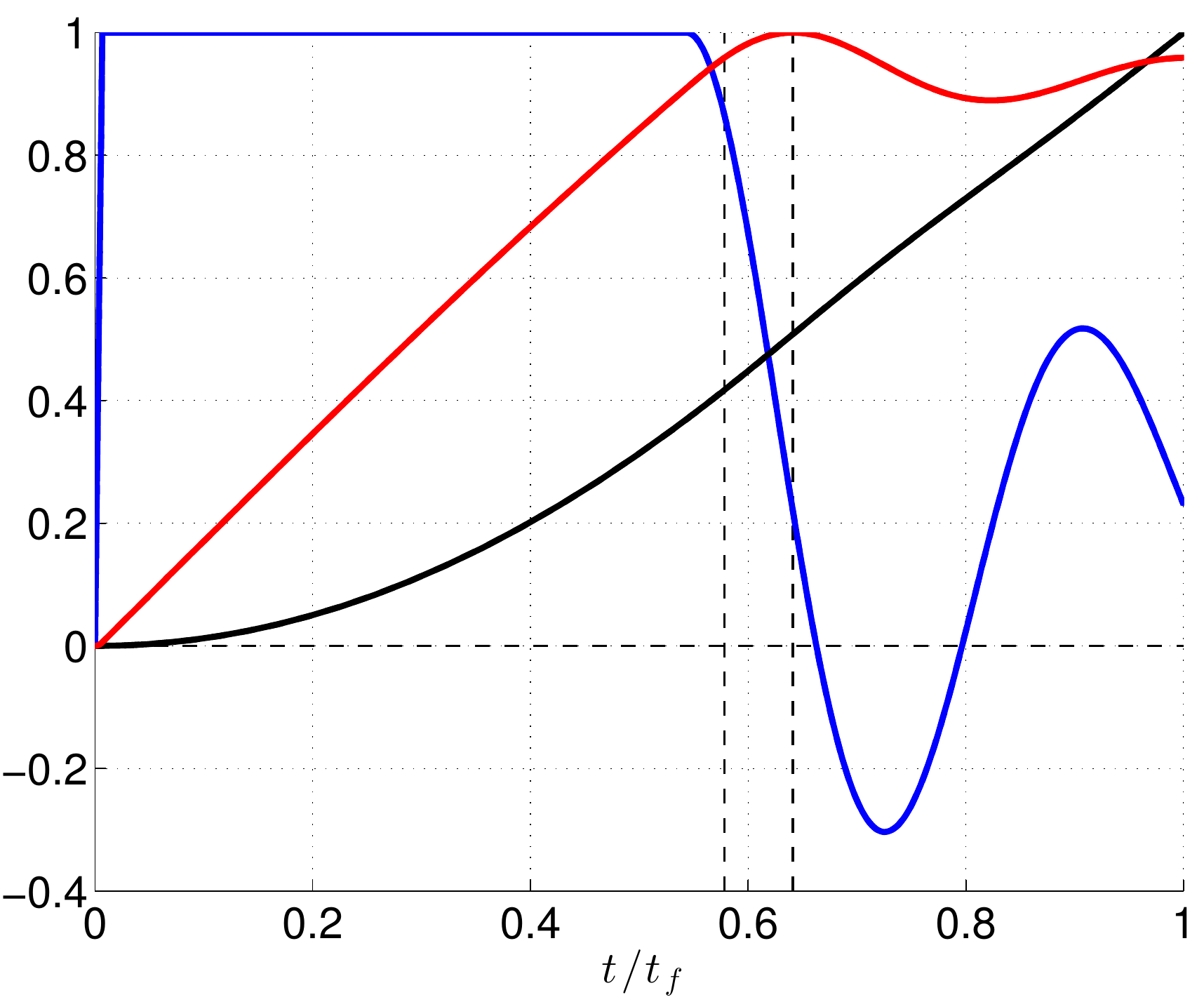}
                    \pxcoordinate{0.61*\x}{0.985*\y}{A1}; \draw (A1) node {\small $\tau_1$};
                    \pxcoordinate{0.67*\x}{0.985*\y}{A2}; \draw (A2) node {\small $\tau_2$};
                \end{tikzgraphics}
                \hspace{3em}
                \def\x{483}
                \def\y{411}
                \begin{tikzgraphics}{\sizeFig\textwidth}{\x}{\y}{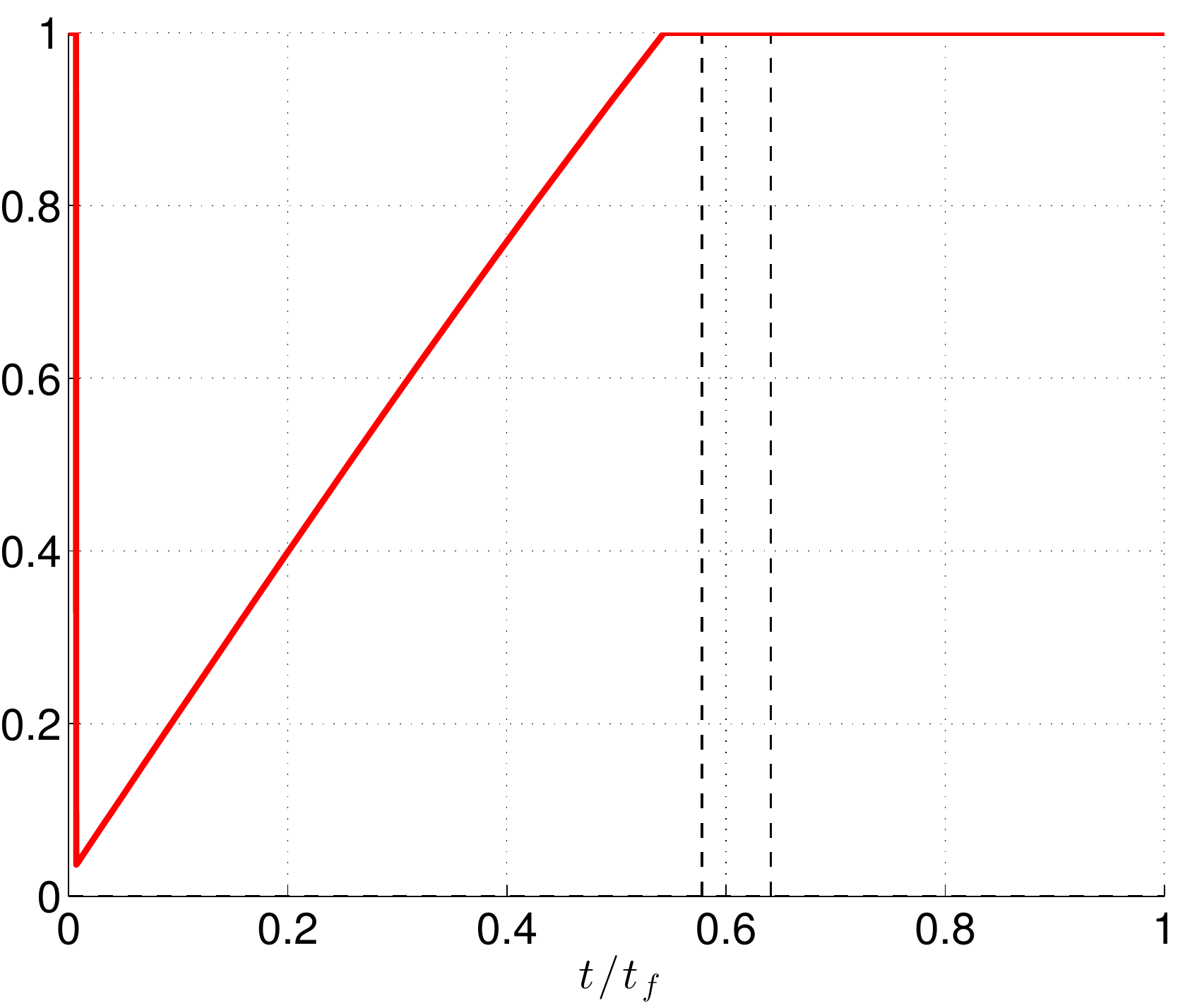}
                    \pxcoordinate{0.60*\x}{0.985*\y}{A1}; \draw (A1) node {\small $\tau_1$};
                    \pxcoordinate{0.66*\x}{0.985*\y}{A2}; \draw (A2) node {\small $\tau_2$};
                \end{tikzgraphics}
                \medskip

                \vspace{-1em}
                \def\x{493}
                \def\y{411}
                \begin{tikzgraphics}{\sizeFig\textwidth}{\x}{\y}{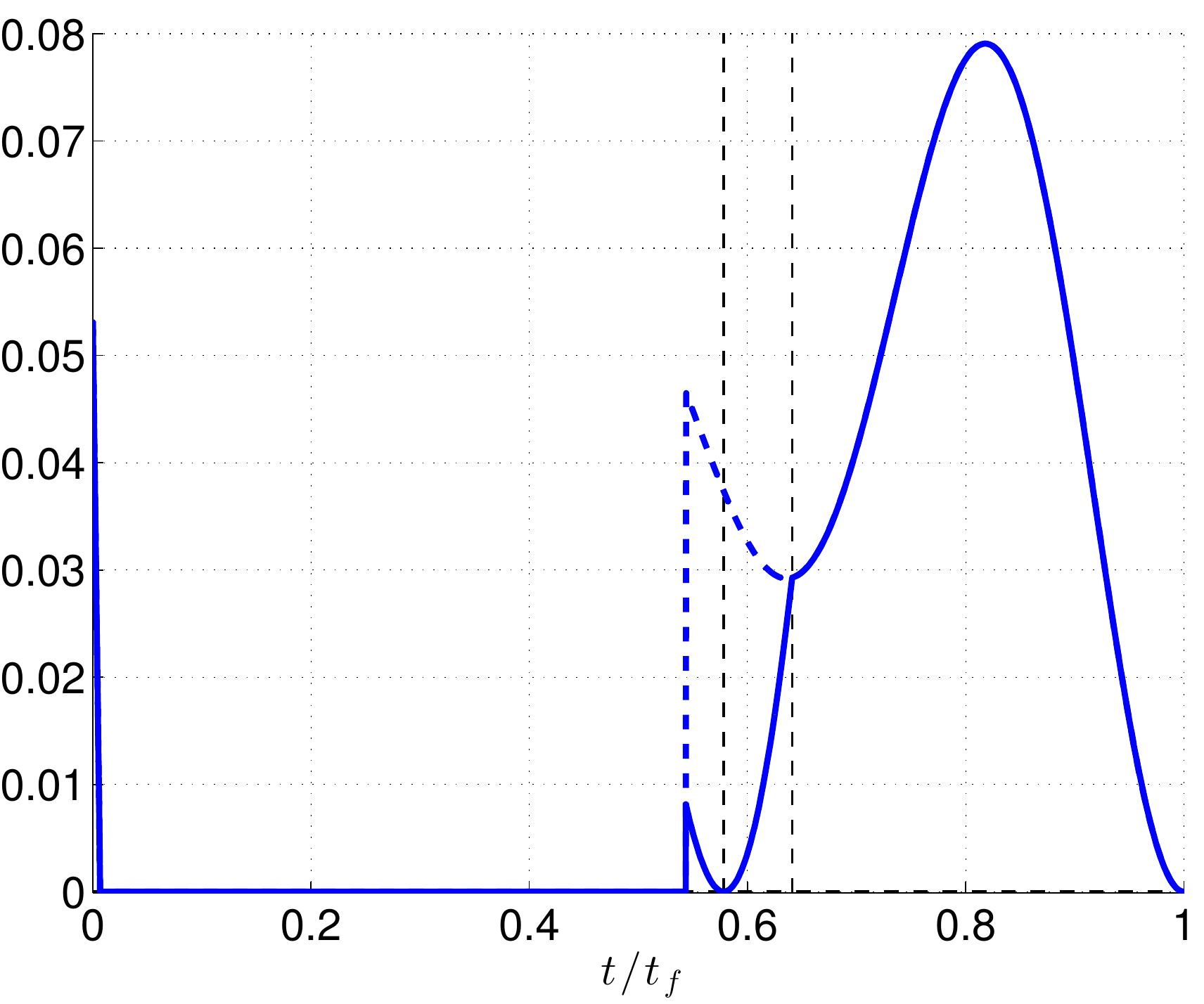}
                    \pxcoordinate{0.61*\x}{0.985*\y}{A1}; \draw (A1) node {\small $\tau_1$};
                    \pxcoordinate{0.67*\x}{0.985*\y}{A2}; \draw (A2) node {\small $\tau_2$};
                \end{tikzgraphics}
                \hspace{3em}
                \def\x{483}
                \def\y{411}
                \begin{tikzgraphics}{\sizeFig\textwidth}{\x}{\y}{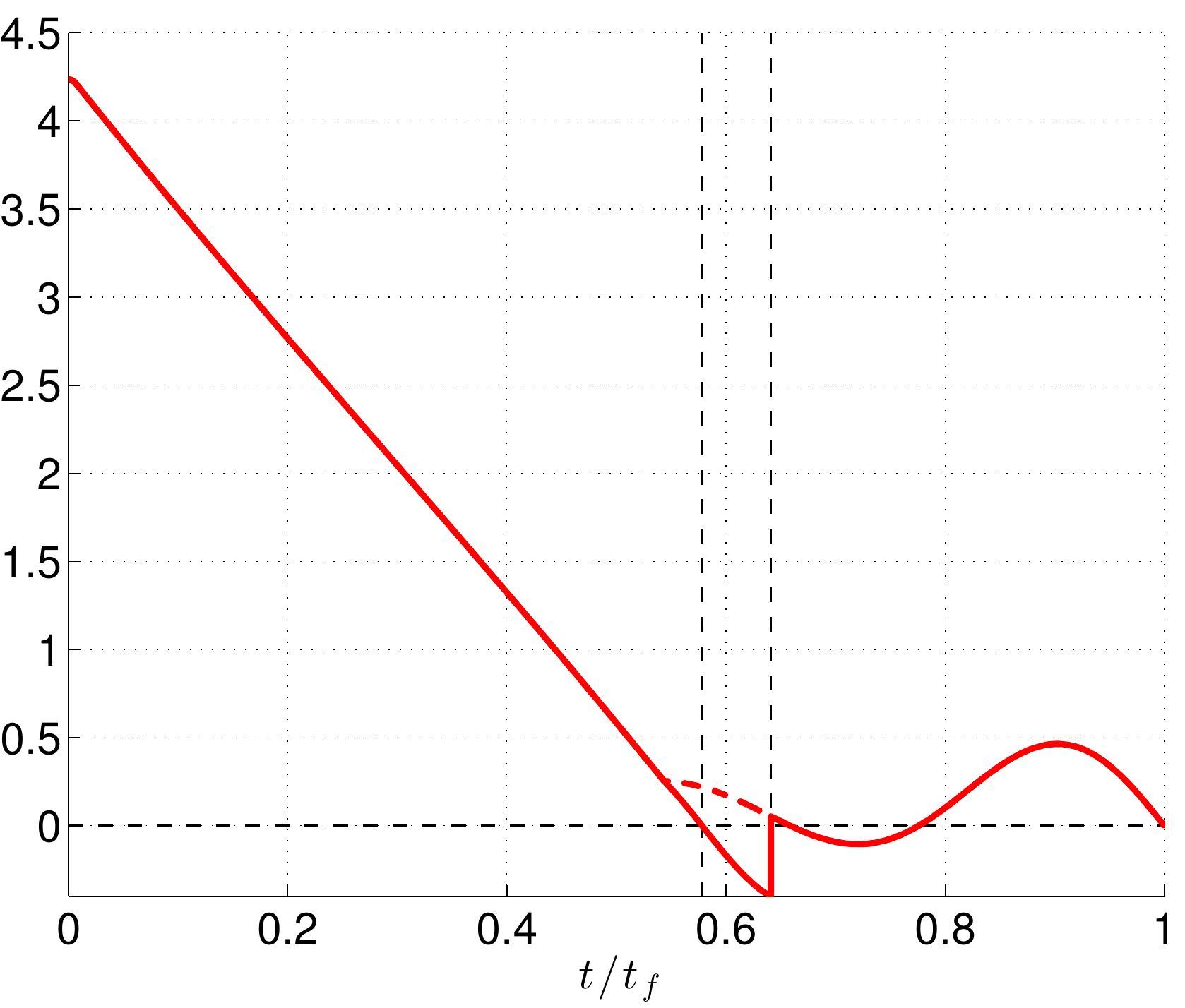}
                    \pxcoordinate{0.60*\x}{0.985*\y}{A1}; \draw (A1) node {\small $\tau_1$};
                    \pxcoordinate{0.66*\x}{0.985*\y}{A2}; \draw (A2) node {\small $\tau_2$};
                \end{tikzgraphics}
                \vspace{-1em}
                \caption{\textbf{Trajectories $\gamma_+ \gamma_{c_1} \gamma_+^{{c_3}}$ and $\gamma_+ \gamma_{c_1} \gamma_+^{{H_1},{c_3}}$}.
                The BC-extremal $\gamma_+ \gamma_{c_1} \gamma_+^{{c_3}}$ at $\vmax = \vmax^{c_3}$ from $\Hom_2^{(b)}$
                is represented by dashed lines
                while the one $\gamma_+ \gamma_{c_1} \gamma_+^{{H_1},{c_3}}$ from $\Hom_3$ is represented by solid lines.
                The contact with the boundary $C_3$ is at time $\tau_2$, while the contact of order 2 with $\Sigma_1^0$ is at time $\tau_1 < \tau_2$.
                (Top-Left) The state trajectories (compare with Figure~\ref{fig:gamma+}).
                (Top-Right) The control.
                (Bottom-Left) The component $p_1(\cdot)$.
                (Bottom-Right) The component $p_3(\cdot)$.}
                \label{fig:limiteCaseHom2Hom3}
            \end{figure}

        \subsubsection{Homotopy $\Hom_4$ and intermediate case $\gamma_+ \gamma_{c_1} \gamma_+ \gamma_- \gamma_{c_3}$,
                        with $u_{c_3}(\cdot) \equiv +1$, between $\Hom_3$ and $\Hom_4$}
        \label{sec:Hom4}

            The BC-extremal at $\vmax = \vmax^{\gamma_{c_3}}$ from $\Hom_3$ is $\gamma_+ \gamma_{c_1} \gamma_+ \gamma_- \gamma_+^{c_3}$ with $c_3(\cdot) \equiv 0$ along the
            last bang arc. For this specific value of $\vmax$, the boundary control $u_{c_3}$ is admissible if we replace the last bang arc by an arc $\gamma_{c_3}$.
            For $\vmax = \vmax^{\gamma_{c_3}} - \varepsilon$, $\varepsilon > 0$ small, the trajectory becomes $\gamma_+ \gamma_{c_1} \gamma_+ \gamma_- \gamma_{c_3}$ and the
            boundary control along the arc $\gamma_{c_3}$ is strictly admissible, \ie assumption \ref{Hyp3} holds.
            Here again, we use the last solution from $\Hom_3$ to initialize and solve 
            $S_4(y_4) = 0$, $y_4  \coloneqq  (p_0,t_f,t_1,t_2,\nu_2,t_3,t_4,\nu_4,z_1,z_2,z_3,z_4),$
            with $\vmax = \vmax^{\gamma_{c_3}}$. Then we solve $\Hom_4(y_4,\vmax) = 0$
            for $\vmax \in \intervalleff{\vmax^{+}}{\vmax^{\gamma_{c_3}}}$ starting from
            $\vmax = \vmax^{\gamma_{c_3}}$. The homotopy process has to stop when $t_2 \ge t_3$ (or $\nu_2 \ge 0$), which gives $\vmax^{+} \approx 64.1641$.
            Figure~\ref{fig:homotopy4_times} shows the initial, the switching, the junction and the final times along the path of zeros of $\Hom_4$,
            while the jumps $\nu_2$ and $\nu_4$ are portrayed in Figure~\ref{fig:homotopy4_jumps}. Figure~\ref{fig:limiteCaseHom3Hom4} gives details about
            the limit case when $\vmax = \vmax^{+}$.

            \def\sizeFig{0.26}
            \begin{figure}[ht!]
                \def\posg{0.14}
                \def\posc{0.915}
                \centering
                \def\x{492}
                \def\y{403}
                \begin{tikzgraphics}{\sizeFig\textwidth}{\x}{\y}{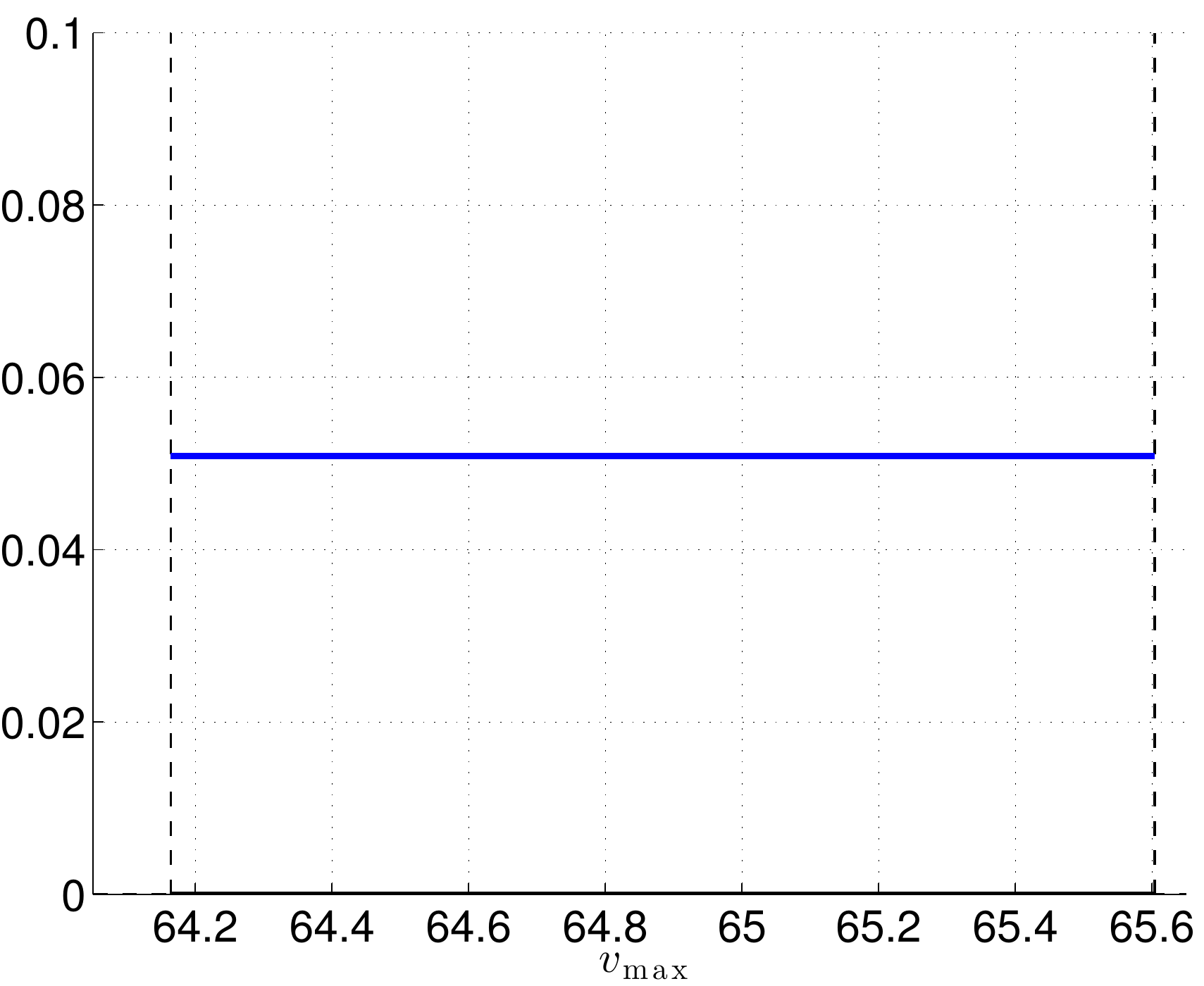}
                    \pxcoordinate{\posg*\x}{-0.03*\y}{A}; \draw (A) node {$\vmax^{+}$};
                    \pxcoordinate{\posc*\x}{-0.03*\y}{A}; \draw (A) node {$\vmax^{\gamma_{c_3}}$};
                \end{tikzgraphics}
                \def\x{492}
                \def\y{393}
                \begin{tikzgraphics}{\sizeFig\textwidth}{\x}{\y}{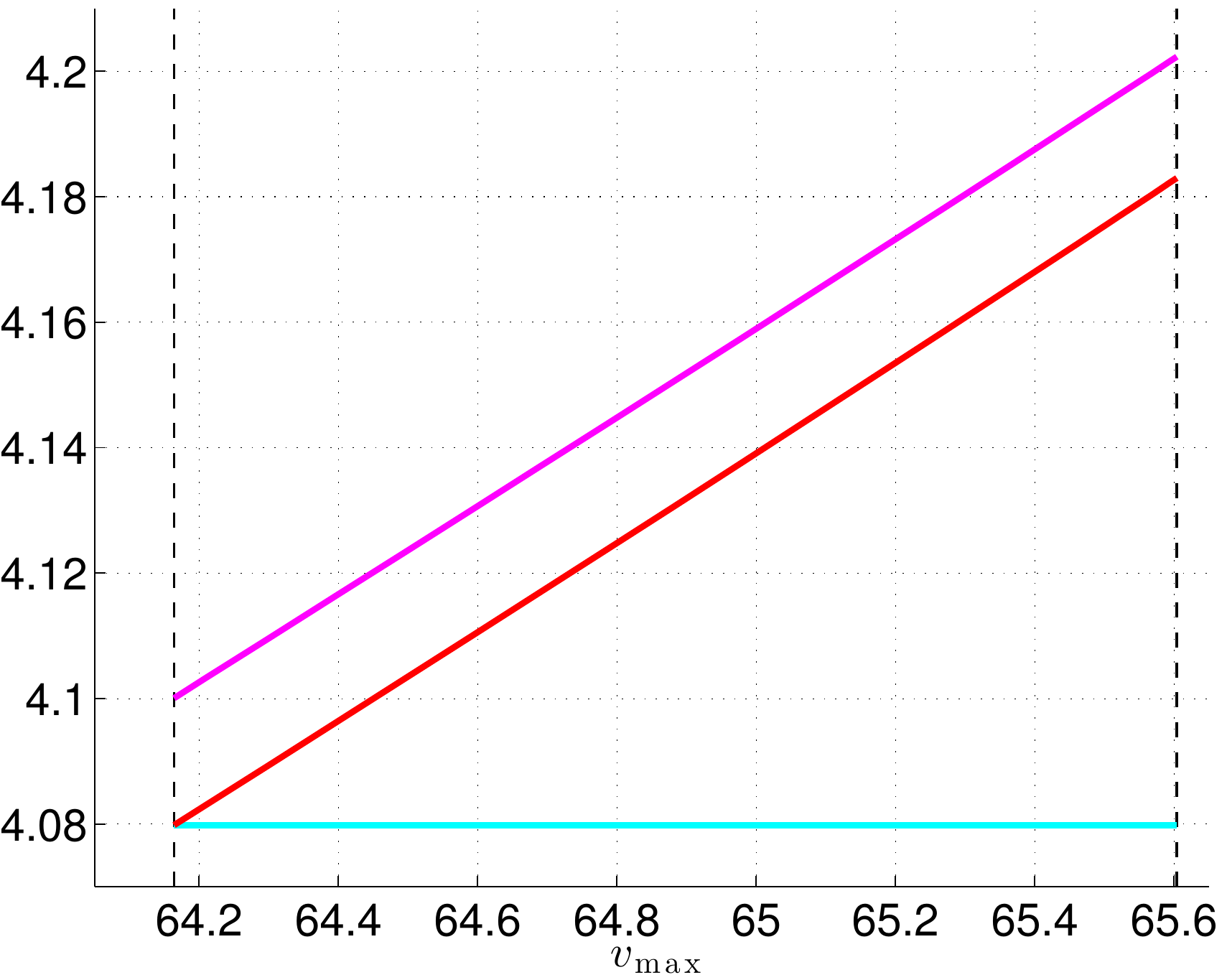}
                    \pxcoordinate{\posg*\x}{-0.06*\y}{A}; \draw (A) node {$\vmax^{+}$};
                    \pxcoordinate{\posc*\x}{-0.06*\y}{A}; \draw (A) node {$\vmax^{\gamma_{c_3}}$};
                \end{tikzgraphics}
                \def\x{492}
                \def\y{393}
                \begin{tikzgraphics}{\sizeFig\textwidth}{\x}{\y}{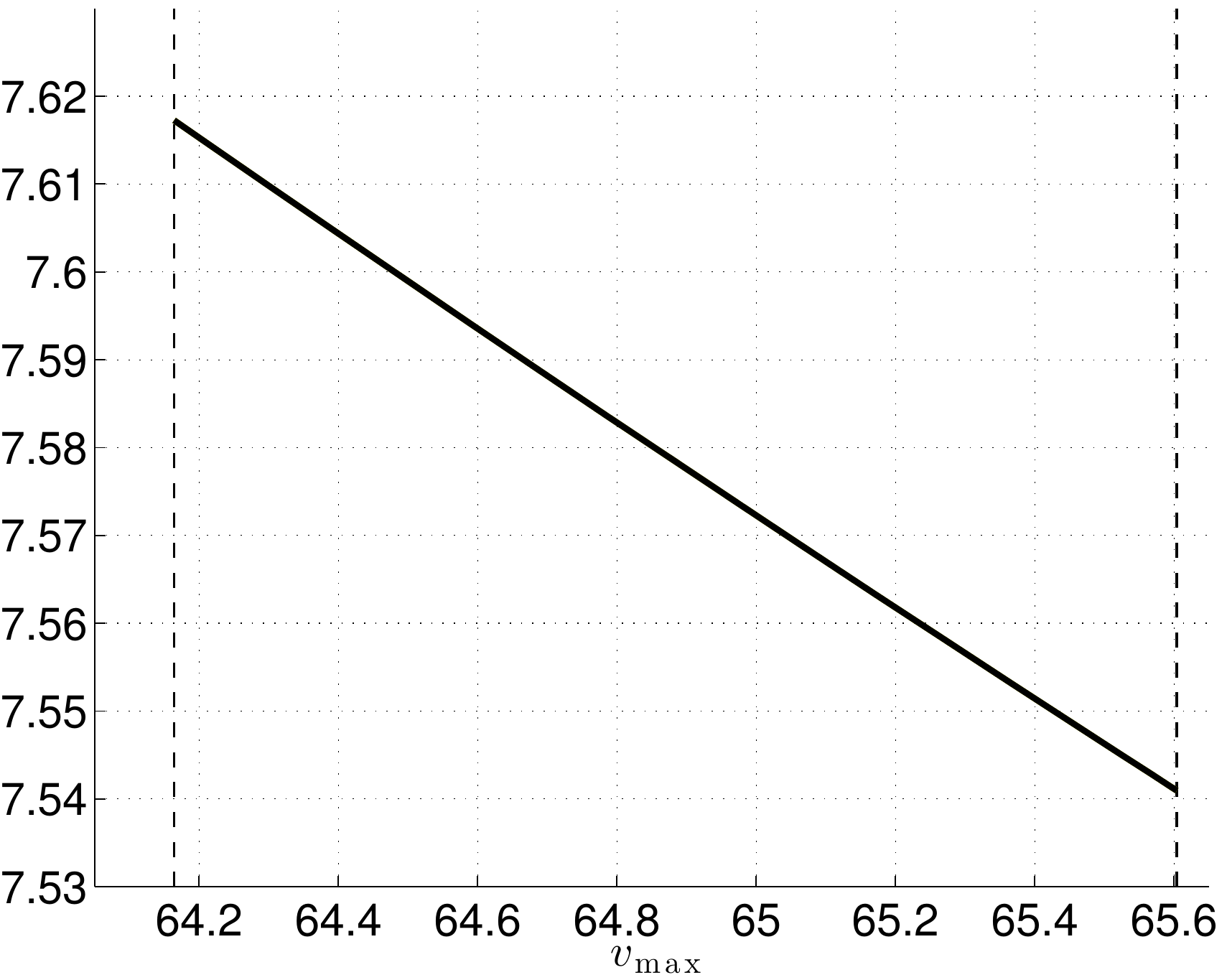}
                    \pxcoordinate{\posg*\x}{-0.06*\y}{A}; \draw (A) node {$\vmax^{+}$};
                    \pxcoordinate{\posc*\x}{-0.06*\y}{A}; \draw (A) node {$\vmax^{\gamma_{c_3}}$};
                \end{tikzgraphics}
                \vspace{-1em}
                \caption{\textbf{Homotopy $\Hom_4$}.  The initial, the switching, the junction and the final times along the path of zeros of $\Hom_4$.
                (Left) The times $t_0 \equiv 0$ and $t_1$.
                (Middle) The times $t_2$, $t_3$ and $t_4$.
                (Right) The final time $t_f$.
                The length $t_3-t_2$ vanishes when $\vmax = \vmax^{+}$.}
                \label{fig:homotopy4_times}
            \end{figure}
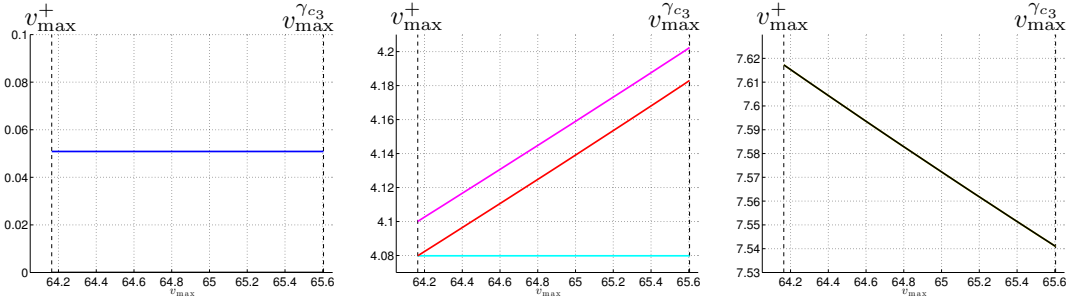

            \def\sizeFig{0.3}
            \begin{figure}[ht!]
                \def\posg{0.155}
                \def\posc{0.95}
                \centering
                \def\x{505}
                \def\y{420}
                \begin{tikzgraphics}{\sizeFig\textwidth}{\x}{\y}{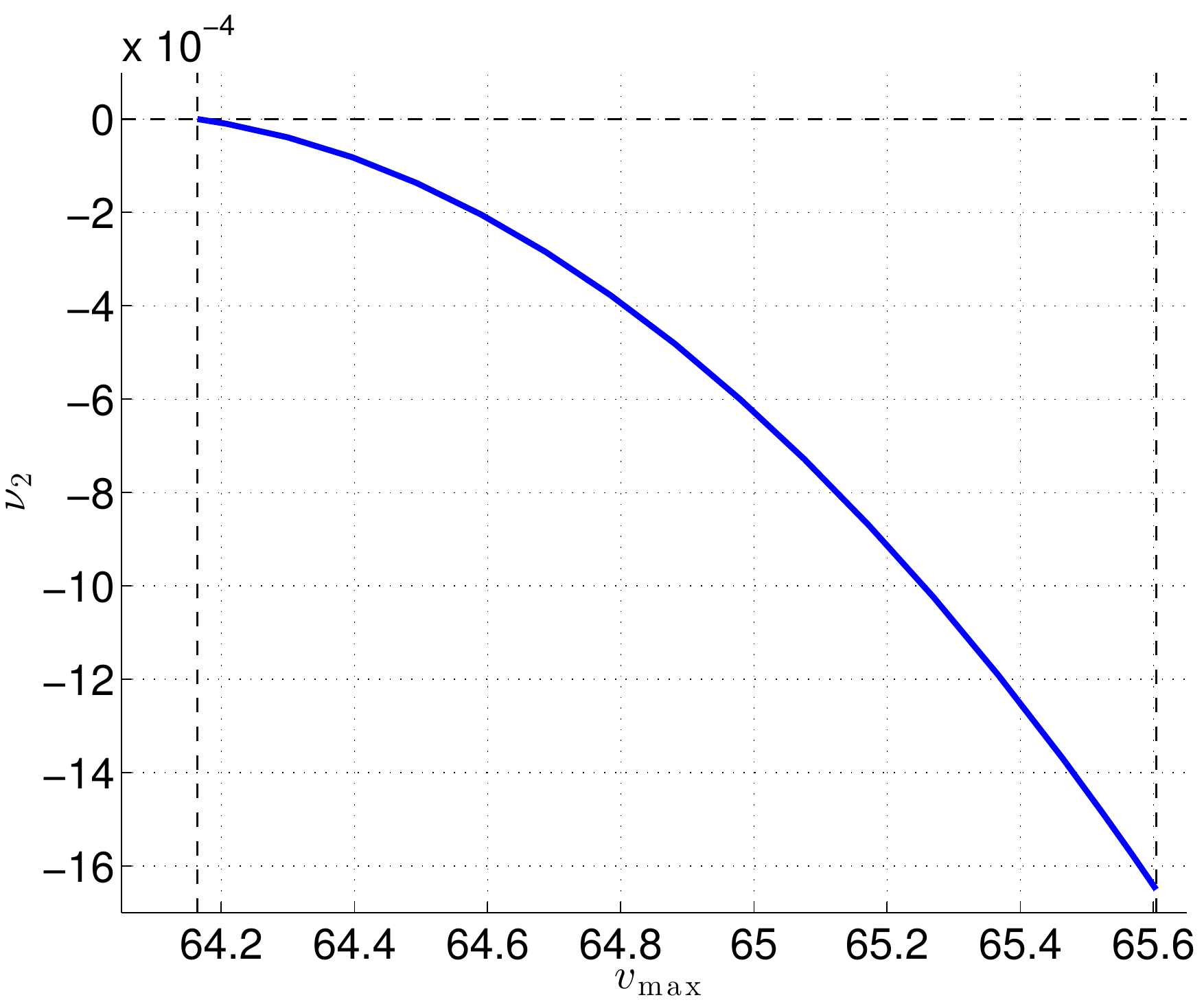}
                    \pxcoordinate{\posg*\x}{-0.03*\y}{A}; \draw (A) node {$\vmax^{+}$};
                    \pxcoordinate{\posc*\x}{-0.03*\y}{A}; \draw (A) node {$\vmax^{\gamma_{c_3}}$};
                \end{tikzgraphics}
                \hspace{3em}
                \def\x{520}
                \def\y{420}
                \begin{tikzgraphics}{\sizeFig\textwidth}{\x}{\y}{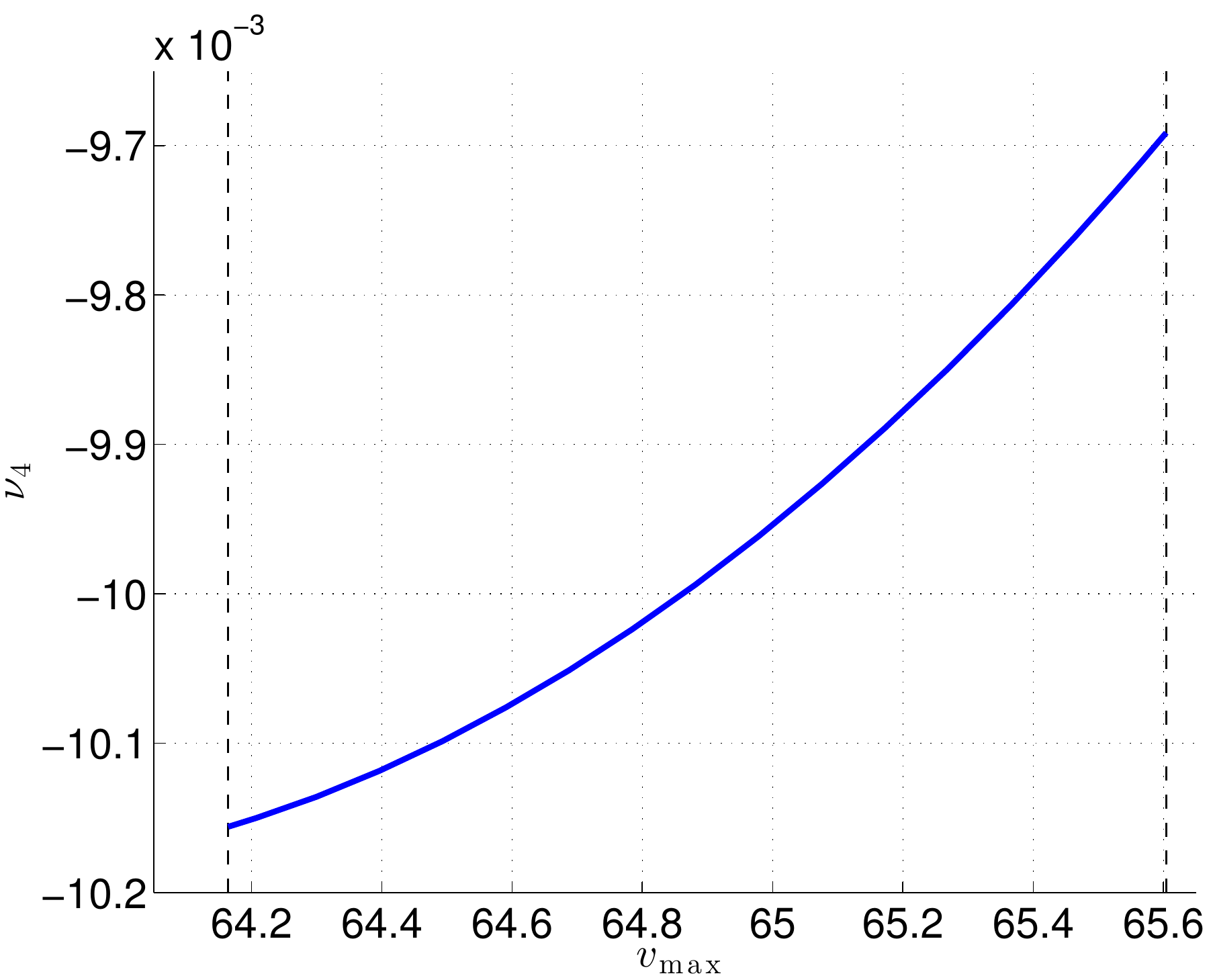}
                    \pxcoordinate{0.18*\x}{-0.05*\y}{A}; \draw (A) node {$\vmax^{+}$};
                    \pxcoordinate{\posc*\x}{-0.05*\y}{A}; \draw (A) node {$\vmax^{\gamma_{c_3}}$};
                \end{tikzgraphics}
                \caption{\textbf{Homotopy $\Hom_4$}. The jumps $\nu_2 \le 0$ (left) and $\nu_4 \le 0$ (right) along the path of zeros of $\Hom_3$.
                Note that the jump $\nu_2$ vanishes when $\vmax = \vmax^{+}$, see lemma \ref{lemma:junction1}. The jump $\nu_4 < 0$ according
                to lemma \ref{lemma:junction2}.
                }
                \label{fig:homotopy4_jumps}
            \end{figure}

            \def\sizeFig{0.3}
            \begin{figure}[ht!]
                \centering
                \includegraphics[width=\sizeFig\textwidth]{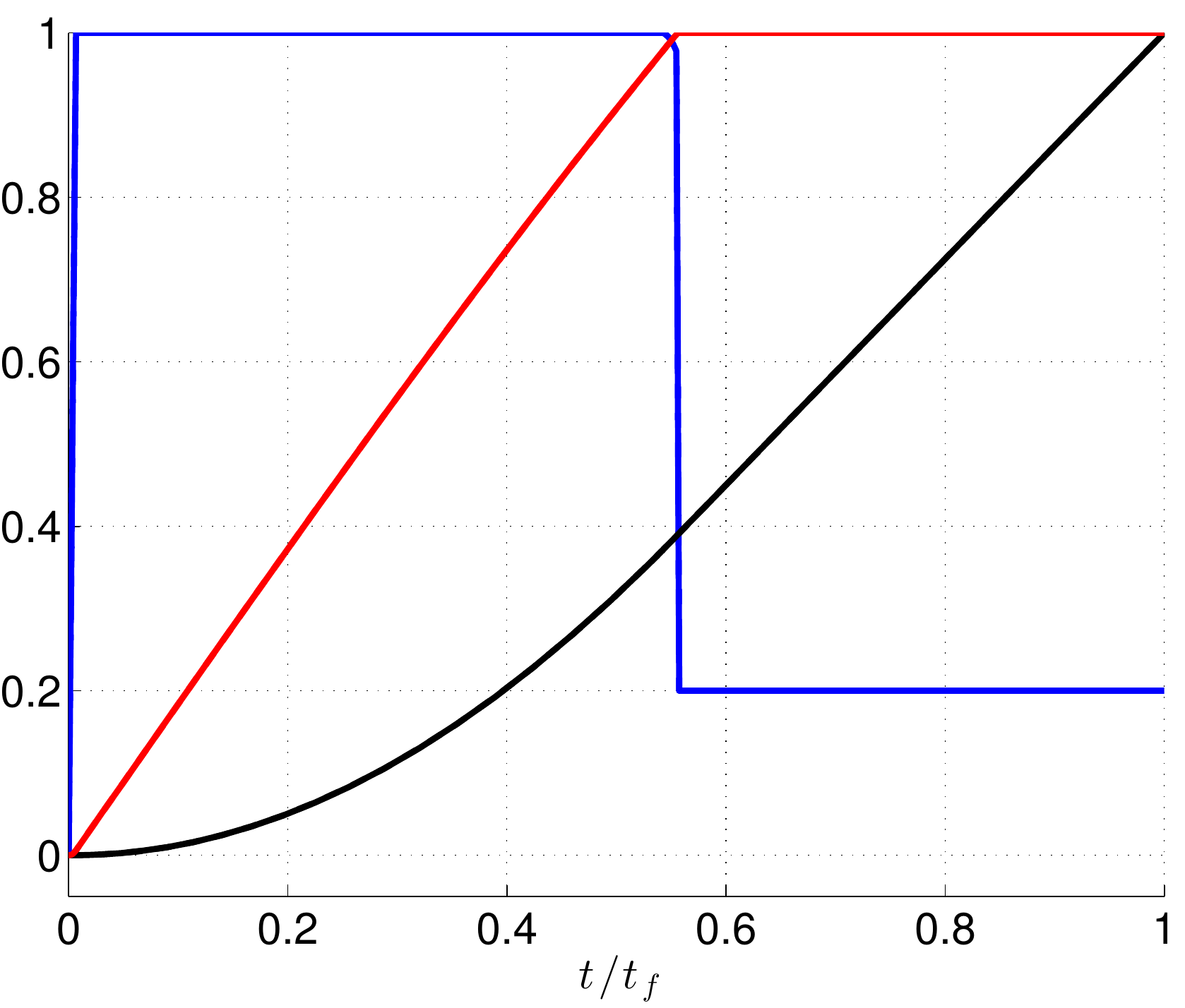}
                \hspace{3em}
                \includegraphics[width=\sizeFig\textwidth]{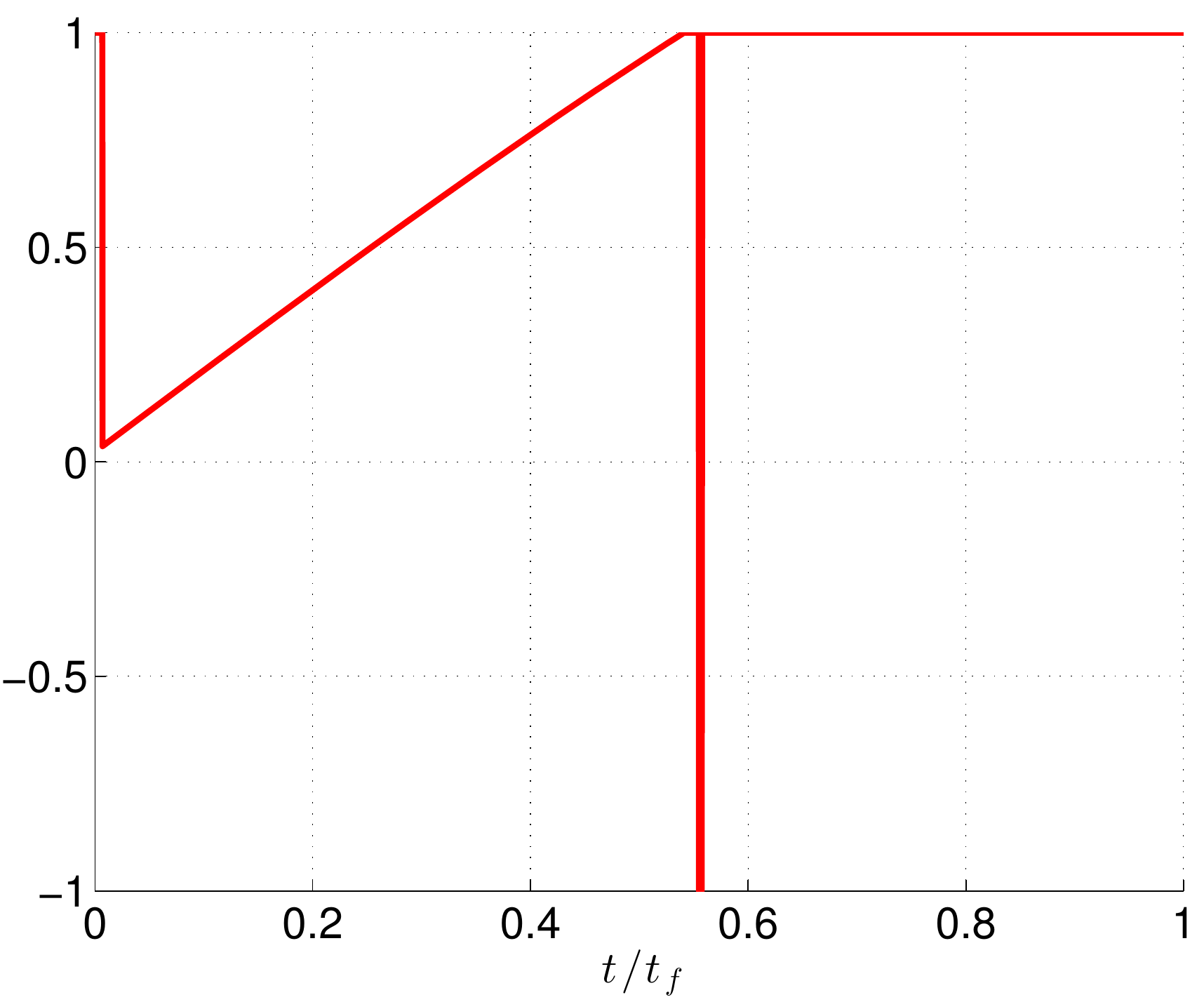}
                \medskip

                \includegraphics[width=\sizeFig\textwidth]{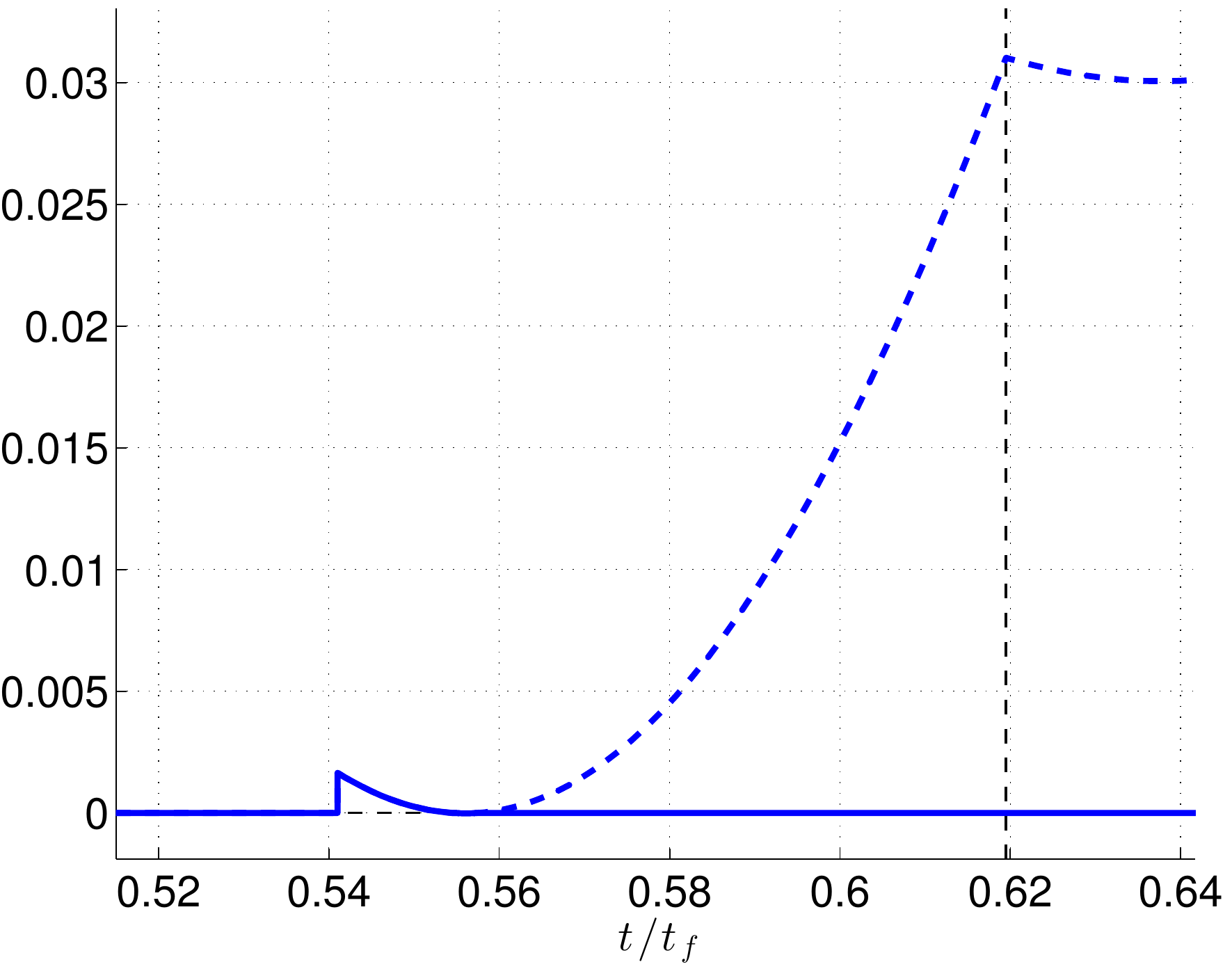}
                \hspace{3em}
                \includegraphics[width=\sizeFig\textwidth]{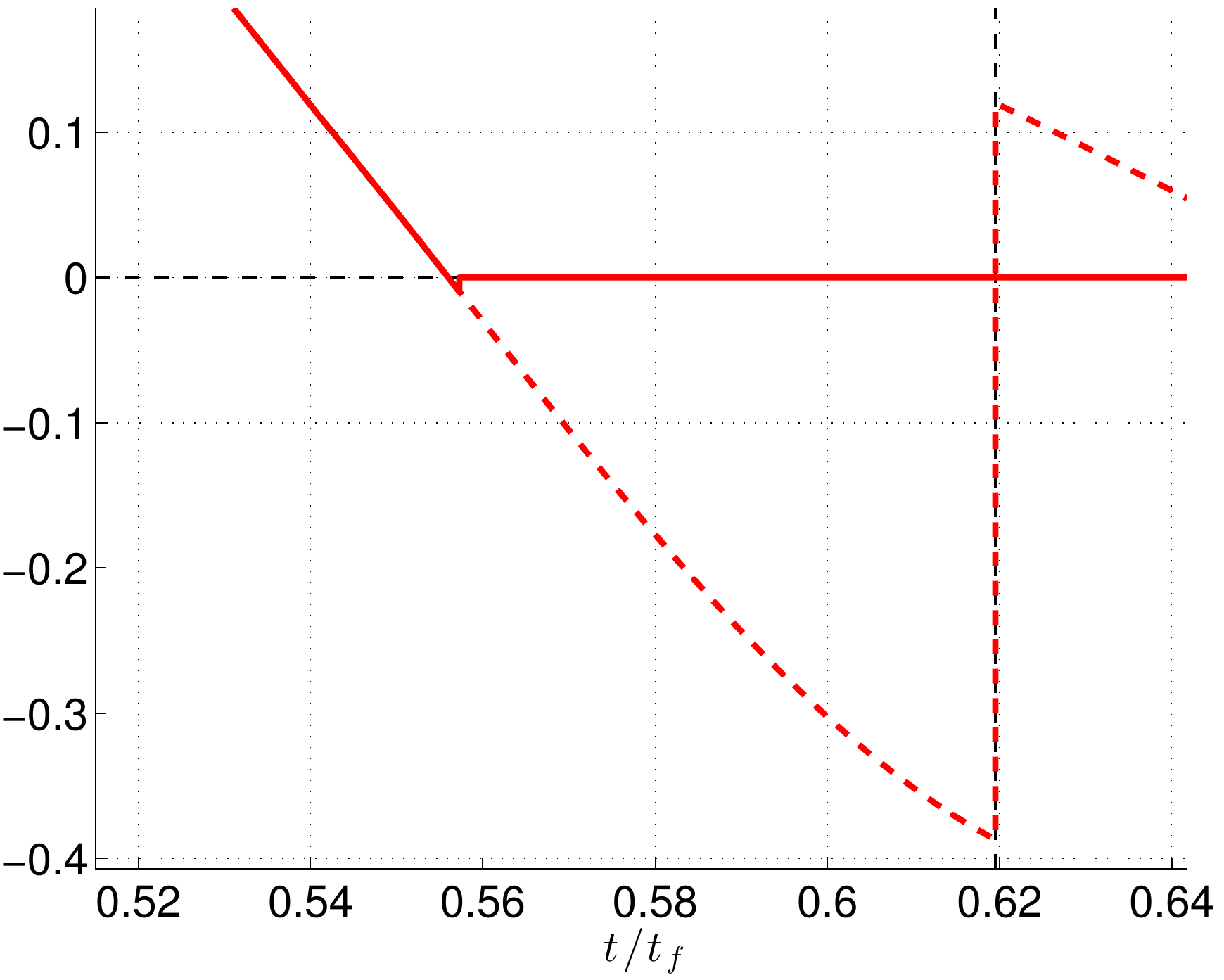}
                \caption{\textbf{Trajectories $\gamma_+ \gamma_{c_1} \gamma_+ \gamma_- \gamma_+^{c_3}$ and $\gamma_+ \gamma_{c_1} \gamma_+ \gamma_- \gamma_{c_3}$,
                with $u_{c_3}(\cdot) = +1$}.
                The BC-extremal $\gamma_+ \gamma_{c_1} \gamma_+ \gamma_- \gamma_+^{c_3}$ at $\vmax = \vmax^{\gamma_{c_3}}$ from $\Hom_3$
                is represented by dashed lines
                while the one $\gamma_+ \gamma_{c_1} \gamma_+ \gamma_- \gamma_{c_3}$ from $\Hom_4$ is represented by solid lines.
                (Top-Left) The state trajectories.
                (Top-Right) The control.
                (Bottom-Left) Zoom on the component $p_1(\cdot)$.
                (Bottom-Right) Zoom on the component $p_3(\cdot)$.
                About the trajectory $\gamma_+ \gamma_{c_1} \gamma_+ \gamma_- \gamma_{c_3}$ with $u_{c_3}(\cdot) = +1$,
                the components $p_1(\cdot)$ and $p_3(\cdot)$ are identically zero along $\gamma_{c_3}$ since
                this extremal is a limit case from $\Hom_4$ for which for each solution, the functions $\Phi$ and $\dot{\Phi}$ are zero along $\gamma_{c_3}$.
                Note that the sine wave shapes of the trajectory at Figure~\ref{fig:limiteCaseHom1Hom2} have changed into square wave shapes.}
                \label{fig:limiteCaseHom3Hom4}
            \end{figure}

        \subsubsection{Homotopy $\Hom_5$ and intermediate case $\gamma_+ \gamma_{c_1} \gamma_+ \gamma_- \gamma_{c_3}$, with $t_2 = t_3$, between $\Hom_4$ and $\Hom_5$}
        \label{sec:Hom5}

        Here is the last homotopy. We use the solution from $\Hom_4$ at $\vmax = \vmax^{+}$ to initialize and solve $S_5(y_5) = 0$,
        $y_5  \coloneqq  (p_0,t_f,t_1,t_2,t_3,\nu_3,z_1,z_2,z_3).$ 
        This is done easily since $t_3$, $\nu_3$ and $z_3$ are given respectively by $t_4$, $\nu_4$ and $z_4$ from
        the solution from $\Hom_4$. In other words, at $\vmax = \vmax^{+}$, the extremals are the same from both the solutions from $\Hom_4$ and $\Hom_5$.
        This is new compare to the others homotopies\footnote{Except obviously between $\Hom_2^{(a)}$ and $\Hom_2^{(b)}$.}, for which the limit cases have different adjoint
        vectors. Figures~\ref{fig:limiteCaseHom4Hom5_1} and \ref{fig:limiteCaseHom4Hom5_2} 
        only shows the components $p_1(\cdot)$ and $p_3(\cdot)$ of the extremal for $\vmax = \vmax^{+}$ since
        the trajectory is quite similar from the one in Figure~\ref{fig:limiteCaseHom3Hom4} and the shape of the control 
        can be guessed from Figure~\ref{fig:limiteCaseHom3Hom4}.
        Figure~\ref{fig:homotopy5_times} shows the initial, the junction and the final times along the path of zeros of $\Hom_5$,
        while the jumps $\nu_3$ and the length $t_3 - t_2$ of the arc $\gamma_-$ are portrayed in Figure~\ref{fig:homotopy5_ecart_saut}.
        One may note that the final time has not a linear shape with respect to $\vmax$. It increases more and more when $\vmax$ decreases,
        contrary to the length of the boundary arc $\gamma_{c_1}$ which decreases linearly. Besides, the
        length of $\gamma_-$ is small but not constant.

            \def\sizeFig{0.3}
            \begin{figure}[ht!]
                \def\posg{0.095}
                \def\posc{0.96}
                \centering
                \def\x{489}
                \def\y{403}
                \begin{tikzgraphics}{\sizeFig\textwidth}{\x}{\y}{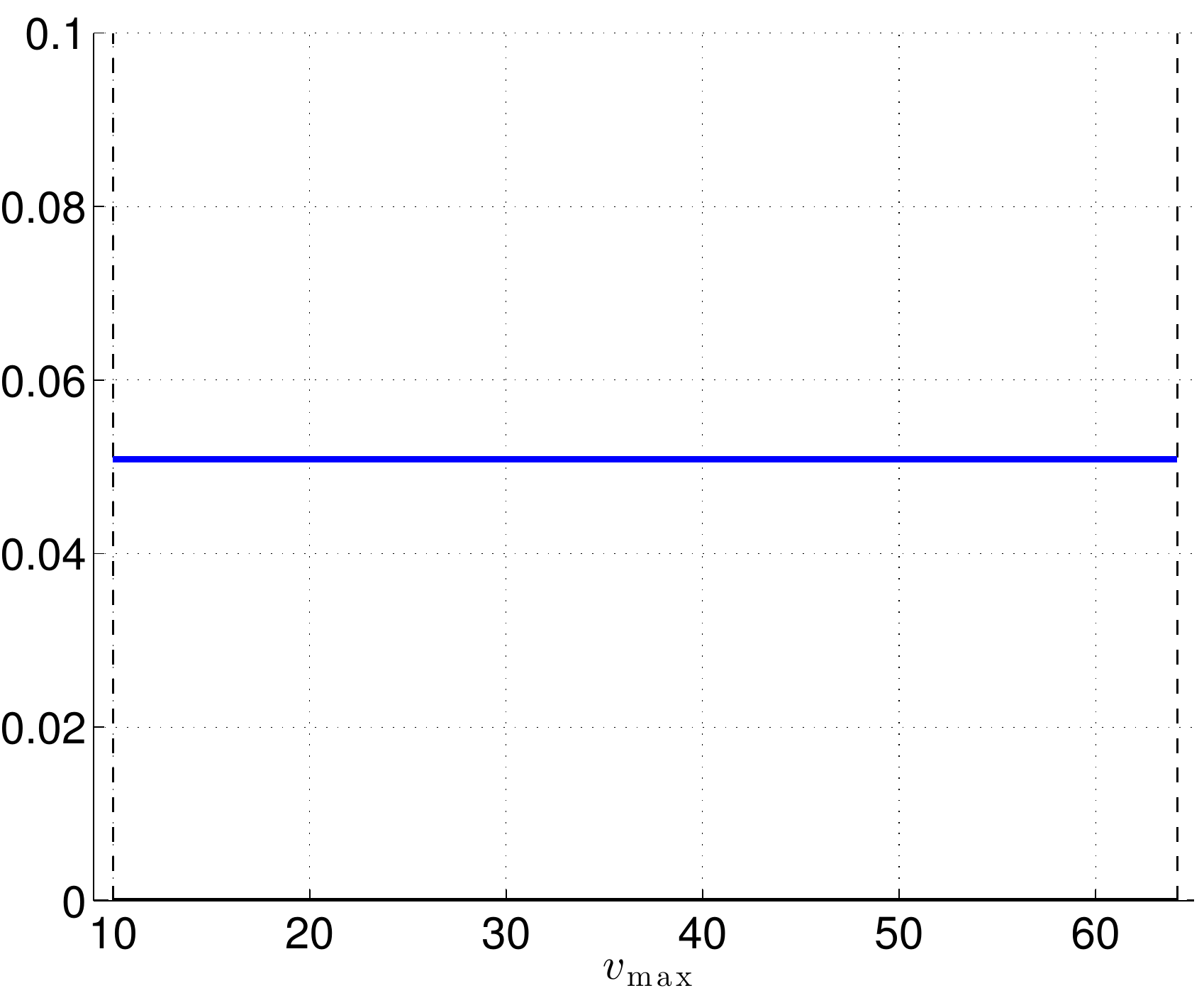}
                    \pxcoordinate{\posg*\x}{-0.03*\y}{A}; \draw (A) node {$10$};
                    \pxcoordinate{\posc*\x}{-0.03*\y}{A}; \draw (A) node {$\vmax^{+}$};
                \end{tikzgraphics}
                \hfill
                \def\x{479}
                \def\y{393}
                \begin{tikzgraphics}{\sizeFig\textwidth}{\x}{\y}{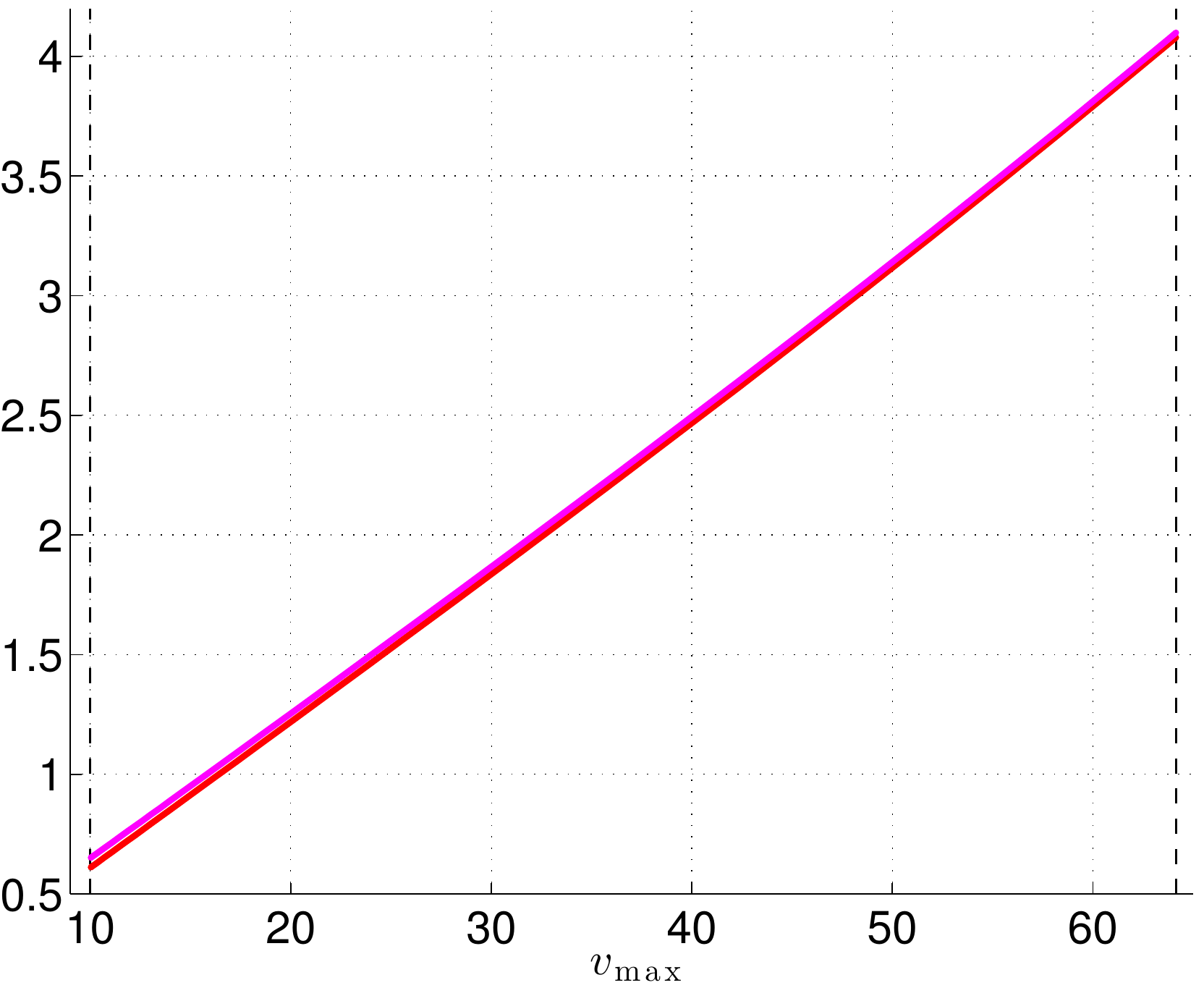}
                    \pxcoordinate{0.08*\x}{-0.04*\y}{A}; \draw (A) node {$10$};
                    \pxcoordinate{\posc*\x}{-0.04*\y}{A}; \draw (A) node {$\vmax^{+}$};
                \end{tikzgraphics}
                \hfill
                \def\x{474}
                \def\y{393}
                \begin{tikzgraphics}{\sizeFig\textwidth}{\x}{\y}{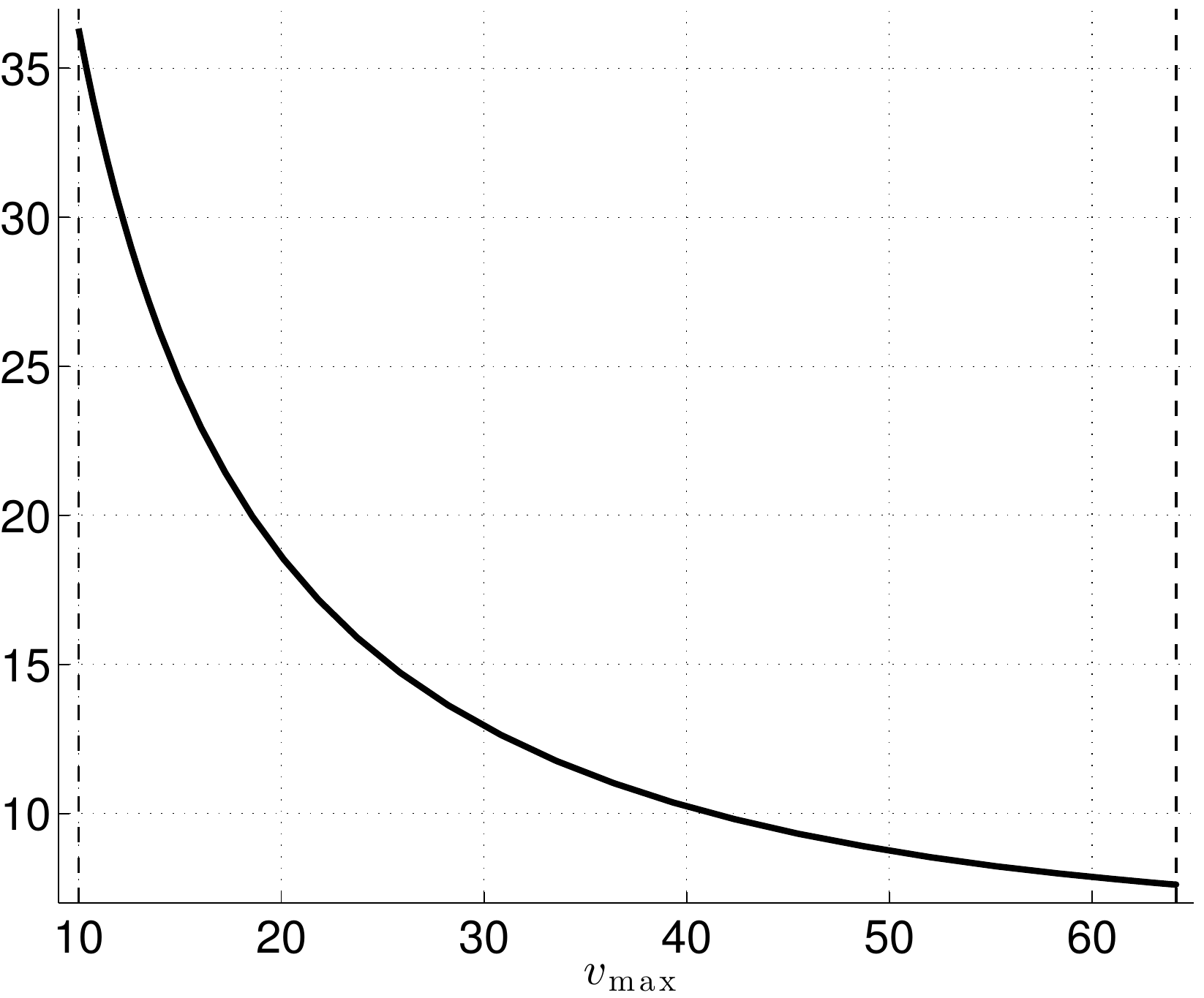}
                    \pxcoordinate{0.07*\x}{-0.04*\y}{A}; \draw (A) node {$10$};
                    \pxcoordinate{\posc*\x}{-0.04*\y}{A}; \draw (A) node {$\vmax^{+}$};
                \end{tikzgraphics}
                \caption{\textbf{Homotopy $\Hom_5$}.
                 The initial, the junction and the final times along the path of zeros of $\Hom_5$.
                (Left) The times $t_0 \equiv 0$ and $t_1$.
                (Middle) The times $t_2$, $t_3$.
                (Right) The final time $t_f$.
                The length $t_3-t_2$ is small but not constant, see Figure~\ref{fig:homotopy5_ecart_saut}. The length of the first bang arc is still constant.}
                \label{fig:homotopy5_times}
            \end{figure}

            \def\sizeFig{0.3}
            \begin{figure}[ht!]
                \def\posg{0.15}
                \def\posc{0.97}
                \centering
                \def\x{518}
                \def\y{393}
                \begin{tikzgraphics}{\sizeFig\textwidth}{\x}{\y}{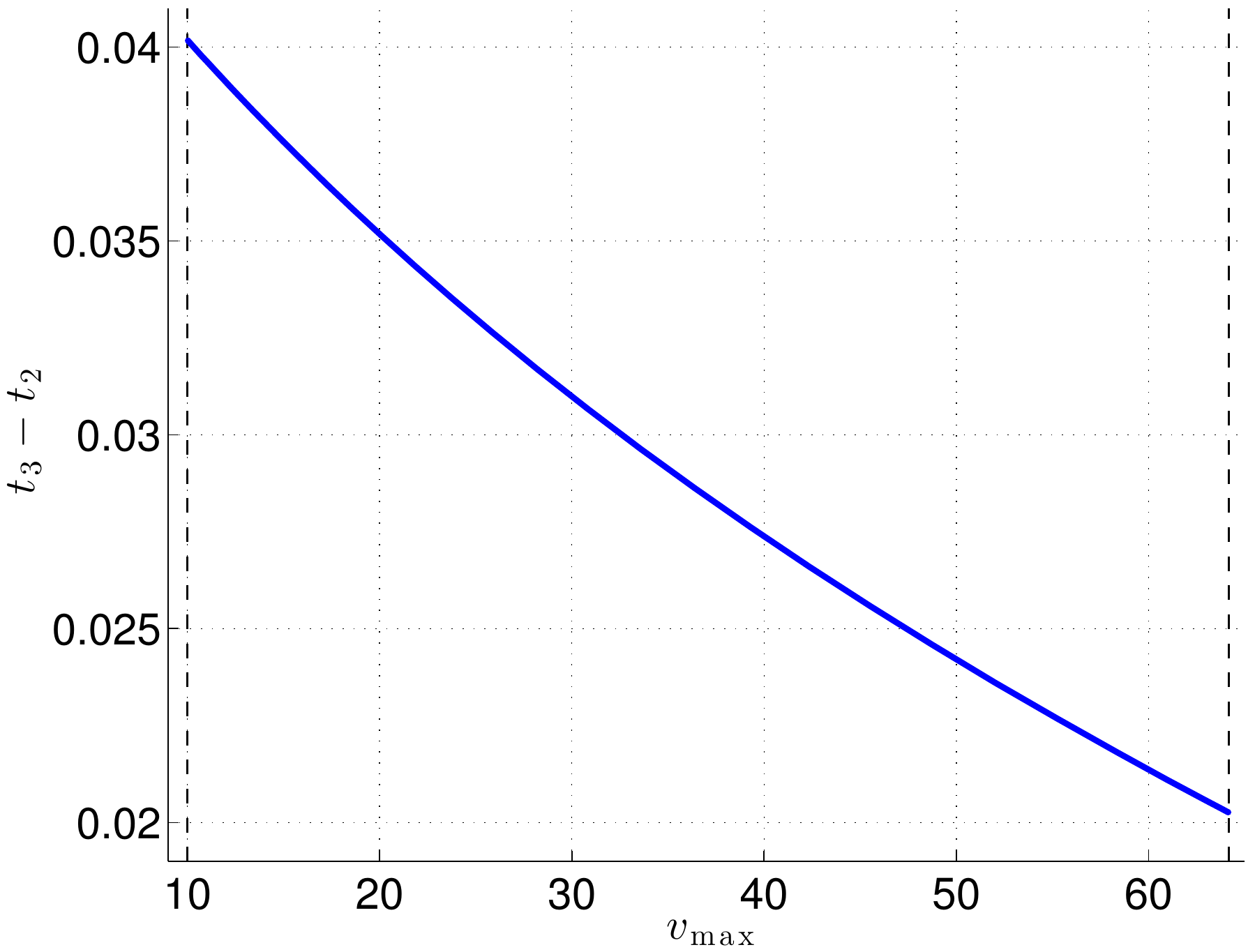}
                    \pxcoordinate{\posg*\x}{-0.04*\y}{A}; \draw (A) node {$10$};
                    \pxcoordinate{\posc*\x}{-0.04*\y}{A}; \draw (A) node {$\vmax^{+}$};
                \end{tikzgraphics}
                \hspace{0em}
                \def\x{526}
                \def\y{403}
                \begin{tikzgraphics}{\sizeFig\textwidth}{\x}{\y}{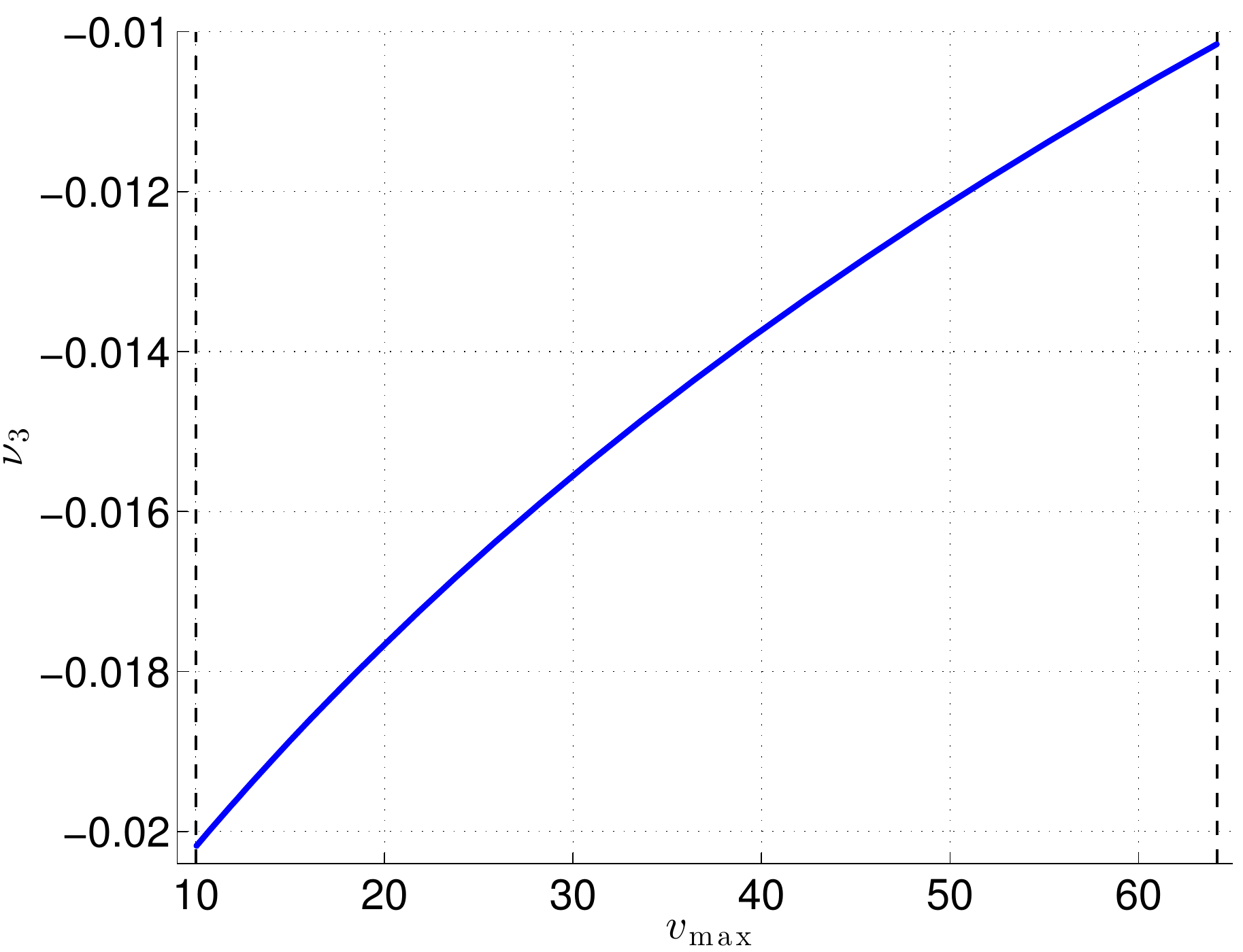}
                    \pxcoordinate{\posg*\x}{-0.03*\y}{A}; \draw (A) node {$10$};
                    \pxcoordinate{\posc*\x}{-0.03*\y}{A}; \draw (A) node {$\vmax^{+}$};
                \end{tikzgraphics}
                \caption{\textbf{Homotopy $\Hom_5$}. The length $t_3-t_2$ (left) and the jump $\nu_3$ (right) along $\Hom_5$.}
                \label{fig:homotopy5_ecart_saut}
            \end{figure}

            \def\sizeFig{0.36}
            \begin{figure}[ht!]
                \centering
                \def\x{493}
                \def\y{411}
                \begin{tikzgraphics}{\sizeFig\textwidth}{\x}{\y}{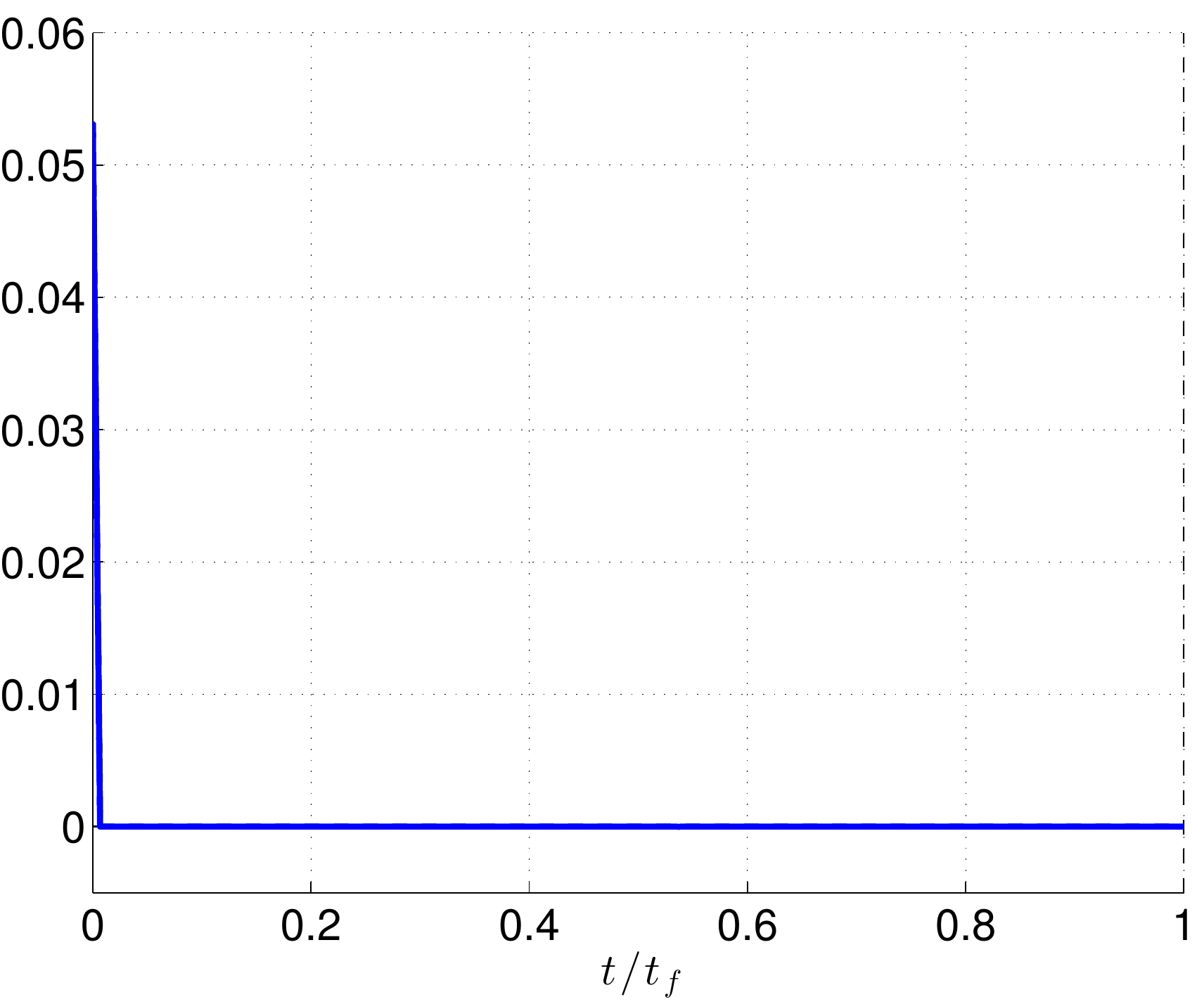}
                    \pxcoordinate{0.555*\x}{0.822*\y}{centre};
                    \draw [thin, dashed] (centre) circle (0.5);
                    \draw [-triangle 45,thin] ([xshift=0.85*\x,yshift=0.3*\y] centre) -- ([xshift=0.03*\x,yshift=0.025*\x] centre);
                \end{tikzgraphics}
                \hspace{-3em}
                \def\x{485}
                \def\y{428}
                \begin{tikzgraphics}{0.26\textwidth}{\x}{\y}{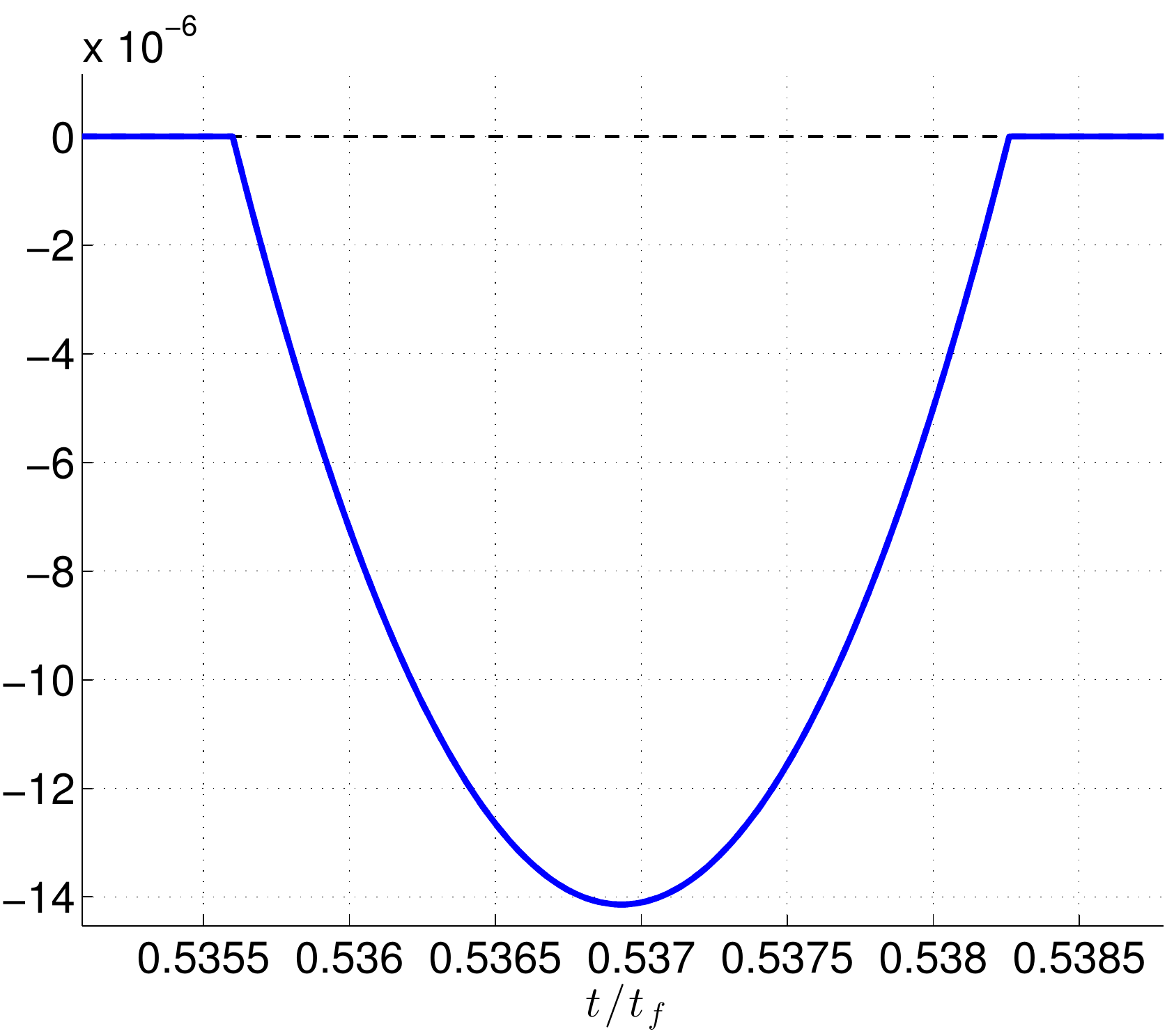}
                    \pxcoordinate{0.53*\x}{0.42*\y}{centre};
                    \draw [thin, dashed] (centre) circle (8);
                \end{tikzgraphics}
                \caption{\textbf{Trajectory $\gamma_+ \gamma_{c_1} \gamma_+ \gamma_- \gamma_{c_3}$, with $t_2 = t_3$.}
                 The BC-extremal $\gamma_+ \gamma_{c_1} \gamma_+ \gamma_- \gamma_{c_3}$ with $t_2 = t_3$ is the same as the
                BC-extremal $\gamma_+ \gamma_{c_1} \gamma_- \gamma_{c_3}$ from $\Hom_5$ at $\vmax = \vmax^{+}$.
                (Left) The component $p_1(\cdot)$ with a zoom (Right) around $\gamma_-$. $p_1(\cdot)$ is zero along $\gamma_{c_1}$ and $\gamma_{c_3}$.}
                \label{fig:limiteCaseHom4Hom5_1}
            \end{figure}

            \def\sizeFig{0.36}
            \begin{figure}[ht!]
                \centering
                \def\x{510}
                \def\y{411}
                \begin{tikzgraphics}{\sizeFig\textwidth}{\x}{\y}{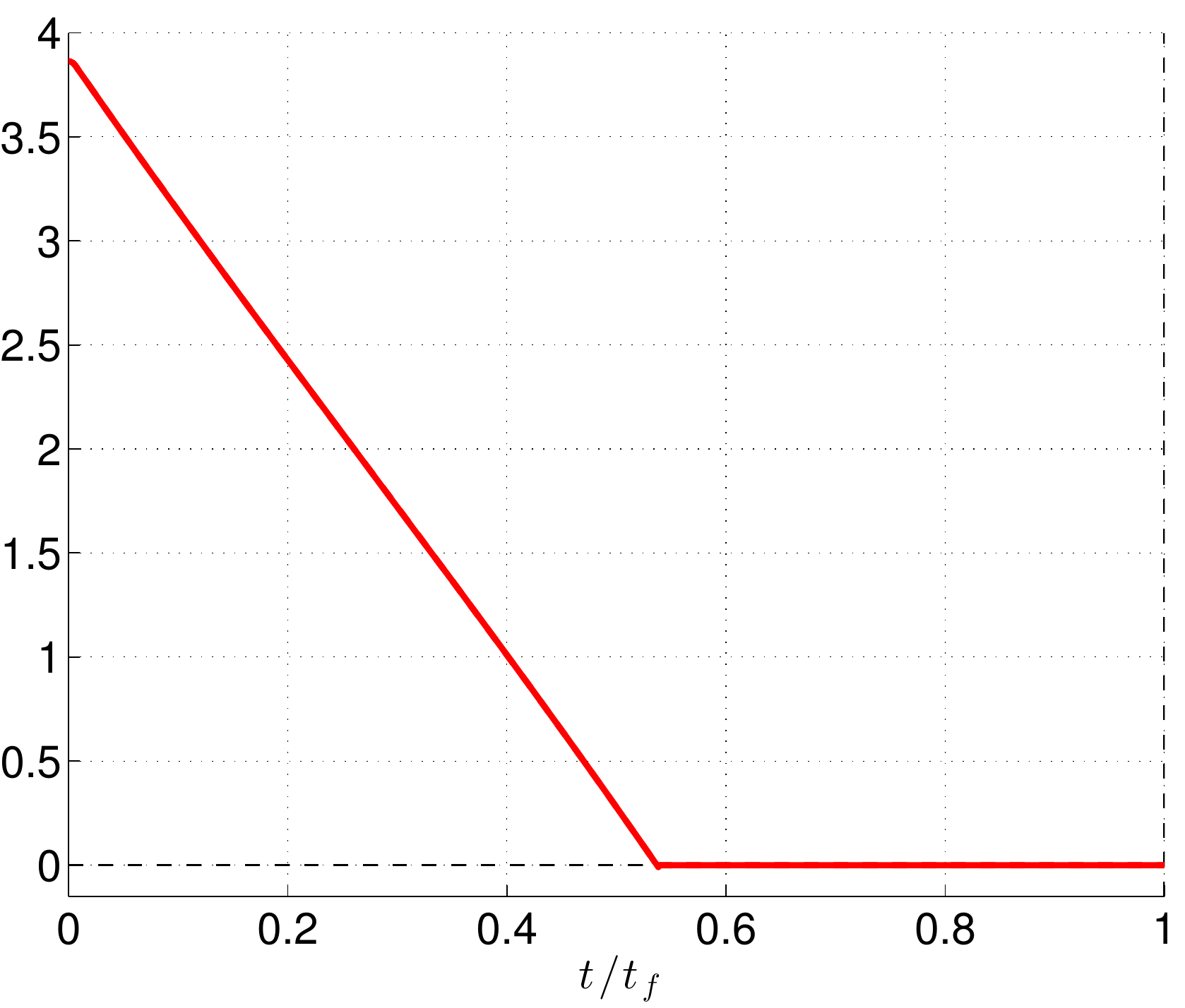}
                    \pxcoordinate{0.555*\x}{0.876*\y}{centre};
                    \draw [thin, dashed] (centre) circle (0.5);
                    \draw [-triangle 45,thin] ([xshift=0.85*\x,yshift=0.3*\y] centre) -- ([xshift=0.03*\x,yshift=0.025*\x] centre);
                \end{tikzgraphics}
                \hspace{-3em}
                \def\x{483}
                \def\y{411}
                \begin{tikzgraphics}{0.26\textwidth}{\x}{\y}{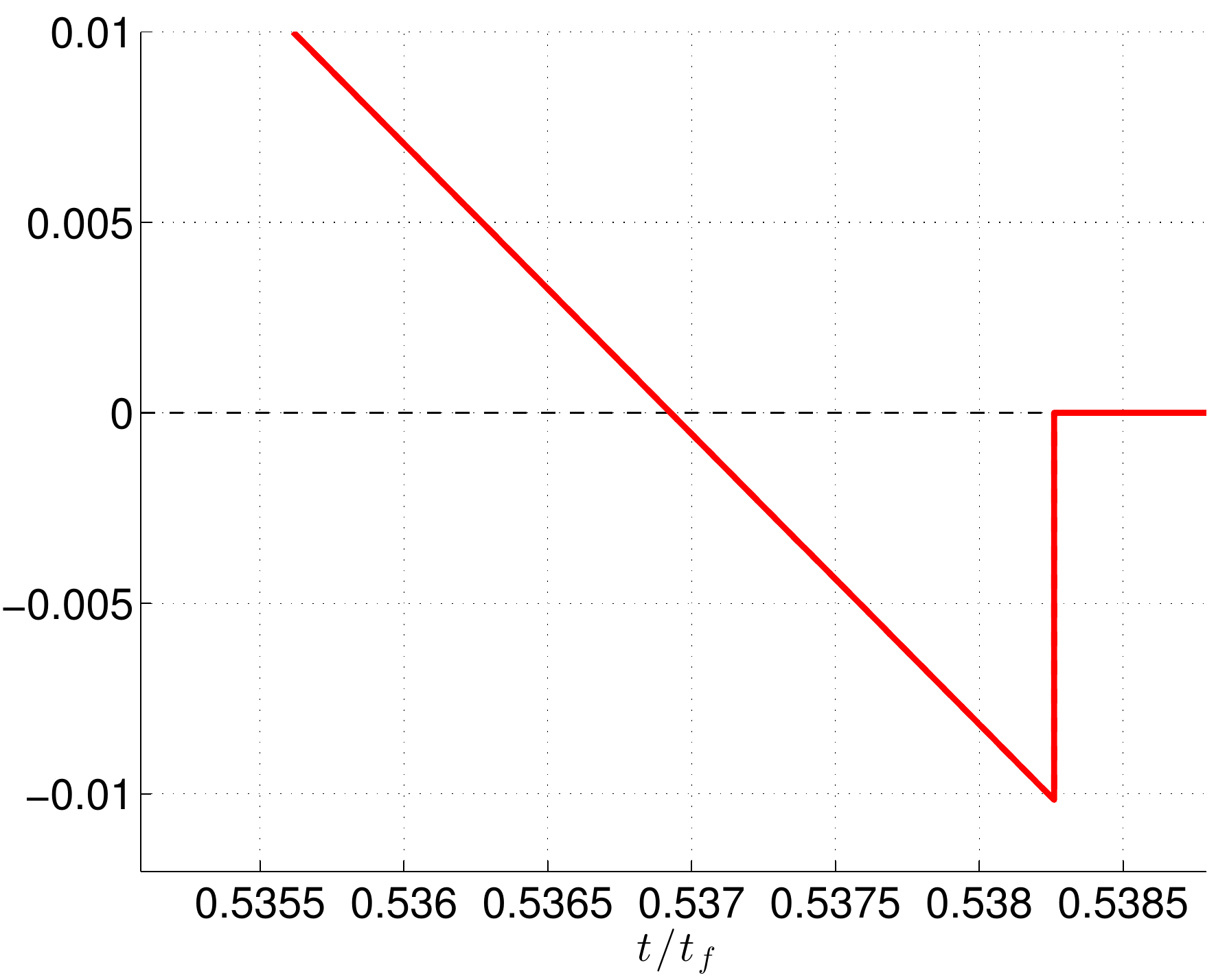}
                    \pxcoordinate{0.535*\x}{0.435*\y}{centre};
                    \draw [thin, dashed] (centre) circle (8);
                \end{tikzgraphics}
                \caption{\textbf{Trajectory $\gamma_+ \gamma_{c_1} \gamma_+ \gamma_- \gamma_{c_3}$, with $t_2 = t_3$.}
                 The BC-extremal $\gamma_+ \gamma_{c_1} \gamma_+ \gamma_- \gamma_{c_3}$ with $t_2 = t_3$ is the same as the
                BC-extremal $\gamma_+ \gamma_{c_1} \gamma_- \gamma_{c_3}$ from $\Hom_5$ at $\vmax = \vmax^{+}$.
                (Left) The component $p_3(\cdot)$ with a zoom (Right) around $\gamma_-$.
                The component $p_3(\cdot)$ is zero only along $\gamma_{c_3}$.}
                \label{fig:limiteCaseHom4Hom5_2}
            \end{figure}

        \subsubsection{Norm of homotopic functions along the paths of zeros for $\vmax \in \intervalleff{10}{110}$}
        \label{sec:normS}

            Finding zeros of shooting functions (or homotopic functions) guarantees the complete characterization of the associated BC-extremals. 
            We represent on Figure~\ref{fig:normS} the norm along the paths of zeros of the following homotopic functions:
            $\Hom_{5}$, $\Hom_{4}$, $\Hom_{3}$ and $\Hom_{2}^{(b)}$,
            respectively for $\vmax$ in $\intervalleff{10}{\vmax^+}$, $\intervalleff{\vmax^+}{\vmax^{\gamma_{c_3}}}$,
            $\intervalleff{\vmax^{\gamma_{c_3}}}{\vmax^{c_3}}$ and $\intervalleff{\vmax^{c_3}}{110}$.
            At the end of each homotopy we perform a correction using the shooting methods which explains the discontinuities at $\vmax^+$, $\vmax^{\gamma_{c_3}}$
            and $\vmax^{c_3}$.
            We can see that the corrections give very accurate solutions while along the homotopies, the norm increases quickly before reaching asymptotic values.
            We can notice two different behaviours. Along $\Hom_{5}^{-1}(0)$, the asymptotic value is around $2e^{-9}$ which is consistent with the tolerances
            of $1e^{-10}$ given to the Runge-Kutta method (\dopri). Along the other homotopies, the asymptotic value is greater 
            but around $4e^{-5}$ which is accurate enough.

            \def\sizeFig{0.58}
            \begin{figure}[ht!]
                \centering
                \def\x{560}
                \def\y{420}
                \begin{tikzgraphics}{\sizeFig\textwidth}{\x}{\y}{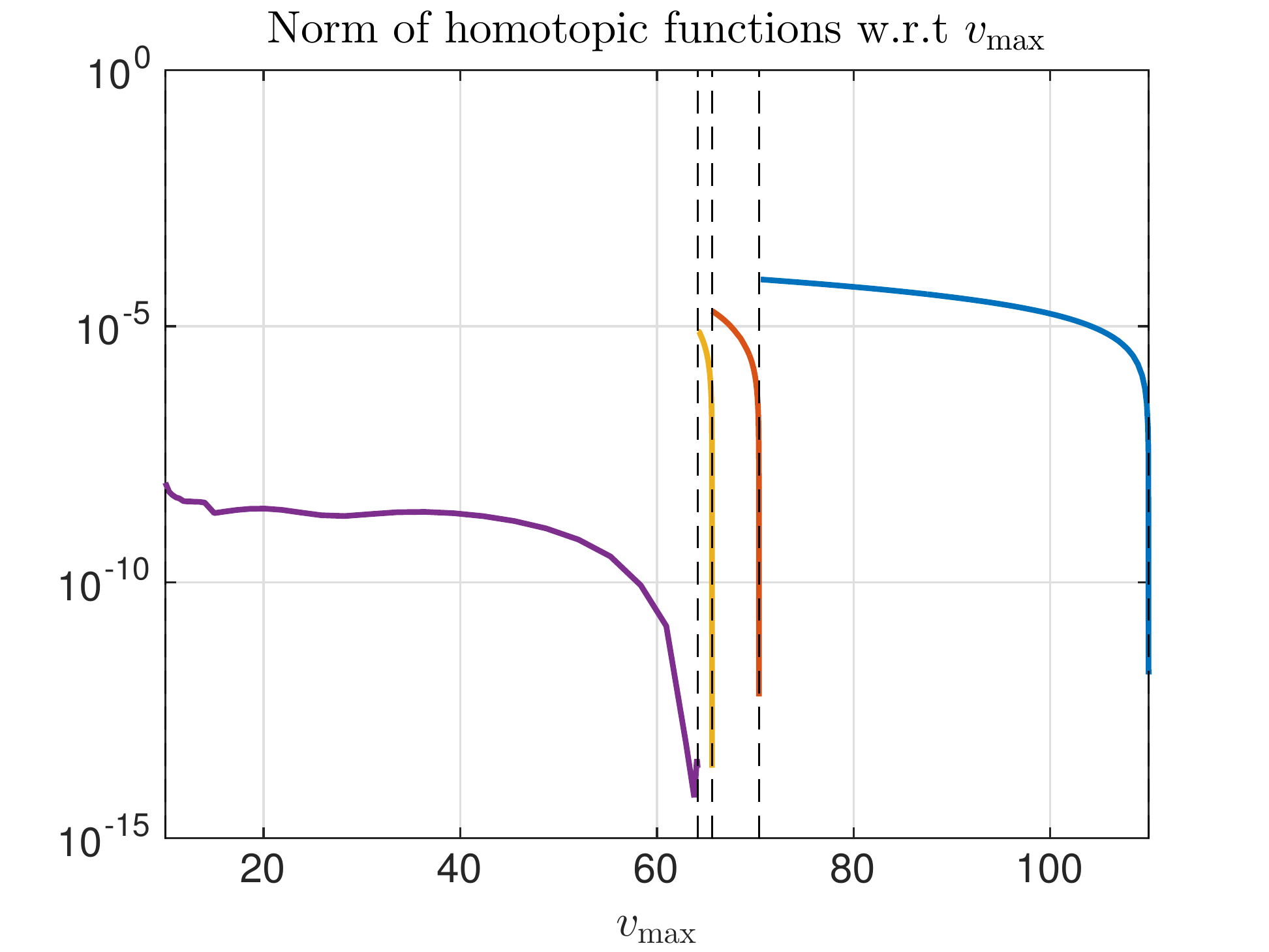}

                    \pxcoordinate{0.6*\x}{0.877*\y}{A}                                                                     ;
                    \def\xs{0.08*\x}
                    \def\ys{0.08*\y}
                    \draw  ([xshift=\xs,yshift=\ys] A) node{{\small $\vmax^{c_3}$}}                                   ;
                    \draw [->, thick, gray] ([xshift=\xs-0.4*\xs,yshift=\ys-0.4*\ys] A) -- (A)                              ;

                    \pxcoordinate{0.564*\x}{0.884*\y}{A}                                                                     ;
                    \def\xs{0.08*\x}
                    \def\ys{-0.08*\y}
                    \draw  ([xshift=\xs,yshift=\ys] A) node{{\small $\vmax^{\gamma_{c_3}}$}}                          ;
                    \draw [->, thick, gray] ([xshift=\xs-0.4*\xs,yshift=\ys-0.4*\ys] A) -- (A)                              ;

                    \pxcoordinate{0.547*\x}{0.877*\y}{A}                                                                     ;
                    \def\xs{-0.08*\x}
                    \def\ys{ 0.08*\y}
                    \draw  ([xshift=\xs,yshift=\ys] A) node{{\small $\vmax^+$}}                          ;
                    \draw [->, thick, gray] ([xshift=\xs-0.4*\xs,yshift=\ys-0.4*\ys] A) -- (A)                              ;

                \end{tikzgraphics}
                \caption{Norm along the paths of zeros of the following homotopic functions: $\Hom_{5}$, $\Hom_{4}$, $\Hom_{3}$ and $\Hom_{2}^{(b)}$,
            respectively for $\vmax$ in $\intervalleff{10}{\vmax^+}$, $\intervalleff{\vmax^+}{\vmax^{\gamma_{c_3}}}$,
            $\intervalleff{\vmax^{\gamma_{c_3}}}{\vmax^{c_3}}$ and $\intervalleff{\vmax^{c_3}}{110}$.}
                \label{fig:normS}
            \end{figure}
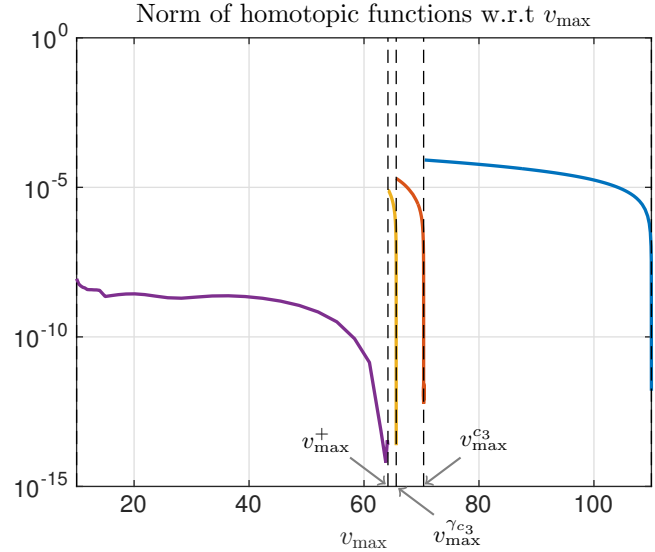

    \subsection{Synthesis with respect to $\imax$ and $\vmax$}
    \label{sec:Numericalsynthesis}

    We present in Figure~\ref{fig:synthesis} the best sub-optimal synthesis we obtain using all the techniques and results established throughout this paper.
    See section \ref{sec:method} for the procedure to construct this sub-optimal synthesis.
    We also need the proposition \ref{prop:gc1} to get the results presented in Figure~\ref{fig:synthesisZooms}.
    We make this synthesis for $\imax \ge 10$, $\vmax \ge 10$ and $\af = 100$.

        \begin{prpstn}
            Let consider a boundary arc $\gamma_{c_1}$ defined on $\intervalleff{t_1}{t_2}$ with $(z(\cdot), u_{c_1}(\cdot), \eta_{c_1}(\cdot))$
            its associated extremal. Let assume there exists
            $\tau \in \intervalleff{t_1}{t_2}$ such that $H_1(z(\tau)) = H_{01}(z(\tau)) = 0$ and that \ref{Hyp3} holds along $\gamma_{_1}$.
            Then either $\tau = t_1$ or $\tau = t_2$.
            \label{prop:gc1}
        \end{prpstn}

        \begin{proof}
                Let assume $\tau \in \intervalleoo{t_1}{t_2}$. Along $\gamma_{c_1}$, $\eta_{c_1} = H_{01} / (F_1 \cdot c_1) \le 0$ with
                $F_1 \cdot c_1 \equiv \cste_7 > 0$, so $H_{01} \le 0$.
                At time $\tau$, $H_{01}(z(\tau)) = 0$ and
                $\dot{H}_{01}(z(\tau)) = H_{001}(z(\tau)) - \eta_{c_1}(z(\tau)) (F_{01}\cdot c_1)(x(\tau)) = H_{001}(z(\tau)) \ne 0$, since
                $ \dim \vect(F_1(x(\tau)),F_{01}(x(\tau)),F_{001}(x(\tau))) = 3$ and $p(\tau) \ne 0$. As a consequence, the sign of $H_{01}$
                changes when crossing $t = \tau$, which is not possible.
        \end{proof}

       \begin{rmrk}
                To build this synthesis we do not deal with the possible local minima since the aim of this paper is not to prove the global optimality of the
                synthesis but to describe how we can build it starting from one point, \ie one BC-extremal for a given value of $(\imax,\vmax)$.
                Here we start from the optimal trajectory $\gamma_+$.
                To construct this synthesis we catch the changes along the homotopies and determine the new strategy using the theoretical results.
                However, if we consider a general optimal control problem depending on a parameter $\lambda$, then to get a better synthesis, for a given $\overline{\lambda}$
                we should compare the cost associated to each component of $h^{-1}(\{0\}) \cap \{ \lambda=\overline{\lambda} \}$, for each
                homotopic function $h$. This approach is crucial when the optimal control problem has for example many local solutions, see \cite{ACTA}.
       \end{rmrk}

    \def\sizeFig{0.65}
    \def\x{522}
    \def\y{407}
    \begin{figure}[ht!]
        \begin{tikzgraphics}{\sizeFig\textwidth}{\x}{\y}{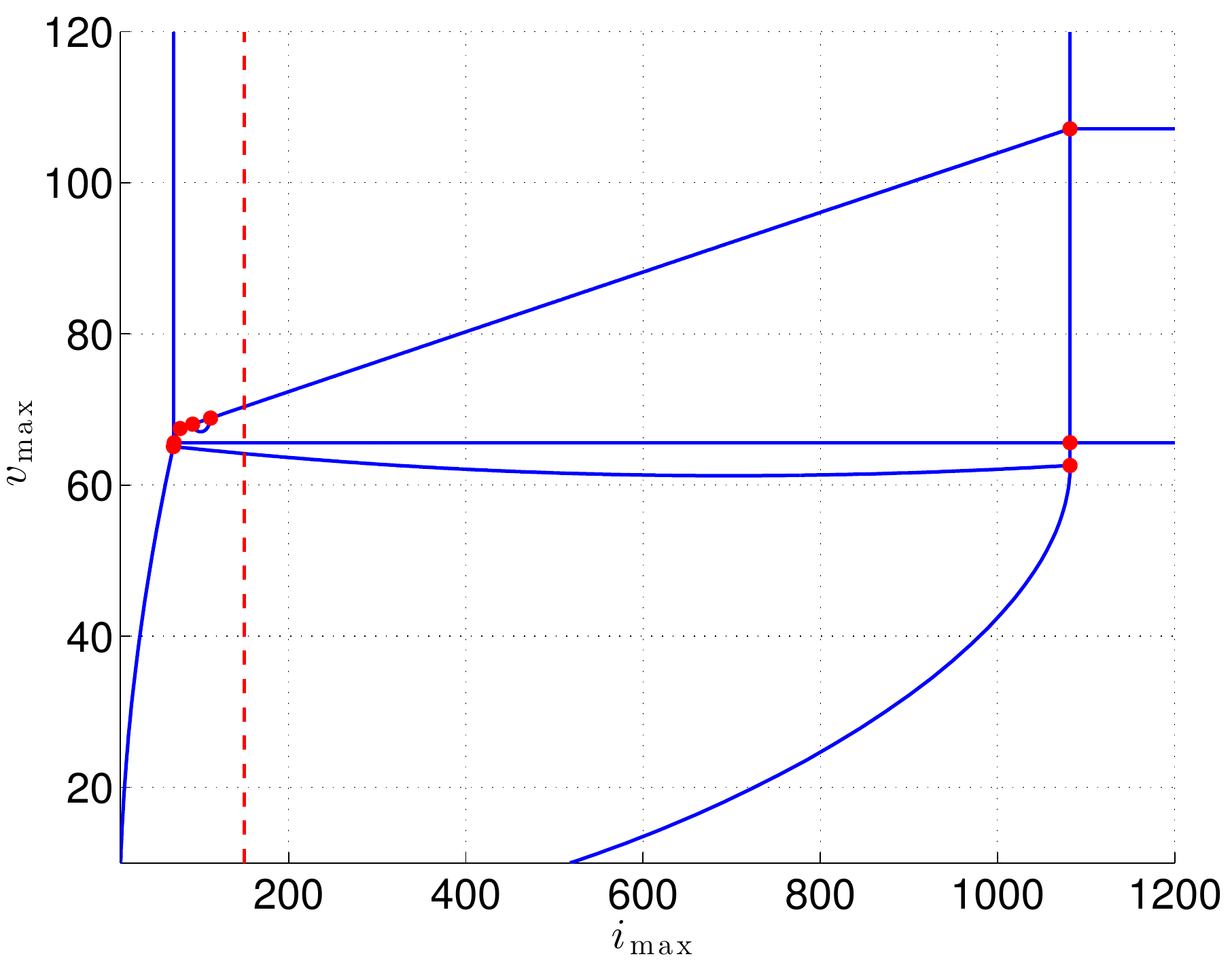}

            \pxcoordinate{0.1985*\x}{0.898*\y}{A}                                                                   ;
            \def\xs{-0.04*\x}
            \def\ys{-0.04*\y}
            \draw [red] ([xshift=\xs,yshift=\ys] A) node{{ $150$}}                                                  ;
            \draw [->, thick, red] ([xshift=\xs-0.4*\xs,yshift=\ys-0.4*\ys] A) -- (A)                               ;

            \pxcoordinate{0.869*\x}{0.08*\y}{A}                                                                     ;
            \def\xs{-0.04*\x}
            \def\ys{0.02*\y}
            \draw [gray] ([xshift=\xs,yshift=\ys] A) node{{\small $\gamma_+^{c_1}$}}                                ;
            \draw [->, thick, gray] ([xshift=\xs-0.4*\xs,yshift=\ys-0.4*\ys] A) -- (A)                              ;

            \pxcoordinate{0.91*\x}{0.134*\y}{A}                                                                     ;
            \def\xs{0.02*\x}
            \def\ys{-0.06*\y}
            \draw [gray] ([xshift=\xs,yshift=\ys] A) node{{\small $\gamma_+^{H_1,c_3}$}}                            ;
            \draw [->, thick, gray] ([xshift=\xs-0.4*\xs,yshift=\ys-0.4*\ys] A) -- (A)                              ;

            \pxcoordinate{0.91*\x}{0.46*\y}{A}                                                                      ;
            \def\xs{0.02*\x}
            \def\ys{0.04*\y}
            \draw [gray] ([xshift=\xs,yshift=\ys] A) node{{\small $u_{c_3} \equiv +1$}}                             ;
            \draw [->, thick, gray] ([xshift=\xs-0.4*\xs,yshift=\ys-0.4*\ys] A) -- (A)                              ;

            \pxcoordinate{0.7*\x}{0.46*\y}{A}                                                                       ;
            \def\xs{0.02*\x}
            \def\ys{0.04*\y}
            \draw [gray] ([xshift=\xs,yshift=\ys] A) node{{\small $u_{c_3} \equiv +1$}}                             ;
            \draw [->, thick, gray] ([xshift=\xs-0.4*\xs,yshift=\ys-0.4*\ys] A) -- (A)                              ;

            \pxcoordinate{0.7*\x}{0.757*\y}{A}                                                                      ;
            \def\xs{-0.04*\x}
            \def\ys{0.06*\y}
            \draw [gray] ([xshift=\xs,yshift=\ys] A) node{{\small $\gamma_+ \delta_{c_1} \gamma_- \gamma_{c_3}$}}   ;
            \draw [->, thick, gray] ([xshift=\xs-0.4*\xs,yshift=\ys-0.4*\ys] A) -- (A)                              ;

            \pxcoordinate{0.869*\x}{0.25*\y}{A}                                                                     ;
            \def\xs{-0.06*\x}
            \def\ys{0.03*\y}
            \draw [gray] ([xshift=\xs,yshift=\ys] A) node{{\small $\gamma_+^{c_1} \gamma_- \gamma_+^{c_3}$}}        ;
            \draw [->, thick, gray] ([xshift=\xs-0.4*\xs,yshift=\ys-0.6*\ys] A) -- (A)                              ;

            \pxcoordinate{0.4*\x}{0.335*\y}{A}                                                                      ;
            \def\xs{-0.02*\x}
            \def\ys{0.06*\y}
            \draw [gray] ([xshift=\xs,yshift=\ys] A) node{{\small $\gamma_+ \gamma_{c_1} \gamma_+^{H_1,c_3}$}}      ;
            \draw [->, thick, gray] ([xshift=\xs-0.4*\xs,yshift=\ys-0.4*\ys] A) -- (A)                              ;

            \pxcoordinate{0.869*\x}{0.473*\y}{A}                                                                    ;
            \def\xs{0.08*\x}
            \def\ys{-0.02*\y}
            \draw [gray] ([xshift=\xs,yshift=\ys] A) node{{\small $\gamma_+^{c_1} \gamma_- \gamma_{c_3}$}}          ;
            \draw [->, thick, gray] ([xshift=\xs-0.6*\xs,yshift=\ys-0.4*\ys] A) -- (A)                              ;

            \pxcoordinate{0.7*\x}{0.495*\y}{A}                                                                       ;
            \def\xs{0.02*\x}
            \def\ys{-0.04*\y}
            \draw [gray] ([xshift=\xs,yshift=\ys] A) node{{\small $u_{c_1}(t_2) = +1$}}                             ;
            \draw [->, thick, gray] ([xshift=\xs-0.4*\xs,yshift=\ys-0.4*\ys] A) -- (A)                              ;

            \pxcoordinate{0.14*\x}{0.05*\y}{A}                                                                       ;
            \def\xs{0.08*\x}
            \def\ys{0.04*\y}
            \draw [gray] ([xshift=\xs,yshift=\ys] A) node{{\small $u_{c_1}(t_f) = +1$}}     ;
            \draw [->, thick, gray] ([xshift=\xs-0.4*\xs,yshift=\ys-0.6*\ys] A) -- (A)                              ;

            \pxcoordinate{0.11*\x}{0.70*\y}{A}                                                                       ;
            \def\xs{0.08*\x}
            \def\ys{-0.04*\y}
            \draw [gray] ([xshift=\xs,yshift=\ys] A) node{{\small $x_{3}(t_f) = +1$}}     ;
            \draw [->, thick, gray] ([xshift=\xs-0.4*\xs,yshift=\ys-0.6*\ys] A) -- (A)                              ;


            \pxnode{0.91*\x}{0.08*\y}{A}   {{\tiny $\circ$}}                                                        ;
            \def\xs{0.02*\x}
            \def\ys{0.02*\y}
            \draw [black] ([xshift=\xs,yshift=\ys] A) node{\large $\gamma_+$}                                              ;

            \pxnode{0.91*\x}{0.3*\y}{A}   {{\tiny $\circ$}}                                                         ;
            \def\xs{0.065*\x}
            \def\ys{0.02*\y}
            \draw [black] ([xshift=\xs,yshift=\ys] A) node{\large $\gamma_+ \gamma_- \gamma_+^{c_3}$}                      ;

            \pxnode{0.5*\x}{0.15*\y}{A}   {{\tiny $\circ$}}                                                         ;
            \def\xs{0.06*\x}
            \def\ys{0.02*\y}
            \draw [black] ([xshift=\xs,yshift=\ys] A) node{\large $\gamma_+ \gamma_{c_1} \gamma_+$}                        ;

            \pxnode{0.12*\x}{0.15*\y}{A}   {{\tiny $\circ$}}                                                        ;
            \def\xs{-0.02*\x}
            \def\ys{0.03*\y}
            \draw [black] ([xshift=\xs,yshift=\ys] A) node{\large $\gamma_+ \gamma_{c_1}$}                                 ;

            \pxnode{0.6*\x}{0.35*\y}{A}   {{\tiny $\circ$}}                                                         ;
            \def\xs{0.105*\x}
            \def\ys{0.025*\y}
            \draw [black] ([xshift=\xs,yshift=\ys] A) node{\large $\gamma_+ \gamma_{c_1} \gamma_+ \gamma_- \gamma_+^{c_3}$};

            \pxnode{0.5*\x}{0.478*\y}{A}   {{\tiny $\circ$}}                                                        ;
            \def\xs{0.11*\x}
            \def\ys{0.00*\y}
            \draw [black] ([xshift=\xs,yshift=\ys] A) node{\normalsize $\gamma_+ \gamma_{c_1} \gamma_+ \gamma_- \gamma_{c_3}$}  ;

            \pxnode{0.4*\x}{0.65*\y}{A}   {{\tiny $\circ$}}                                                         ;
            \def\xs{0.09*\x}
            \def\ys{0.02*\y}
            \draw [black] ([xshift=\xs,yshift=\ys] A) node{\large $\gamma_+ \gamma_{c_1} \gamma_- \gamma_{c_3}$}           ;

            \pxnode{0.82*\x}{0.75*\y}{A}   {{\tiny $\circ$}}                                                         ;
            \def\xs{0.065*\x}
            \def\ys{0.02*\y}
            \draw [black] ([xshift=\xs,yshift=\ys] A) node{\large $\gamma_+ \gamma_- \gamma_{c_3}$}           ;

            \pxcoordinate{0.17*\x}{0.45*\y}{centre};
            \draw [thin] (centre) circle (1.0);
            \draw [-triangle 45,thin] ([xshift=0.10*\x,yshift=-0.08*\y] centre) -- ([xshift=0.04*\x,yshift=-0.04*\x] centre);
            \draw [black] ([xshift=0.12*\x,yshift=-0.10*\y] centre) node{See Fig.~\ref{fig:synthesisZooms}}           ;

        \end{tikzgraphics}

        \caption{\textbf{Synthesis with respect to $\imax$ and $\vmax$}.
        The structures are displayed on the graph.
        The sequence labeled $\gamma_+ \delta_{c_1} \gamma_-$ means that there is a contact point with the constraint set
        $C_1$ at the switching time between $\gamma_+$ and $\gamma_-$.
        The blue lines represent structures of degree 2 while at the red points we have structures of degree 3 or 4.}
        \label{fig:synthesis}
    \end{figure}

    \def\sizeFig{0.5}
    \def\x{492}
    \def\y{407}
    \begin{figure}[ht!]
        \begin{tikzgraphics}{\sizeFig\textwidth}{\x}{\y}{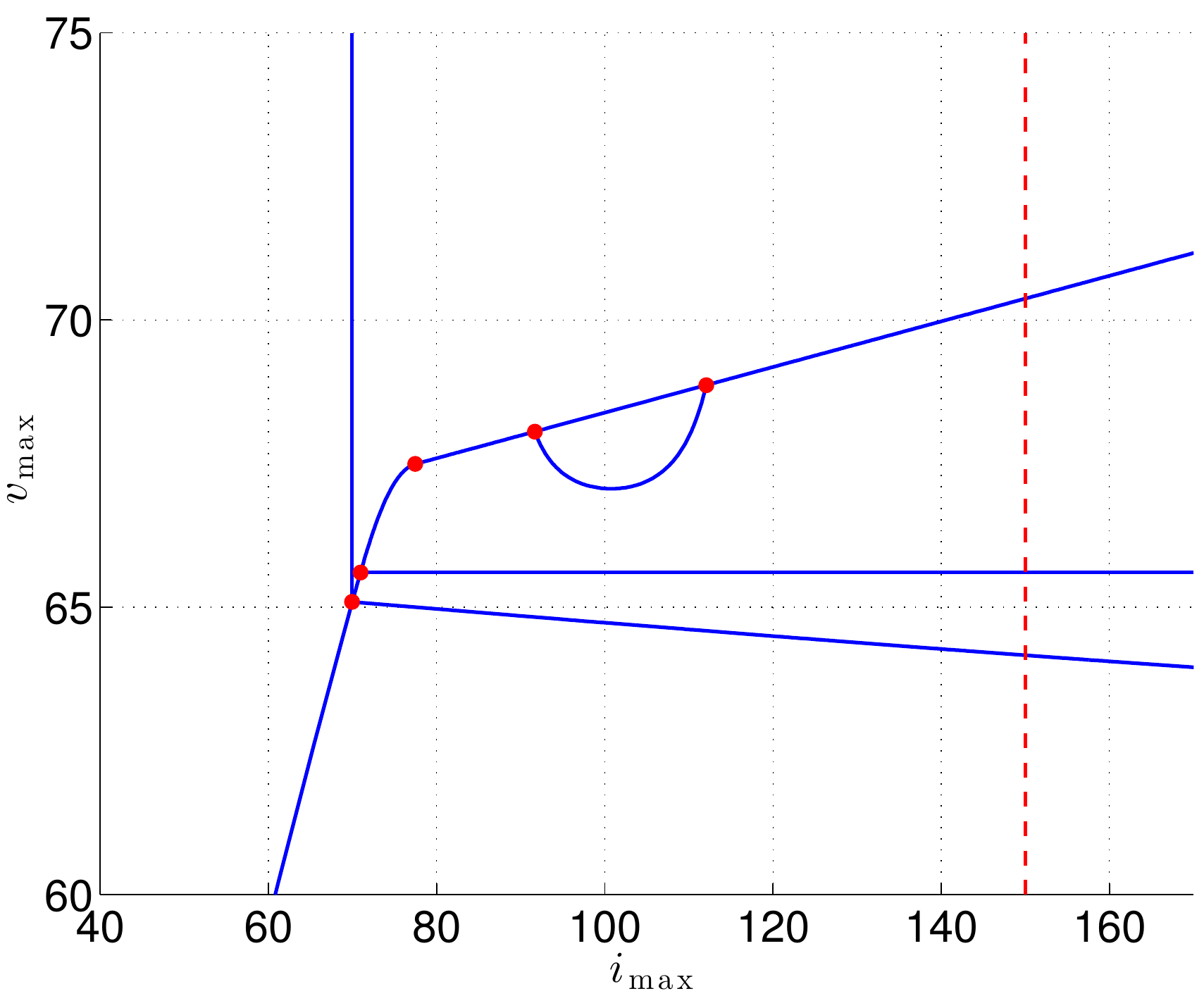}

            \pxcoordinate{0.29*\x}{0.125*\y}{A}                                                                       ;
            \def\xs{0.14*\x}
            \def\ys{0.06*\y}
            \draw [gray] ([xshift=\xs,yshift=\ys] A) node{{\small $u_{c_1}(t_f) = +1$}}     ;
            \draw [->, thick, gray] ([xshift=\xs-0.6*\xs,yshift=\ys-0.6*\ys] A) -- (A)                              ;

            \pxcoordinate{0.25*\x}{0.80*\y}{A}                                                                       ;
            \def\xs{-0.12*\x}
            \def\ys{0.04*\y}
            \draw [gray] ([xshift=\xs,yshift=\ys] A) node{{\small $x_{3}(t_f) = +1$}}     ;
            \draw [->, thick, gray] ([xshift=\xs-0.4*\xs,yshift=\ys-0.6*\ys] A) -- (A)                              ;

            \pxcoordinate{0.9*\x}{0.575*\y}{A}                                                                       ;
            \def\xs{0.02*\x}
            \def\ys{0.06*\y}
            \draw [gray] ([xshift=\xs,yshift=\ys] A) node{{\small $u_{c_3} \equiv +1$}}                             ;
            \draw [->, thick, gray] ([xshift=\xs-0.4*\xs,yshift=\ys-0.4*\ys] A) -- (A)                              ;

            \pxcoordinate{0.9*\x}{0.663*\y}{A}                                                                       ;
            \def\xs{0.02*\x}
            \def\ys{-0.08*\y}
            \draw [gray] ([xshift=\xs,yshift=\ys] A) node{{\small $u_{c_1}(t_2) = +1$}}                             ;
            \draw [->, thick, gray] ([xshift=\xs-0.4*\xs,yshift=\ys-0.5*\ys] A) -- (A)                              ;

            \pxcoordinate{0.9*\x}{0.285*\y}{A}                                                                      ;
            \def\xs{0.04*\x}
            \def\ys{0.08*\y}
            \draw [gray] ([xshift=\xs,yshift=\ys] A) node{{\small $\gamma_+ \gamma_{c_1} \gamma_+^{H_1,c_3}$}}      ;
            \draw [->, thick, gray] ([xshift=\xs-0.4*\xs,yshift=\ys-0.4*\ys] A) -- (A)                              ;

            \pxcoordinate{0.587*\x}{0.388*\y}{A}                                                                      ;
            \def\xs{0.06*\x}
            \def\ys{0.10*\y}
            \draw [gray] ([xshift=\xs,yshift=\ys] A) node{{\small $\gamma_+ \gamma_{c_1} \delta_{H_1} \gamma_+^{c_3}$}}   ;
            \draw [->, thick, gray] ([xshift=\xs-0.4*\xs,yshift=\ys-0.4*\ys] A) -- (A)                              ;

            \pxcoordinate{0.4*\x}{0.449*\y}{A}                                                                      ;
            \def\xs{-0.14*\x}
            \def\ys{0.02*\y}
            \draw [gray] ([xshift=\xs,yshift=\ys] A) node{{\small $\gamma_+ \gamma_{c_1} \gamma_+^{H_1,c_3}$}}      ;
            \draw [->, thick, gray] ([xshift=\xs-0.7*\xs,yshift=\ys-0.2*\ys] A) -- (A)                              ;

            \pxcoordinate{0.321*\x}{0.5*\y}{A}                                                                      ;
            \def\xs{-0.18*\x}
            \def\ys{-0.03*\y}
            \draw [gray] ([xshift=\xs,yshift=\ys] A) node{{\small $\gamma_+ \gamma_{c_1} \gamma_+ \delta_{H_1,c_3}$}}      ;
            \draw [->, thick, gray] ([xshift=\xs-0.6*\xs,yshift=\ys-0.2*\ys] A) -- (A)                              ;

            \pxcoordinate{0.48*\x}{0.483*\y}{A}                                                                       ;
            \def\xs{-0.02*\x}
            \def\ys{-0.05*\y}
            \draw [gray] ([xshift=\xs,yshift=\ys] A) node{{\small $u_{c_1}(t_2) = +1$}}                             ;
            \draw [->, thick, gray] ([xshift=\xs-0.4*\xs,yshift=\ys-0.5*\ys] A) -- (A)                              ;

            \pxcoordinate{0.52*\x}{0.41*\y}{A}                                                                      ;
            \def\xs{-0.10*\x}
            \def\ys{0.06*\y}
            \draw [gray] ([xshift=\xs,yshift=\ys] A) node{{\small $\gamma_+ \gamma_{c_1} \gamma_- \gamma_+^{c_1,c_3}$}}      ;
            \draw [->, thick, gray] ([xshift=\xs-0.8*\xs,yshift=\ys-0.4*\ys] A) -- (A)                              ;

            \pxcoordinate{0.52*\x}{0.45*\y}{A}                                                                      ;
            \def\xs{0.18*\x}
            \def\ys{-0.02*\y}
            \draw [gray] ([xshift=\xs,yshift=\ys] A) node{{\small $\gamma_+ \gamma_{c_1} \gamma_- \gamma_+^{c_3}$}}      ;
            \draw [->, thick, gray] ([xshift=\xs-0.6*\xs,yshift=\ys-0.4*\ys] A) -- (A)                              ;


            \pxnode{0.2*\x}{0.15*\y}{A}   {{\tiny $\circ$}}                                                        ;
            \def\xs{-0.02*\x}
            \def\ys{0.04*\y}
            \draw [black] ([xshift=\xs,yshift=\ys] A) node{\large $\gamma_+ \gamma_{c_1}$}                                 ;

            \pxnode{0.6*\x}{0.15*\y}{A}   {{\tiny $\circ$}}                                                         ;
            \def\xs{0.10*\x}
            \def\ys{0.04*\y}
            \draw [black] ([xshift=\xs,yshift=\ys] A) node{\large $\gamma_+ \gamma_{c_1} \gamma_+$}                        ;

            \pxnode{0.82*\x}{0.4*\y}{A}   {{\tiny $\circ$}}                                                         ;
            \def\xs{0.16*\x}
            \def\ys{0.04*\y}
            \draw [black] ([xshift=\xs,yshift=\ys] A) node{\large $\gamma_+ \gamma_{c_1} \gamma_+ \gamma_- \gamma_+^{c_3}$};

            \pxnode{0.6*\x}{0.61*\y}{A}   {{\tiny $\circ$}}                                                        ;
            \def\xs{0.16*\x}
            \def\ys{0.00*\y}
            \draw [black] ([xshift=\xs,yshift=\ys] A) node{\large $\gamma_+ \gamma_{c_1} \gamma_+ \gamma_- \gamma_{c_3}$}  ;

            \pxnode{0.45*\x}{0.75*\y}{A}   {{\tiny $\circ$}}                                                         ;
            \def\xs{0.125*\x}
            \def\ys{-0.04*\y}
            \draw [black] ([xshift=\xs,yshift=\ys] A) node{\large $\gamma_+ \gamma_{c_1} \gamma_- \gamma_{c_3}$}           ;

        \end{tikzgraphics}
        \caption{\textbf{Zoom of the synthesis}.
        The structures are displayed on the graph. The sequence labeled $\gamma_{c_1} \delta_{H_1}$ means that there is a contact point of order 2
        with the switching surface at the exit-time of the boundary arc.}
        \label{fig:synthesisZooms}
    \end{figure}
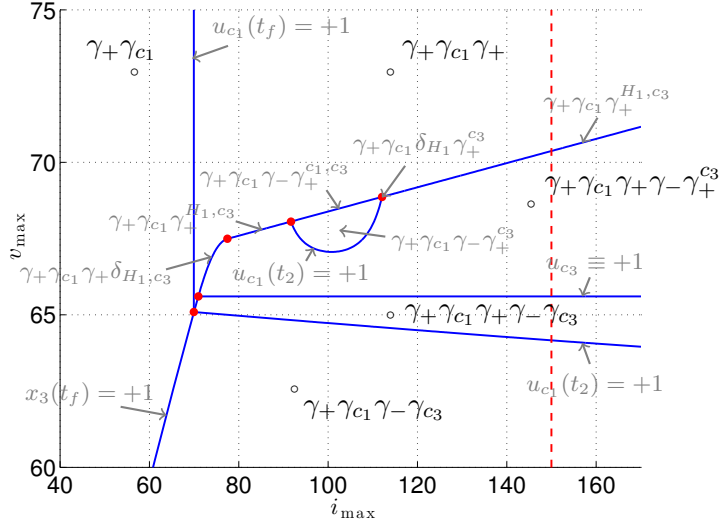

\section{Conclusion}
\label{sec:conclusion}

    In this paper, we have presented the minimum time control problem of an electric vehicle which have been modeled as a Mayer problem in optimal control,
    with affine dynamics with respect to the control (scalar) and with two state constraints, one of order 1 and the other of order 2, see
    sections \ref{sec:Model} and \ref{sec:ApplicationControleContraint}.
    We have tackled first the state unconstrained problem from the application of the Pontryagin Maximum Principle, see section
    \ref{sec:PMP}, which gives necessary conditions of
    optimality. Then we used in section \ref{sec:LieBracketConf} the Lie bracket configuration to show that the only existing extremals are
    bang-bang with finitely many switchings (see section \ref{sec:SingularRegularExtremals}) and we proved in \ref{sec:bangOptimal}
    that the optimal trajectory is of the form $\gamma_+$.
    %
    %
    Nevertheless, the classification of bang-bang extremals, presented in section \ref{sec:Classification}, was an efficient tool to describe the behavior
    of some possible extremals (here for the state constrained case) and to get better insight into such extremals, see section \ref{sec:homotopieC3}.
    The analysis of problem \eqref{ref:OCP_tf_min}, \ie with the state constraints, consisted first in computing the boundary controls and the multipliers
    associated to the state constraints, for both the first and second-order constraints. On the other hand, we gave new junction conditions,
    see section \ref{sec:BoundaryCase}, which have been used to define multiple shooting functions in sections \ref{sec:Shooting2} and \ref{sec:Shooting4}.
    Moreover, we gave a local time minimal synthesis, see section \ref{sec:timeMinimalSynthesis}, that we encountered in the numerical simulations in sections
    \ref{sec:Hom4} and \ref{sec:Hom5}. Finally, in section \ref{sec:Numericalsynthesis}, we obtained a synthesis of solutions which
    satisfy the necessary conditions of optimality (\ie BC-extremals) with respect to the parameters $\imax$ and $\vmax$, with $x_1(0)$, $x_3(0)$, $\af$ fixed,
    and for one specific car.

\begin{figure}[ht!]
        \includegraphics[width=\tfig\textwidth]{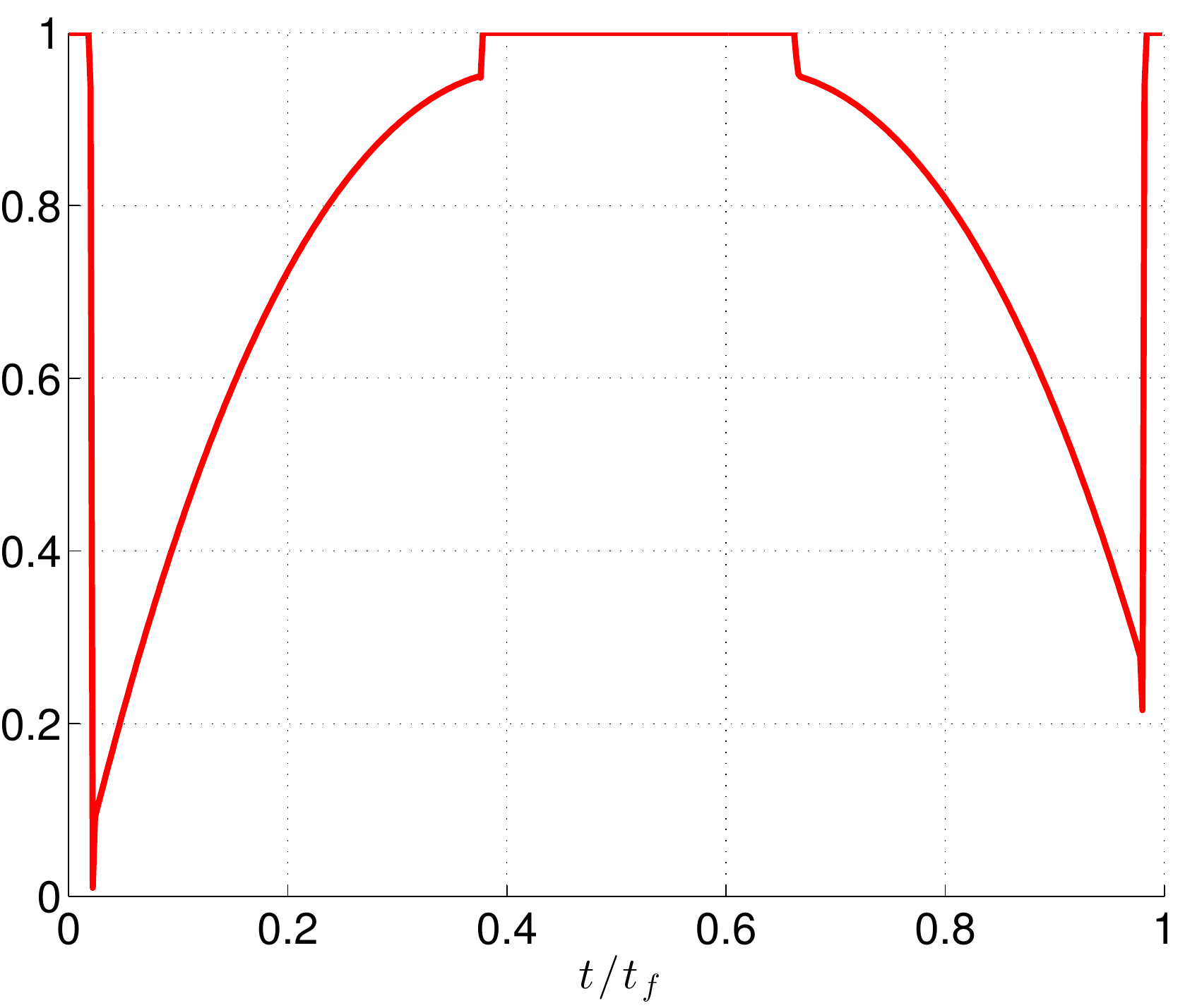}%
        \hspace{3em}
        \includegraphics[width=\tfig\textwidth]{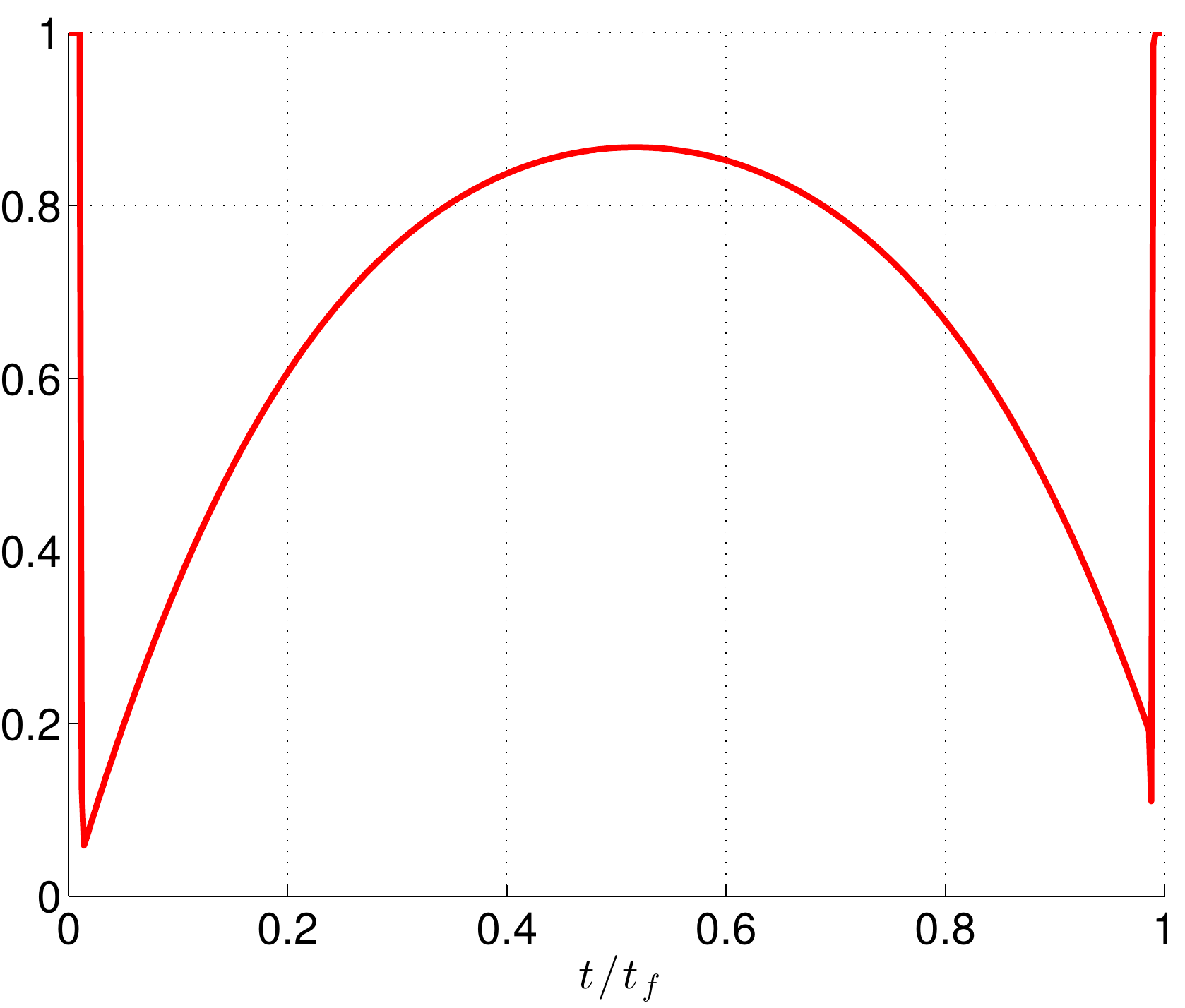}%
        \caption{\textbf{Problem \eqref{ref:OCP_E_min}}.
        Optimal control for problem \eqref{ref:OCP_E_min} for $t_f = 1.25\, T_\mathrm{min}$ (left) and $t_f = 1.50\, T_\mathrm{min}$ (right), with
        $T_\mathrm{min} \approx 5.6156$, $\af = 100$, $\imax = 1100$ and $\vmax = 110$. The structure is $\gamma_+ \gamma_s \gamma_+ \gamma_s \gamma_+$
        for $t_f = 1.25\, T_\mathrm{min}$ and $\gamma_+ \gamma_s \gamma_+$ for $t_f = 1.50\, T_\mathrm{min}$. The solutions have been
        computed with the \bocop\ \cite{Bocop} software.}
        \label{fig:Energie_min_controle}
    \end{figure}

    The same techniques could be applied to other problems, such as the problem of the minimization of the electric energy \eqref{ref:OCP_E_min}. In this case,
    the final time $t_f$ is fixed and the structure of the optimal trajectories depends on $t_f$, see Figure~\ref{fig:Energie_min_controle}.
    The dynamics associated to problem \eqref{ref:OCP_E_min} is bilinear on the state and the control variables and this leads to more complex optimal
    structures with singular extremals. Even though the state unconstrained problem \eqref{ref:OCP_E_min} is more complex than the time minimal case,
    the same tools such as Lie bracket configuration and local classification can be used.
    The main difference comes from the fact that the final time is now a parameter and one may study its influence on the optimal trajectories. This was already done in
    \cite{ACTA}, where a state unconstrained and affine scalar control problem of Mayer form was analyzed taking into account the influence of the final time.
    Moreover, the state constrained case should be clarified with the same tools as those presented in this paper.

\bibliographystyle{model1-num-names}


\begin{acknowledgement}
The author is thankful to Bernard Bonnard, Jean-Baptiste Caillau and Joseph Gergaud for their geometric and numerical influences on this work.
\end{acknowledgement}

\end{document}